\documentclass[12pt]{article}

\usepackage{amsmath,amssymb,amsthm,amsfonts}
\usepackage[utf8]{inputenc}
\usepackage[T1]{fontenc}
\usepackage[french,english]{babel}
\usepackage{xcolor}
\usepackage{soul}
\usepackage{subcaption}
\usepackage{caption}
\usepackage{floatflt}
\usepackage{dsfont}
\usepackage{multirow}
\usepackage{mathrsfs}
\usepackage{diagbox}
\usepackage{ulem}
\usepackage{verbatim}
\usepackage{moreverb}
\usepackage{listings}
\usepackage{comment}
\usepackage{enumerate}
\usepackage{cite}
\usepackage[top=2cm,bottom=2cm,right=2cm,left=2cm]{geometry}
\usepackage{hyperref}
\usepackage{authblk}
\usepackage{graphicx}

\newtheorem{Theorem}{Theorem}[section]
\newtheorem {Corollary}{Corollary}[section]
\newtheorem {Proposition}{Proposition}[section]
\newtheorem{Rem}{Remark}[section]

\DeclareMathOperator{\pen}{pen}

\title{Choosing the penalty in nonparametric regression: short and long-range dependence}

\author[1]{Emmanuel Caron Parte}
\affil[1]{Avignon Universit\'e, Laboratoire de Math\'ematiques d’Avignon EA2151, 84000 Avignon, France,  \url{emmanuel.caron-parte@univ-avignon.fr}}
\author[2]{Jérôme Dedecker}
\affil[2]{Universit\'e Paris Cité, CNRS, MAP5, F-75006 Paris, France,  \url{jerome.dedecker@u-pariscite.fr}}
\author[3]{Bertrand Michel}
\affil[3]{Ecole Centrale Nantes, Nantes Universit\'e, Laboratoire de Math\'ematiques Jean Leray UMR CNRS 6629,
France, \url{bertrand.michel@ec-nantes.fr}}

\date{}

\begin{document}

\maketitle

\begin{abstract}
    
In this work, we study the one-dimensional regression problem under random design and Gaussian errors. Our framework is very general: we make no prior assumptions about the design (which may be nonstationary and exhibit short or long-range dependence), nor do we assume that the errors are homoscedastic. We examine in detail the cases where the error process exhibits short or long-range  dependence.
We adopt a least-squares penalized strategy using piecewise polynomials to estimate the regression function, following the framework of Baron, Birgé and Massart \cite{Barronetal1999}. We derive explicit penalties, up to calibration constants, to obtain adaptive estimators for which we establish risk bounds.
Since these penalties depend on the dependence properties of the error process, which are unknown in practice, 
we  propose several adaptations of the \textit{dimension jump} calibration algorithm to make our procedures fully data-driven.
\end{abstract}


\smallskip


\section{Introduction}

In this article, we focus on the model
\begin{equation}
\label{base_model}
    Y_i=f^*(X_i) + \varepsilon_i \, ,
\end{equation}
in which one observes the couples $(X_i, Y_i)_{1 \leq i \leq n}$ taking values in ${\mathbb R}^2$.  Let ${\mathbf X}_n=(X_i)_{1 \leq i \leq n}$ and $ {\mathbf e}_n=(\varepsilon_i)_{1 \leq i \leq n}$.
We assume that the  conditional  distribution  of ${\mathbf e}_n$ given ${\mathbf X}_n$ is a normal distribution ${\mathcal N}_n(0, \Sigma({\mathbf X}_n))$. Our goal is to estimate the function  $f^*$ on a given interval $I$ using piecewise polynomials defined on a regular partition of $I$, without making any assumptions about the regularity of $f^*$, and to obtain non-asymptotic bounds for the conditional mean-square error (conditional MSE)
$$
{\mathbb E} \left [ \sum_{i=1}^n (\hat f (X_i)- f^*(X_i))^2 {\mathbf 1}_{X_i \in I} \Big |{\mathbf X}_n \right] \, .
$$

We will use penalized estimators, following the approach described in detail in \cite{Barronetal1999} and \cite{birge2001gaussian} (see also \cite{Baraud2000} for the case of fixed-design regression).
In this context, we start with the theoretical (and in practice uncomputable) penalty proposed in \cite{caron2021gaussian} by working conditionally on ${\mathbf X}_n$. Our first objective is to provide explicit penalties, up to calibration constants, in order to obtain adaptive estimators that are computable in practice.

Our first important result is as follows: when the error process ${\mathbf e}_n$ is short-memory (in the broadest possible sense, based on the spectral radius of the  matrix $\Sigma({\mathbf X}_n)$), we show that the penalty has the same form as in the case where the errors are independent and identically distributed  (i.i.d.), and this holds true without any assumptions about the design ${\mathbf X}_n$. To calibrate the unknown constant in the penalty term, we can use the dimension-jump algorithm implemented in the  R package \texttt{capushe}. We verify through simulations that the penalized estimator obtained by this method performs well. Note that the variance matrix $\Sigma({\mathbf X}_n)$ is allowed to depend on the full process ${\mathbf X}_n$, and that our results are therefore valid in a very general heteroscedastic framework. 

In the case of long-memory error processes, things are, of course, more complicated.  When the conditional autocovariance function of the errors decays as $k^{-\gamma}$ for $\gamma$ in $(0,1)$ (see \eqref{decaygamma}), we show that  the penalty can be chosen to be proportional to $m^\gamma$, where $m$ is the size of the partition of $I$ (see \eqref{pen2_ld}).  We thus recover the penalty proposed in \cite{caron2021gaussian} for  fixed-designs, which leads to optimal rates in that case. We also show, using simulations, that this penalty performs well - provided that $\gamma$ is known - for a wide variety of designs. However, we also observe that it is suboptimal, for example, when the design ${\mathbf X}_n$ is i.i.d. and independent of ${\mathbf e}_n$: in this case, the standard penalty seems to perform better in certain situations. We provide some explanations for this phenomenon in Section \ref{sec:main_longrange2}, where we propose a heuristic that allows us to conjecture that the standard penalty is preferable under a condition involving both the dependency structure of the design and that of the error process (see Proposition \ref{prop_ld_improved}).

In practice, unless the data come from a known sampling scheme, the user does not know whether the error process has short memory, what the value of the $\gamma$ parameter controlling the decay of the autocovariances is, or what the nature of the design is. After presenting our theoretical results and illustrating them with numerous simulations, the question that arises is whether we can propose an effective, fully data-driven penalty. To answer this question, we first show - using simulations - that if we were to directly observe the errors, we could propose an effective penalty - denoted pen$_\varepsilon$ - inspired by the initial theoretical penalty (see Section \ref{subsec:FullyDataDriven}). Finally, we propose a two-step procedure: we select an initial model using a penalty that incorporates the Hurst coefficient of the process $(Y_i)_{1 \leq i \leq n}$, and then we use the residuals from this initial model to calculate the penalty. This two-step procedure appears to be very robust and often yields results close to the penalty pen$_\varepsilon$, reaching its limits only when both the error process and the design are  long-memory processes with very slow decay of correlations.

The paper is organized as follows. Section~\ref{sec:regPiece} presents the regression context and the notation used throughout the paper, while also introducing the dependence coefficients for random designs that are considered in subsequent sections. Section~\ref{sec3} provides the main results of the paper, beginning with risk bounds for short-memory error processes, including the case of heteroscedastic regression. Risk bounds for long-range dependence are then presented in the latter part of the section. Section~\ref{sec:main_longrange2} examines the optimality of the proposed penalty by studying the special case where the error process is a stationary Gaussian sequence independent of the design, and where the design belongs to a particular class of dependent sequences. Several penalization approaches are presented and compared in Section~\ref{sec:simus}. Finally, proofs of all results are provided in the appendix.

\section{Regression on a given interval via Piecewise polynomials}
\label{sec:regPiece}
\subsection{Statistical model and notation}

We consider the model \eqref{base_model}.
Our aim is to estimate the regression function $f^*$ on a given interval $I=[a,b]$ based on the observations $(X_i, Y_i)$ for which $X_i \in I$. We shall only consider estimators based on piecewise polynomials (although the results of Sections~\ref{sec:general_pen} and \ref{sec:main_shortrange} are true for general projection estimators). Let us define this class of estimators in more details. 

We consider a regular partition 
${\mathcal Q}$ of $I$ of size $m$ (with $m\leq n$). For $j \in \{1, \ldots , m\}$, let $Q_j$ be the $j$-th element of this partition. 
Let $r$ be a non-negative integer (representing the degree of the piecewise polynomials), and let $V_{j}$ be the subspace of  ${\mathbb R}^n$ (with dimension at most $r+1$) generated by the  $r+1$ rows 
\begin{equation} 
\label{vect_piecewise}
  v_j^k = (X_{1}^k  {\bf 1}_{X_{1} \in Q_{j}}, \cdots, X_{n}^k  {\bf 1}_{X_{n} \in Q_{j}})^t
\end{equation}
with $k \in \{0, ..., r\}$, i.e.
\begin{equation}\label{Vj}
   V_j = \mathrm{Vect}(v_j^0, v_j^1, \ldots, v_j^r) \, . 
\end{equation} 
Let then
$$S_{m}=V_{1} \oplus \cdots \oplus V_{m}.$$

The estimator $\hat f_m$ is then defined as the projection of the vector $Y=(Y_1, \ldots, Y_n)^t$ on $S_m$ with respect to the euclidean norm on ${\mathbb R}^n$. Note that $(\hat f_m)_i$ is expressed as a polynomial of degree $r$ in the variable $X_i$. Hence, in the following we use the notation: for any $i\in \{1, \dots, n\}$, $(\hat f_m)_i=\hat f_m(X_i)$.

Let $n({\mathbf X})= \sum_{i=1}^n {\mathbf 1}_{X_i \in I}$. For any $g$ from ${\mathbb R}$ to ${\mathbb R}$, let 
$$\|g\|_{n({\mathbf X})}^2= \frac{1}{n({\mathbf X})\vee 1}\sum_{i=1}^n (g(X_i))^2{\mathbf 1}_{X_i \in I} \, .$$

By convention, If $n({\mathbf X})=0$ we set $\|g\|_{n({\mathbf X})}^2=0$. Let also 
\begin{equation}\label{defAI}
A_I=\{i \in \{1, \ldots , n \} \ \text{s.t.} \ X_i \in I \}
\end{equation}
(in such a way that $n({\mathbf X})=\text{Card}(A_I)$). Note that the estimator $\hat f_m$ can be also defined as a projection estimator with respect to the norm $\|\cdot \|_{n({\mathbf X})}$ on an appropriate subspace of ${\mathbb R}^{n({\mathbf X})}$, which is deduced from $S_m$ by keeping only the coordinate $i \in A_I$ in the vectors introduced in \eqref{vect_piecewise}. With a slight abuse of notation, we still denote by $S_m$ this subspace of ${\mathbb R}^{n({\mathbf X})}$.

As in~\cite{caron2021gaussian}, we shall consider penalized estimators $\hat f_{\hat m}$, with
\begin{equation}
\label{hatm}
\hat m\in \text{argmin}_{m \in \{1,\ldots, n\}} \left \{  \frac{1}{n({\mathbf X})\vee 1}\sum_{i=1}^n (Y_i-\hat f_m(X_i))^2{\mathbf 1}_{X_i \in I}
 + \text{pen}(m) \right  \} \, ,
\end{equation}
where pen : $\{1, \ldots, n \}$ to a ${\mathbb R}^+$ is a penalty function. We shall work conditionally to the design ${\mathbf X}_n$, and describe different penalty functions, derived from the general form of the penalty given in Theorem 2.1 of~\cite{caron2021gaussian}.

We shall specify the rates of convergence of our projection estimators when $f^*$ belongs to a class of regular functions (regularity being defined on the interval $I$). Let $s > 0, L>0$, and let $k$ be the unique non-negative integer such that $s \in (k, k+1]$. We say that a function $g$ from $I$ to $\mathbb{R}$ belongs to the class ${\mathcal H}_s(L)$ when $g$ is $k$-times differentiable on $I$ and its $k$-th derivative $g^{(k)}$ satisfies 
$$
|g^{(k)}(x)-g^{(k)}(y)|\leq L |x-y|^{s-k} \, .
$$

\subsection{Dependence coefficients for  random designs} \label{subsec:coeff}
The main results of Sections \ref{sec:main_shortrange} and \ref{sec:main_longrange} (Theorem \ref{genth} and \ref{ThLD} respectively) are given conditionally to the design ${\mathbf X}_n$, and are valid without any assumptions on the design. However, when looking for some examples where some explicit rates can be computed, we need to make some assumptions on the design (this is also true for the special cases considered
in Section \ref{sec:main_longrange2}). 

In this subsection, we describe the dependence coefficients that we shall use, when the variables ${\mathbf X}_n$ come from a strictly stationary sequence $(X_i)_{i \geq 1}$. Let ${\mathcal I}$ be the collection of all non empty intervals in ${\mathbb R}$ (without loss of generality, on can restrict ourselves to closed intervals).
The dependence coefficients $(\alpha(k))_{k \geq 0}$ of the sequence $(X_i)_{i \geq 1}$ are defined as follows: 
\begin{equation}
\label{defalpha}
\alpha(k)= \sup_{I \in {\mathcal I}} \left |{\mathbb P}(\{X_1 \in I \}\cap  \{X_{k+1} \in I\})- ({\mathbb P}(X_1 \in I))^2 \right | \, .
\end{equation}
     Note that the coefficient $\alpha(k)$ is always smaller than the corresponding $\alpha$-mixing coefficient of  the sequence $(X_i)_{i \geq 1}$ (see \cite{Rosenblatt1956}). It is also smaller than $2\alpha_X(k)$, where $\alpha_X(k)$ is the dependence coefficient defined in \cite{DedeckerPrieur2005} (we refer to this paper for many examples for which these coefficients can be computed). 
The dependence coefficients $(\rho(k))_{k \geq 0}$ of the sequence $(X_i)_{i \geq 1}$ are defined as follows: 
\begin{equation}
\label{defrho}
\rho(k)= \sup_{I \in {\mathcal I}} \left | \frac{\text{Cov}({\bf 1}_{X_1 \in I}, {\bf 1}_{X_{k+1} \in I})}{\text{Var}({\bf 1}_{X_1 \in I})}\right | \, ,
\end{equation}
where we use the convention $0/0=0$ in \eqref{defrho}. Note that we always have that $\alpha(k) \leq \rho(k)$, and that 
the coefficient $\rho(k)$ is always smaller than the corresponding $\rho$-mixing coefficient of  the sequence $(X_i)_{i  \geq 1}$. Another interesting example is the following: if $(Y_k)_{k \geq 1}$ is a stationary Gaussian sequence with auto-covariance $(\gamma_Y(k))_{k \geq 0}$, and if $X_k=f(Y_k)$ for some measurable function $f$, then 
$$\rho(k) \leq \frac{\left |\gamma_Y(k) \right |}{\gamma_Y(0)} \, . 
$$
This follows from \cite{Lancaster1957} (see  \cite{Demboetal2001}  for more details). 

\section{Bounding the conditional MSE of penalized estimators} \label{sec3}
\subsection{A general penalty}
\label{sec:general_pen}


We denote by $
\rho_{{\mathbf e} |{\mathbf X}}$ the spectral radius of the matrix $\Sigma({\mathbf X}_n)$, and
 by $
\rho_{{\mathbf e} |{\mathbf X}, I}$ the spectral radius of the matrix $\Sigma_I({\mathbf X}_n)$, where $\Sigma_I({\mathbf X}_n)$ is the sub-matrix of  $\Sigma({\mathbf X}_n)$ in which  the raws and columns corresponding to the indexes $i\notin A_I$  have been suppressed (so that
$\Sigma_I({\mathbf X}_n)$ is a $n({\mathbf X})\times n({\mathbf X})$ matrix). By convention, If $n({\mathbf X})=0$ we set $\rho_{{\mathbf e} |{\mathbf X}, I}=0$.

In Theorem 2.1 of~\cite{caron2021gaussian}, we give a general penalty in a Gaussian context. By conditioning on the design ${\mathbf X}_n$, the result of this theorem asserts that one can choose any penalty such that
\begin{equation}
\label{GenPen}
\text{pen}(m) \geq \frac{K}{n({\mathbf X})\vee 1} \left (\sqrt{\text{tr}\left ( \text{Proj}_{S_m} \Sigma_I({\mathbf X}_n)\right ) + \rho_{{\mathbf e} |{\mathbf X}, I} } + \sqrt{\rho_{{\mathbf e} |{\mathbf X}, I}} \sqrt{2\log \left ( \frac 1 {\pi_m} \right )} 
\right )^2 \, ,
\end{equation}
where $K>1$, Proj$_{S_m}$ is the $n({\mathbf X})\times n({\mathbf X})$ matrix corresponding to the orthogonal projection on $S_m$, and $(\pi_m)_{m \in \{1, \ldots n\}}$ are some weights such that $\pi_1+\cdots + \pi_n = 1$. Here, we shall always choose the weights $\pi_m$ of order $m^{-2}$, so that $\log (\pi_m^{-1})$ is of order $\log(m)$.

\subsection{Short range dependence} 
\label{sec:main_shortrange}

From the discussion on Page 4831 in~\cite{caron2021gaussian}, and since dim$(S_m) \leq (r+1)m$, we see that one can choose a penalty of the form 
\begin{equation}
\label{pen1_sd}
\text{pen}(m)= \kappa  \rho_{{\mathbf e} | {\mathbf X},I}\frac{m}{n({\mathbf X}) \vee 1}
\end{equation}
where $\kappa >K_r$, for a constant $K_r$ depending only of the degree $r$. Hence, if we define $\hat m$ via \eqref{hatm} and the penalty \eqref{pen1_sd}, we infer from Theorem 2.1 of~\cite{caron2021gaussian} that:

\begin{Theorem} 
\label{genth}
Let $\kappa>K_r$, and let $\hat m$ be defined via \eqref{hatm} and the penalty \eqref{pen1_sd}. Then, almost surely, 
\begin{equation}
\mathbb{E} \left[ \left \| f^{\ast} - \hat{f}_{\hat m} \right \|_{n({\mathbf X})}^{2} \Big |{\mathbf X}_n\right] \leq C \left( \inf_{m} \left\{ \mathbb{E} \left[ \left \| f^{\ast} - \hat{f}_{m}  \right \|_{n({\mathbf X})}^{2} \Big | {\mathbf X}_n \right] + \kappa \rho_{{\mathbf e} | {\mathbf X},I}\frac{m}{n({\mathbf X}) \vee 1} \right\} + \frac {\rho_{{\mathbf e} | {\mathbf X},I}}{n({\mathbf X})\vee 1} \right) \, ,
\end{equation}
where $C>1$ depends only on $\kappa$.
\end{Theorem}

Note that the penalty \eqref{pen1_sd} depends on the spectral radius $\rho_{{\mathbf e} | {\mathbf X},I}$, which is, of course, unknown in practice. This is why it is necessary to use an algorithm that allows us to calibrate the penalty constant in front of the term in $m/n({\mathbf X})$. We will use the dimension jump algorithm described in Section 5.1.

Based on this general result, we can easily obtain convergence rates when $f^*{\bf 1}_I$ belongs to the class $H_s(L)$.

\begin{Corollary}
\label{cor1_sd}
Let $\kappa>K_r$, and let $\hat m$ be defined via \eqref{hatm} and the penalty \eqref{pen1_sd}.
Assume that $f^*{\bf 1}_I$ belongs to the class ${\mathcal H}_s(L)$, for some $s \in (0, r+1]$. Then there exist  positive constants $C_{r,L}$ and $C_r$ such that 
$$
\mathbb{E} \left[ \left \| f^{\ast} - \hat{f}_{\hat m} \right \|_{n({\mathbf X})}^{2} \Big |{\mathbf X}_n\right] \leq \frac {C_{r,L}\rho_{{\mathbf e} | {\mathbf X}, I}^{\frac{2s}{2s+1}}}{(n({\mathbf X})\vee 1)^{\frac{2s}{2s+1}}}+ 
\frac {C_r\rho_{{\mathbf e} | {\mathbf X},I}}{n({\mathbf X})\vee 1} \ \text{almost surely}.
$$
If moreover 
the variables $X_i$'s take their values in the interval $I$, then $n({\mathbf X})=n$ and
\begin{equation}\label{11}
\mathbb{E} \left[ \left \| f^{\ast} - \hat{f}_{\hat m} \right \|_{n}^{2} \right] \leq \frac {C_{r,L}{\mathbb E}\big (\rho_{{\mathbf e} | {\mathbf X}}^{\frac{2s}{2s+1}}\big )}{n^{\frac{2s}{2s+1}}} + \frac {C_r{\mathbb E}(\rho_{{\mathbf e} | {\mathbf X}})}{n}\, .
\end{equation}
\end{Corollary}

The upper bound \eqref{11} gives a rate of order $n^{-\frac{2s}{2s+1}}$ for $\mathbb{E} [  \| f^{\ast} - \hat{f}_{\hat m}  \|_{n({\mathbf X})}^{2} ]$ provided ${\mathbb E}(\rho_{{\mathbf e} | {\mathbf X}}) \leq c$ for some $c>0$ and the $X_i$'s take their values in $I$. Now, if $\rho_{{\mathbf e} | {\mathbf X}, I} \leq \rho
$ for some $\rho >0$, we can obtain a similar control by assuming that  $(X_i)_{i \geq 1}$ is a strictly stationary $\alpha$-dependent sequence with summable coefficients and ${\mathbb P}(X_1 \in I)>0$. 

\begin{Corollary}\label{cor2_sd_alpha} Let $\kappa>K_r$, and let $\hat m$ be defined via \eqref{hatm} and the penalty \eqref{pen1_sd}.
Assume that $f^*{\bf 1}_I$ belongs to the class ${\mathcal H}_s(L)$, for some $s \in (0, r+1]$, and that the sequence $(X_i)_{i \geq 1}$ is a strictly stationary  $\alpha$-dependent sequence such that ${\mathbb P}(X_1 \in I)>0$.  Assume moreover that there exist a positive constant $\rho$ such that $\rho_{{\mathbf e} | {\mathbf X}, I} \leq \rho
$ almost surely. 
\begin{enumerate}
\item If $\sum_{k>0} \alpha (k) < \infty$, then
$$
\mathbb{E} \left[ \left \| f^{\ast} - \hat{f}_{\hat m} \right \|_{n({\mathbf X})}^{2} \right] = O\left(n^{\frac{-2s}{2s+1}}\right)\, .
$$
\item If $\alpha(k)= O(k^{-a})$ for some $a \in (0,1)$, then
$$
\mathbb{E} \left[ \left \| f^{\ast} - \hat{f}_{\hat m} \right \|_{n({\mathbf X})}^{2} \right] = O\left( n^{-\min\left (a,\frac{2s}{2s+1} \right )}\right)\, .
$$
\end{enumerate}
\end{Corollary}

\subsection{Heteroscedastic regresssion for i.i.d. couples}
\label{subsec:heterosceda}
Corollary \ref{cor1_sd} has an interesting application to heteroscedastic regression. Assume that the sequence $(X_i, \varepsilon_i)_{1 \leq i \leq n}$ is i.i.d., and let $(X, \varepsilon)$ be a random couple with the same distribution as $(X_1, \varepsilon_1)$. Assume also that the conditional distribution of $\varepsilon$ given $X$ is a normal distribution ${\mathcal N}(0, \sigma^2(X))$. In this context, we have that 
$$
\rho_{{\mathbf e}| {\mathbf X},I}= \max_{1 \leq i \leq n} \sigma^2(X_i) {\mathbf 1}_{X_i \in I} \, .
$$
As a consequence of Corollary \ref{cor1_sd}, we get:
\begin{Corollary}
\label{cor_heteroiid} Let $\kappa>K_r$, and let $\hat m$ be defined via \eqref{hatm} and the penalty \eqref{pen1_sd}.
Assume that $f^*{\bf 1}_I$ belongs to the class ${\mathcal H}_s(L)$, for some $s \in (0, r+1]$. Assume also that ${\mathbb P}(X_1 \in I)>0$. Then, under the assumptions of this subsection, 
\begin{enumerate}
\item If $\sigma^2(X){\mathbf 1}_{X \in I}$ belongs to ${\mathbb L}^p$ for $p>1$,
$$
\mathbb{E} \left[ \left \| f^{\ast} - \hat{f}_{\hat m} \right \|_{n({\mathbf X})}^{2} \right] = o\left( n^{\frac{-2s(p-1)}{(2s+1)p}}\right)\, .
$$
\item If $\sigma^2(X){\mathbf 1}_{X \in I}$ has an exponential moment, 
$$
\mathbb{E} \left[ \left \| f^{\ast} - \hat{f}_{\hat m} \right \|_{n({\mathbf X})}^{2} \right] = o\left( (\log (n))^{\frac{2s}{2s+1}} n^{\frac{-2s}{2s+1}}\right)\, .
$$
\item If $\|\sigma^2(X){\mathbf 1}_{X \in I}\|_\infty < \infty $, 
$$
\mathbb{E} \left[ \left \| f^{\ast} - \hat{f}_{\hat m} \right \|_{n({\mathbf X})}^{2} \right] = O\left(  n^{\frac{-2s}{2s+1}}\right)\, .
$$
\end{enumerate}
\end{Corollary}

\begin{Rem}
Note that, if the function $x \mapsto \sigma^2(x)$ is continuous on $I$, then the condition of Item 3 of Corollary \ref{cor_heteroiid} is satisfied. 
\end{Rem}

\subsection{Long range dependence} \label{sec:main_longrange}

  We assume now that, for any $(i, j) \in A_I\times A_I$
\begin{equation}
\label{decaygamma}
 |\mathrm{Cov}(\varepsilon_i, \varepsilon_{j}| {\mathbf X}_n)| \leq \frac{c_I({\mathbf X}_n)}{ (|i-j|+1)^{\gamma}}  \, ,
\end{equation}
where $\gamma \in (0,1)$ and $c_I({\mathbf X}_n) < \infty$ almost surely. If 
the variables $X_i$'s take their values in the interval $I$, we simply write $c_I({\mathbf X}_n)=c({\mathbf X}_n)$. By convention, If $n({\mathbf X})=0$ we set $c_I({\mathbf X}_n)=0$.

If we only assume that \eqref{decaygamma} holds, the penalty \eqref{pen1_sd} may be suboptimal (as in the case of regular fixed design, see Section 3.2 in~\cite{caron2021gaussian}) because the spectral radius may not be bounded in $n$ (this situation is what we call 
``long range dependence''). We shall see that another possible penalty is given by
\begin{equation}
\label{pen2_ld}
\text{pen}(m)= \kappa c_I({\mathbf X}_n) \left (\frac{m}{n({\mathbf X})\vee 1}\right )^\gamma
\end{equation}
where $\kappa >K_{r, \gamma}$, for a constant $K_{r, \gamma}$ depending only of the degree $r$ and on the exponent $\gamma$. Hence, if we define $\hat m$ via \eqref{hatm} and the penalty \eqref{pen2_ld}, we can prove the following result:

\begin{Theorem} 
\label{ThLD} 
Let $\kappa>K_{r, \gamma}$, and let  $\hat m$ be defined via \eqref{hatm} and  the penalty \eqref{pen2_ld}. Then, almost surely, 
\begin{multline}
\mathbb{E} \left[ \left \| f^{\ast} - \hat{f}_{\hat m} \right \|_{n({\mathbf X})}^{2}  \Big |{\mathbf X}_n\right] \\ 
\leq C \left( \inf_{m} \left\{ \mathbb{E} \left[ \left \| f^{\ast} - \hat{f}_{m}  \right \|_{n({\mathbf X})}^{2} \Big | {\mathbf X}_n \right] + \kappa c_I({\mathbf X}_n) \left (\frac{m}{n({\mathbf X})\vee 1}\right )^\gamma \right\} + \frac {\rho_{{\mathbf e} | {\mathbf X},I}}{n({\mathbf X})\vee 1} \right) \, ,
\end{multline}
where $C>1$ depends only on $\kappa$.
\end{Theorem}
 
 As in the previous sections, we can easily obtain convergence rates when $f^*{\bf 1}_I$ belongs to the class $H_s(L)$. 
 
\begin{Corollary}
\label{cor1_ld}
Let $\kappa>K_{r, \gamma}$, and let  $\hat m$ be defined via \eqref{hatm} and  the penalty \eqref{pen2_ld}.
Assume that $f^*{\bf 1}_I$ belongs to the class ${\mathcal H}_s(L)$, for some $s \in (0, r+1]$. Then there exist  positive constants $C_{r,\gamma, L}$ and $C_{r, \gamma}$ such that 
$$
\mathbb{E} \left[ \left \| f^{\ast} - \hat{f}_{\hat m} \right \|_{n({\mathbf X})}^{2} \Big |{\mathbf X}_n\right] \leq \frac {C_{r, \gamma, L} c_I({\mathbf X}_n)^{\frac{2s}{2s+\gamma}}}{(n({\mathbf X})\vee 1)^{\frac{2\gamma s}{2s+\gamma}}} + \frac {C_{r, \gamma}\rho_{{\mathbf e} | {\mathbf X}, I}}{n({\mathbf X})\vee 1} 
\ \text{almost surely}.
$$
If moreover 
the variables $X_i$'s take their values in the interval $I$, then $n({\mathbf X})=n$ and
$$
\mathbb{E} \left[ \left \| f^{\ast} - \hat{f}_{\hat m} \right \|_{n}^{2} \right] \leq  
\frac {C_{r, \gamma, L} \mathbb{E} \big ( c_I({\mathbf X}_n)^{\frac{2s}{2s+\gamma}}\big)}{n^{\frac{2\gamma s}{2s+\gamma}}} + \frac {C_{r, \gamma}\mathbb{E} (\rho_{{\mathbf e} | {\mathbf X}, I})}{n} 
\, .
$$
\end{Corollary}

Under the same assumptions regarding  $(X_i)_{i \geq 1}$ as in Corollary \ref{cor2_sd_alpha}, we obtain the following result.

\begin{Corollary} 
\label{cor2_ld_alpha} Let $\kappa>K_{r, \gamma}$, and let  $\hat m$ be defined via \eqref{hatm} and  the penalty \eqref{pen2_ld}.
Assume that $f^*{\bf 1}_I$ belongs to the class ${\mathcal H}_s(L)$, for some $s \in (0, r+1]$, and that the sequence $(X_i)_{i \geq 1}$ is a strictly stationary  $\alpha$-dependent sequence such that ${\mathbb P}(X_1 \in I)>0$.  Assume moreover that there exist a positive constant $c$ such that $c_I({\mathbf X}_n) \leq  c
$ almost surely. 
\begin{enumerate}
\item If $\sum_{k>0} \alpha (k) < \infty$, then
$$
\mathbb{E} \left[ \left \| f^{\ast} - \hat{f}_{\hat m} \right \|_{n({\mathbf X})}^{2}  \right] = O\left(n^{\frac{-2\gamma s}{2s+\gamma}}\right)\, .
$$
\item If $\alpha(k)= O(k^{-a})$ for some $a \in (0,1)$, then
$$
\mathbb{E} \left[ \left \| f^{\ast} - \hat{f}_{\hat m} \right \|_{n({\mathbf X})}^{2}  \right] = O\left( n^{-\min \left(a,\frac{2 \gamma s}{2s+\gamma} \right )}\right)\, .
$$
\end{enumerate}
\end{Corollary}

\section{Discussion on the penalty}
\label{sec:main_longrange2}

The penalty \eqref{pen2_ld} is appropriate in some cases (such as in case of regular fixed design, see Section 3.2 in~\cite{caron2021gaussian}), but it may be suboptimal for some random designs (see Experiment 11 in Section \ref{subsec533}). In this section, we examine the penalty in more detail in the case of regular regressograms (i.e. when $r=0$), and when the error process $(\varepsilon_i)_{i \geq 1}$ is a stationary Gaussian sequence  independent of the design $(X_i)_{i \geq 1}$. For the sake of simplicity, we also assume that the $X_i$'s take their values in the interval $I$.

Let $\Sigma_n$ be the covariance matrix of the vector $(\varepsilon_1, \ldots, \varepsilon_n)^t$ and $\rho_{{\mathbf e},n}$ be the spectral radius of $\Sigma_n$.
According to Theorem 2.1 of~\cite{caron2021gaussian}), we can choose a penalty satisfying
\begin{equation}
\label{GenPen_ld}
\text{pen}(m) \geq \frac{K}{n} \left (\sqrt{\text{tr}\left ( \text{Proj}_{S_m} \Sigma_n \right ) + \rho_{{\mathbf e},n}} + \sqrt{\rho_{{\mathbf e},n}} \sqrt{2\log \left ( \frac 1 {\pi_m} \right )}
\right )^2 \, ,
\end{equation}
where $K>1$. 
Note first that
$$
\text{tr}\left ( \text{Proj}_{S_m} \Sigma_n \right )=n{\mathbb E} \left (  \left \|\text{Proj}_{S_{m}}(\varepsilon) \right \|_n^2   \big | {\mathbf X}_n \right )\, .
$$
Since $r=0$, we have (with the usual convention that $0/0=0$)
$$
\|\text{Proj}_{S_{m}}(\varepsilon)\|_n^2=  \frac 1 n\sum_{i=1}^{m}   \ell_{i}(\bar \varepsilon_{i})^2 \, ,
$$
where
$$
   \bar \varepsilon_{i} = \frac{1}{\ell_{i}}  \sum_{k=1}^{n}  \varepsilon_{k} 
   {\bf 1}_{X_{k} \in Q_{i}}\, , \quad  \ell_{i}=  \sum_{k=1}^{n}  
   {\bf 1}_{X_{k} \in Q_{i}} \, .
$$
Let $\gamma_{\varepsilon}(k)= \mathrm{Cov}(\varepsilon_1, \varepsilon_{1+k})$ be the auto-covariance of $(\varepsilon_i)_{i \geq 1}$. It follows that 
$$
\text{tr}\left ( \text{Proj}_{S_m} \Sigma_n \right )=n{\mathbb E} \left (  \left \|\text{Proj}_{S_{m}}(\varepsilon) \right \|_n^2   \big | {\mathbf X}_n \right ) = \sum_{i=1}^m \frac{1}{\ell_i} \sum_{j=1}^n \sum_{k=1}^n \gamma_{\varepsilon}(j-k) {\bf 1}_{X_j \in Q_i} {\bf 1}_{X_k \in Q_i} 
\, .
$$
We see that the main term of the penalty $\text{tr}\left ( \text{Proj}_{S_m} \Sigma_n \right )$ can be computed from the data if we know the expression of the auto-covariance $\gamma_\varepsilon$. Moreover, if  $\sup_{n >0}\rho_{{\mathbf e},n} \leq \rho < \infty$, then $\text{tr}\left ( \text{Proj}_{S_m} \Sigma_n \right )$ is of order $m$ (see the upper bound (2.4) in \cite{caron2021gaussian}), and the penalty pen($m$) can be chosen as pen($m$)=$\kappa m/n$, with  $\kappa> K_\rho$, for $K_\rho$ depending on $\rho$. This is in accordance with the results of Section \ref{sec:main_shortrange}.

Let us now examine the case where 
\begin{equation}
\label{gammak}
|\gamma_{\varepsilon} (k)| \leq \frac{c}{ (k+1)^{\gamma}}  \, ,
\end{equation}
where $\gamma \in (0,1)$ and $c>0$. In that case $
\rho_{{\mathbf e},n} \leq \kappa_{\gamma, c} n^{1-\gamma}
$, and according to Section \ref{sec:main_longrange}, we can choose 
\begin{equation}
\label{pen_ld}
\text{pen}(m)= \kappa  \left (\frac{m}{n}\right )^\gamma
\end{equation}
where $\kappa >K_{\gamma, c}$, for a constant $K_{\gamma, c}$ depending on $(\gamma, c)$. As shown in Section \ref{PrThLD}, this form of penalty arises from a uniform control of  $\text{tr}\left ( \text{Proj}_{S_m} \Sigma_n \right )$, which holds regardless of the structure of $(X_i)_{i \geq 1}$. 
However, the Experiment 11 of Section \ref{subsec533}  shows that this penalty term is sub-optimal in this particular case. 
To understand this phenomenon, we shall give a control of ${\mathbb E}\left(\text{tr}\left ( \text{Proj}_{S_m} \Sigma_n \right )\right )$, which is much more precise than the uniform bound obtained in Section \ref{PrThLD}. This control will be valid under some additional conditions on the sequence $(X_i)_{i \geq 1}$.


\begin{Proposition} 
\label{prop_ld_improved}
Assume that $\gamma_\varepsilon$ satisfies \eqref{gammak}, and that the sequence $(X_i)_{i \geq 1}$ is strictly stationary and $\rho$-dependent,  with  coefficients $(\rho(k))_{k \geq 0}$. Assume moreover that $D=\sum_{k=0} |\gamma_\varepsilon (k)| \rho(k) < \infty$. Then 
\begin{equation*}
\frac 1 n {\mathbb E} \left ( \mathrm{tr}\left ( \mathrm{Proj}_{S_m} \Sigma_n \right ) \right ) ={\mathbb E} \left (  \left \|\mathrm{Proj}_{S_{m}}(\varepsilon) \right \|_n^2 \right ) \leq B\frac{m}{n} + \frac{2K_{\gamma, c}}{n^{\gamma}} \, ,
 \end{equation*}
 where the constant $B$ depends on $(\gamma_{\varepsilon}(0), D)$ and
 $K_{\gamma, c}$ depends on $(\gamma, c)$.
\end{Proposition}

Recall  that, under \eqref{gammak},  $
\rho_{{\mathbf e},n} \leq \kappa_{\gamma, c} n^{1-\gamma}
$,  and that $\log(\pi_m^{-1})$ is of order $\log(m)$.
Although this is just an heuristic, Proposition \ref{prop_ld_improved} suggests that, instead of \eqref{pen_ld}, the correct form of the penalty should be 
\begin{equation}
\label{pen_ld_improved}
\text{pen}(m)= \kappa_1 \frac m n + \kappa_2 \frac{\log(m)}{n^\gamma}
\end{equation}
where $\kappa_1 >C$ and $\kappa_2> C_{\gamma, c}$, for $C$ depending on $(\gamma_{\varepsilon}(0), D)$ and $C_{\gamma, c}$ depending on $(\gamma, c)$. As the simulations show, the penalty \eqref{pen_ld_improved} is indeed appropriate for the Experiment 11 of Section \ref{subsec533} (for which the assumptions of Proposition \ref{prop_ld_improved} hold).

\section{Numerical experiments}
\label{sec:simus}

In this section, we validate the penalties described by the theoretical results of the paper. To calibrate these penalties, we employ the slope heuristic, partially drawing on the protocols outlined in\cite{caron2021gaussian}. To limit the computational cost of the procedures, we present the experiments for regressograms (with $r=0$).

\subsection{Slope heuristics} 
\label{subsec:slopeh}

For the results given in the previous sections, the penalty functions are known, in the best case, up to a multiplicative constant. The aim of the slope heuristics method proposed by Birg\'e and  Massart~ \cite{birge2007minimal} is precisely to calibrate a penalty function for model selection purposes. See \cite{baudry2012slope} and ~\cite{arlot2019minimal} for a general presentation of the method, and also~\cite{caron2021gaussian} for a presentation of the algorithm in a context very similar to that of this paper.
This method has shown very good performances and comes with mathematical guarantees for nonparametric Gaussian regression with i.i.d. error terms,  see \cite{birge2007minimal,arlot2019minimal} and references therein. The slope heuristics have several versions (see ~\cite{arlot2019minimal}). In this paper we use the dimension jump algorithm, which is implemented for instance in the R package {\ttfamily capushe}. 

The aim is to tune the constant $\kappa$ in a penalty of the form $\pen(m) = \kappa \pen_{\tiny \mbox{shape}}(m)$ where $\pen_{\tiny \mbox{shape}}$ is a known penalty shape. In the most standard cases, $\pen_{\tiny \mbox{shape}}$  is the dimension of the model. Let  $\hat{m}(\kappa)$
\[\hat{m}(\kappa) \in \mathrm{argmin}_{m \in \mathcal{M}} \left\{ \left \| Y - \hat{f}_{m} \right \|_n^{2} +  \kappa \pen_{\tiny \mbox{shape}} (m) \right\}\]
 be the model selected by the penalized criterion with constant $\kappa$.
The Dimension Jump algorithm consists of the following steps
\begin{enumerate}
\item  Compute $ \kappa  \mapsto  \hat{m}( \kappa )$,
\item Find the constant ${\hat \kappa}^{dj} > 0$ that corresponds to the highest jump of the function $\kappa \rightarrow  \hat{m}(\kappa)$,
\item Select the model $\hat{m}(2 {\hat \kappa}^{dj})$, 
  $$\hat{m} \in \mathrm{argmin}_{m \in \mathcal{M}} \left\{ \left \| Y - \hat{f}_{m} \right \|_n^{2} + 2 {\hat \kappa}^{dj} \pen_{\tiny \mbox{shape}} (m)\right\} .$$
\end{enumerate}

The penalty shapes discussed earlier in this paper depend on unknown quantities, so the challenge extends beyond simple calibration via a multiplicative factor. To avoid addressing all issues simultaneously, we first consider the idealized setting where the generative models—and thus the penalty shapes—are perfectly known. We then tackle the general case in Section~\ref{subsec:FullyDataDriven}, where the stochastic nature of both the error and design processes remains unknown.

\subsection{Presentation of the experiments}

We simulate $n$ observations according to several  generative model on $[0,1]$ satisfying the equation
\begin{equation}
\label{simus_model}
Y_i  =  f^\ast \left( X_i \right) + \varepsilon_{i}, \quad i = 1 \dots n  ,
\end{equation}
with $ f^\ast $ being one of the two functions
\begin{align*}
f_1^\ast : t \in [0,1]   &\mapsto   3 - 0.1 t + 0.5 t^2 - t^3 + \sin(8 t) \\
\end{align*}
and
\begin{align*}
f_2^\ast  : t \in [0,1]   &\mapsto 
\begin{cases}
1 + x + \sin(10x)  & \text{if } x \leq \dfrac{2}{3}, \\[0.5em]
\sqrt{x-\dfrac{2}{3}} + 1 + \dfrac{2}{3} + \sin\!\left(\dfrac{20}{3}\right) 
& \text{if } x > \dfrac{2}{3}.
\end{cases}
\end{align*}

The aim is to estimate $f^\ast $ based on the observations of $(X_i, Y_i)$ on a regular partition of size $m$, for $m \in \{1, \ldots, \min(n/2,500)\}$. Please note that, for the sake of clarity, the risk curves shown below do not always depict the largest models. 

We use the same function \( f_1^\ast \) as in our previous paper~\cite{caron2021gaussian}, in order to facilitate comparison with the deterministic design setting considered in our previous work.
The second function \( f_2^\ast \)  exhibits a significant singularity for \(x = 2/3 \) and thus allows us to consider a regression function that is far less smooth than in the case of \( f_1^\ast \).

We simulate $n$ observations $(X_i,\varepsilon_i)$ where the design $X_1, \dots,X_n$ and the errors $\varepsilon_1, \dots,\varepsilon_n$ follows various distributions. Moreover, $(X_i)$ and $(\varepsilon_i)$ are assumed to be independent in most cases,  but we also consider the case where they are not (heteroscedastic setting).

We simulate $n$ observations $\varepsilon_i$ according to AR processes, Fractional Gaussian processes and  non Gaussian Markov chains. We also simulate $n$ observations $X_i$ according to deterministic designs, uniform designs, (functions of) Fractional Gaussian processes and  non Gaussian Markov chains. 

Note that we do not always simulate $\varepsilon_i$ according to a Gaussian process. This allows us to evaluate the robustness of the model selection procedure without the Gaussian assumption on the errors. 

We shall consider samples of size  $n=500$ and $n=2000$; for each estimator that will be considered below, the boxplots of the risks will be carried out through 100 independent trials. 

\medskip

We now give more details on the processes we use for the simulations.
\begin{itemize}
\item[•] {\bf Non mixing auto-regressive model.}
We  consider random variables $(X_1, \ldots, X_n)$,  generated according to the simple 
AR(1) equation
$$
\text{for $k \geq 1$,} \quad X_{k+1}=\frac12 \left ( X_k + e_{k+1} \right ) \, ,
$$
where $X_1$ is uniformly distributed over $[0,1]$,
and $(e_i)_{i \geq 2}$ is a sequence of 
i.i.d. random variables with distribution ${\mathcal B}(1/2)$,
independent of $X_1$.  

One can check that the transition kernel of the  chain  $(X_i)_{i \geq 1}$ is
$$
 K(f)(x)= \frac 1 2 \left (f\left (\frac x 2\right )+ f\left (\frac{x+1}{2} \right ) \right ) \, ,
$$
and that the uniform distribution on $[0,1]$ is the unique
invariant distribution by $K$. Hence, the chain $(X_i)_{i \geq 1}$
is strictly stationary. 

It is well known that the chain $(X_i)_{i \geq 1}$ is not $\alpha$-mixing 
in the sense of  \cite{Rosenblatt1956} (see for instance  \cite{Andrews1984}). 
However, one can  prove that the coefficients $\phi_X(k)$ (see Definition 2 in  \cite{DedeckerPrieur2005})
of  $(X_i)_{i \geq 1}$ are such that 
$$
   \phi_X(k) \leq 2^{-k} \, .
$$
Now, the coefficients  $\rho(k)$ defined in Section \ref{subsec:coeff} are such that $\rho(k) \leq 2 \sqrt{2\phi_X(k)}$ (this can be easily deduced from Proposition 2.1 in \cite{Dedecker2004}). It follows that the coefficients  $ \rho(k)$ of $(X_i)_{i \geq 1}$ decrease at an exponential rate. 

In the simulations, the error process is obtained by centering the chain, namely
\begin{equation}
    \varepsilon_i = X_i - \frac{1}{2}.
\label{Nonmix_ar1}    
\end{equation}


\item[•] \textbf{Fractional Gaussian Noise.} 
The Fractional Gaussian Noise (FGN, see for instance \cite{MvN68} and \cite{Ber94}) is a stationary sequence $(\varepsilon_i)_{i \geq 1}$ of zero-mean Gaussian random variables with auto--covariances 
$$
\gamma_\varepsilon (k)= \frac{\sigma^2}{2} \left ( |k+1|^{2H}-2|k|^{2H} +|k-1|^{2H}  \right), \quad \text{for } k\in \mathbb N,
$$
where $\sigma^2= \gamma_\varepsilon(0)= \text{Var}(\varepsilon_i)$,  and $H \in (0, 1)$ is the so-called Hurst parameter. If  $H=1/2$, the sequence $(\varepsilon_i)_{i \geq 1}$ is   a Gaussian white noise with variance $\sigma^2$. For any $H \in (0,1)$ the following asymptotic expansion is valid
$$
    \gamma_\varepsilon (k) \sim \sigma^2 H(2H-1) k^{2(H-1)} \, . 
$$
Consequently, if $H>1/2$, the process is positively correlated and long-range dependent. If $H<1/2$, the process is negatively correlated and   $\sum_{k\geq 0} |\gamma_\varepsilon (k) | < \infty$. Let $\Sigma_n$ be the covariance matrix of the vector $(\varepsilon_1, \ldots, \varepsilon_n)^t$ and $\rho_{{\mathbf e},n}$ be the spectral radius of $\Sigma_n$. If $H<1/2$, then the spectral radius $\rho_{{\mathbf e},n}$ is bounded by $\sum_{k \in {\mathbb Z}}|\gamma_\varepsilon (k) | $ and the process is short-range dependent.

Now, as mentioned in Section \ref{subsec:coeff}, if $X_k=f(\varepsilon_k)$ for some measurable function $f$, then the coefficients $ \rho(k)$ of $(X_i)_{i \geq 1}$ are such that
$$\rho(k) \leq \frac{\left |\gamma_\varepsilon(k) \right |}{\gamma_{\varepsilon}(0)} \, . 
$$
In the simulations we take $\sigma^2=1$, and the design process is obtained by applying the cumulative distribution function $\Phi$ of the ${\mathcal N}(0,1)$-distribution, namely
\begin{equation}
    X_i = \Phi(\varepsilon_i) \, ,
\label{transFGN}    
\end{equation}
so that $X_i$ is uniformly distributed over $[0,1]$.

\item[•] \textbf{Non Gaussian Markov chain.} 
We start from the Markov chain introduced by Doukhan, Massart and Rio \cite{DMR94}. 

Let $a$ be a positive real number, 
let $\nu$ be the probability with density $x \rightarrow (1+a)x^a {\bf 1}_{[0,1]}$ and $\pi$ be the probability with density $x \rightarrow ax^{a-1}{\bf 1}_{[0,1]}$. We define now a strictly stationary Markov chain by specifying its transition
probabilities $K(x,A)$ as follows
\[ K(x,A)=(1-x)\delta_{x}(A)+x\nu(A)\, ,\]
where $\delta_{x}$ denotes the Dirac measure at point $x$. Then $\pi$ is the unique invariant probability measure of the chain with transition probabilities $K(x, \cdot)$. Let $(Z_i)_{i \in {\mathbb Z}}$ be the stationary Markov chain on $[0,1]$ with transition probabilities $K(x, \cdot)$ and invariant distribution $\pi$. 
From  \cite{DMR94}, we know that the $\beta$-mixing coefficients  $\beta_Z(n)$ of the chain are such that
 $\beta_Z(n)\sim \frac{1}{n^a}$. One can easily check than 
\begin{equation}
    X_i=Z_i^a 
    \label{markov_chain_X}
\end{equation}
is uniformly distributed over $[0,1]$ with $\beta_X(n)\sim \frac{1}{n^a}$,
 so that 
\begin{equation}
    \varepsilon_i= Z_i^a -0.5
    \label{markov_chain_eps}
\end{equation}

is a stationary Markov chain (as an invertible function of a stationary Markov chain), with mean zero and mixing coefficient $\beta_\varepsilon (n) \sim \frac{1}{n^a}$. This chain is short range dependent if $a>1$ and long-range dependent if $a \in (0,1)$ (see for instance \cite{DGM18} for a deeper discussion on this subject). 
For $a\in (0,1)$, the corresponding Hurst coefficient is $H=1-a/2$.

Regarding the coefficients defined in Section \ref{subsec:coeff}, the coefficients $\alpha(k)$ of the chain $(X_i)_{i \geq 1}$ are such that $\alpha(k) \leq \beta_X(k)$. However,  one cannot control the coefficients $\rho(k)$ of  $(X_i)_{i \geq 1}$.

\end{itemize}

\paragraph{Penalization approaches.}  
In the experiments, we compare several penalization approaches:
\begin{itemize}
    \item {\bf CDJ}: Classical Dimension Jump method with a penalty shape proportional to the dimension.
    \item {\bf HGiven}: Dimension Jump method with a penalty shape \( m^{2-2H} \), where the Hurst exponent \( H \) of the error process is given.
    \item  $\mathbf{pen_{\varepsilon}}$: Dimension Jump method with a penalty shape which is a proxy for the variance bound $ \mathbb{E} \left[\left \| \text{Proj}_{S_m}(\varepsilon) \right \|^{2} \big | {\mathbf X}_n \right]$, as presented in Section~\ref{subsec:FullyDataDriven}. We display this penalty in all risk-comparison figures for reference. As explained above, this quantity is, of course, not known in practice. We calculate it for the sake of comparison, since in our simulation context, we have access to the observations from the error process.  
    \item {\bf Two-step iid} and {\bf Two-step Whittle}: Fully data-driven penalization procedures, presented in Section~\ref{subsec:FullyDataDriven} below.
\end{itemize}


\subsection{Validation of the penalty when knowing the Hurst coefficient} 

In this section we validate the form of the penalties introduced  in Sections \ref{sec3} and \ref{sec:main_longrange2} above in various settings where the data-generating process is known, particularly when the Hurst coefficients of the error process are available.
We first examine error processes exhibiting short-range dependence, and subsequently investigate long-range dependence under several generative models for the design process.

\subsubsection{Short-range dependent errors}

This section illustrates the results of Section~\ref{sec:main_shortrange}. We consider short-range dependence in the error process while examining several design scenarios, all of which lead to penalties proportional to the model dimension. Note that these experiments extend those of \cite{caron2021gaussian}, which only considered deterministic designs, to settings involving random designs.

 \medskip 

\begin{figure}[htbp]
\centering
\begin{subfigure}{0.48\textwidth}
\centering
\includegraphics[width=\linewidth]{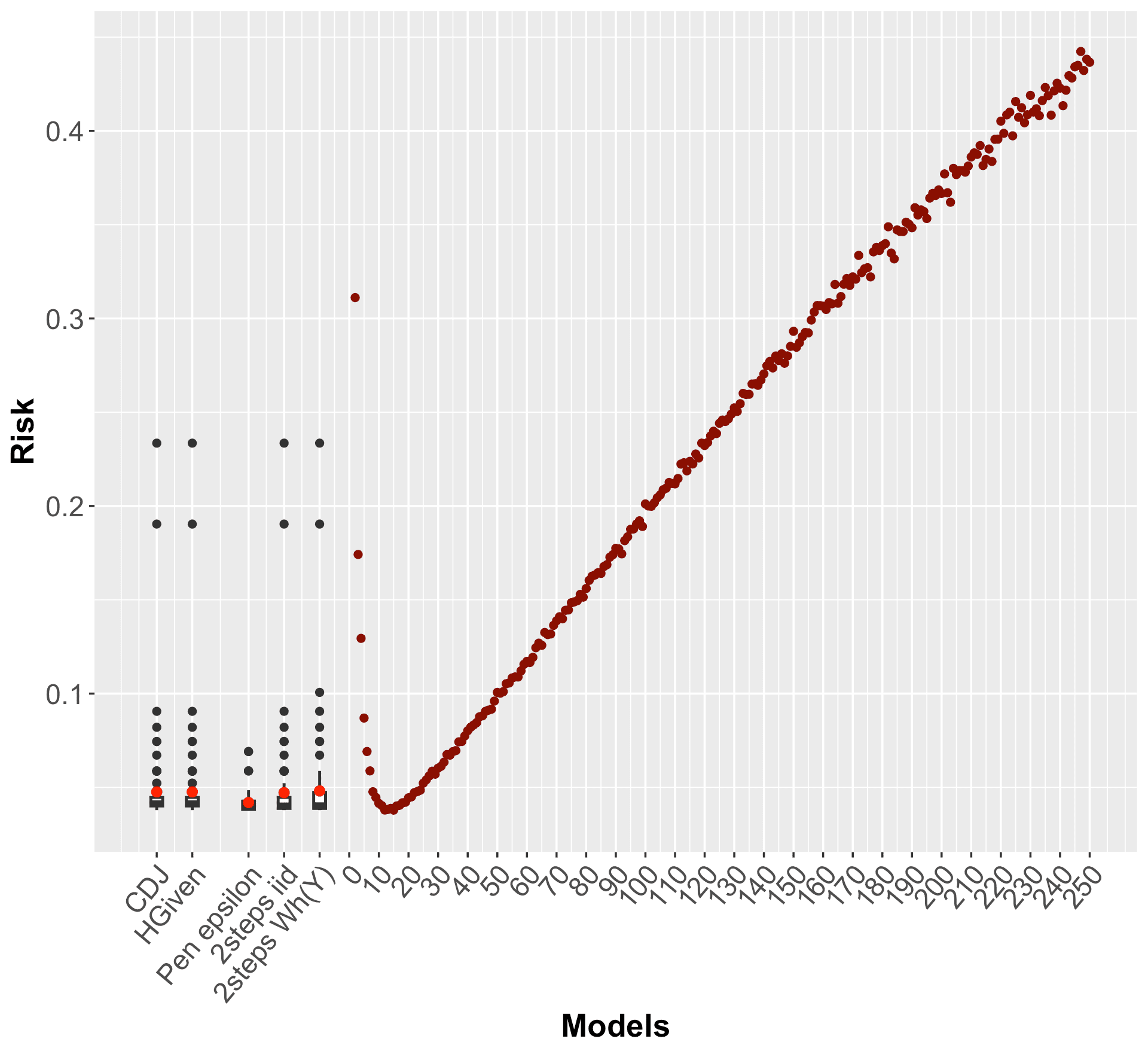}
\caption{$f_1$, $n=500$}
\end{subfigure}
\hfill
\begin{subfigure}{0.48\textwidth}
\centering
\includegraphics[width=\linewidth]{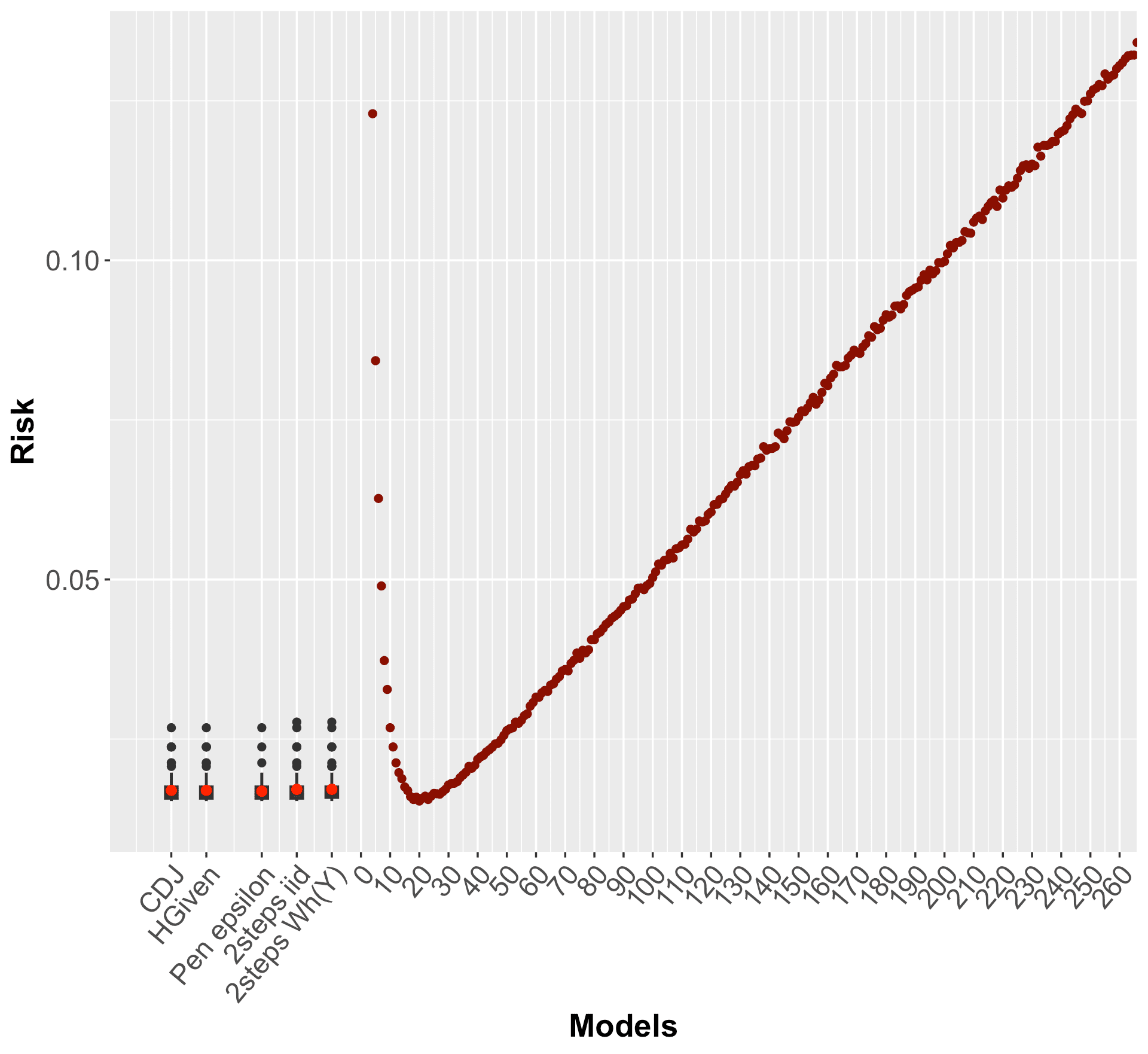}
\caption{$f_1$, $n=2000$}
\end{subfigure}

\vspace{0.3cm}

\begin{subfigure}{0.48\textwidth}
\centering
\includegraphics[width=\linewidth]{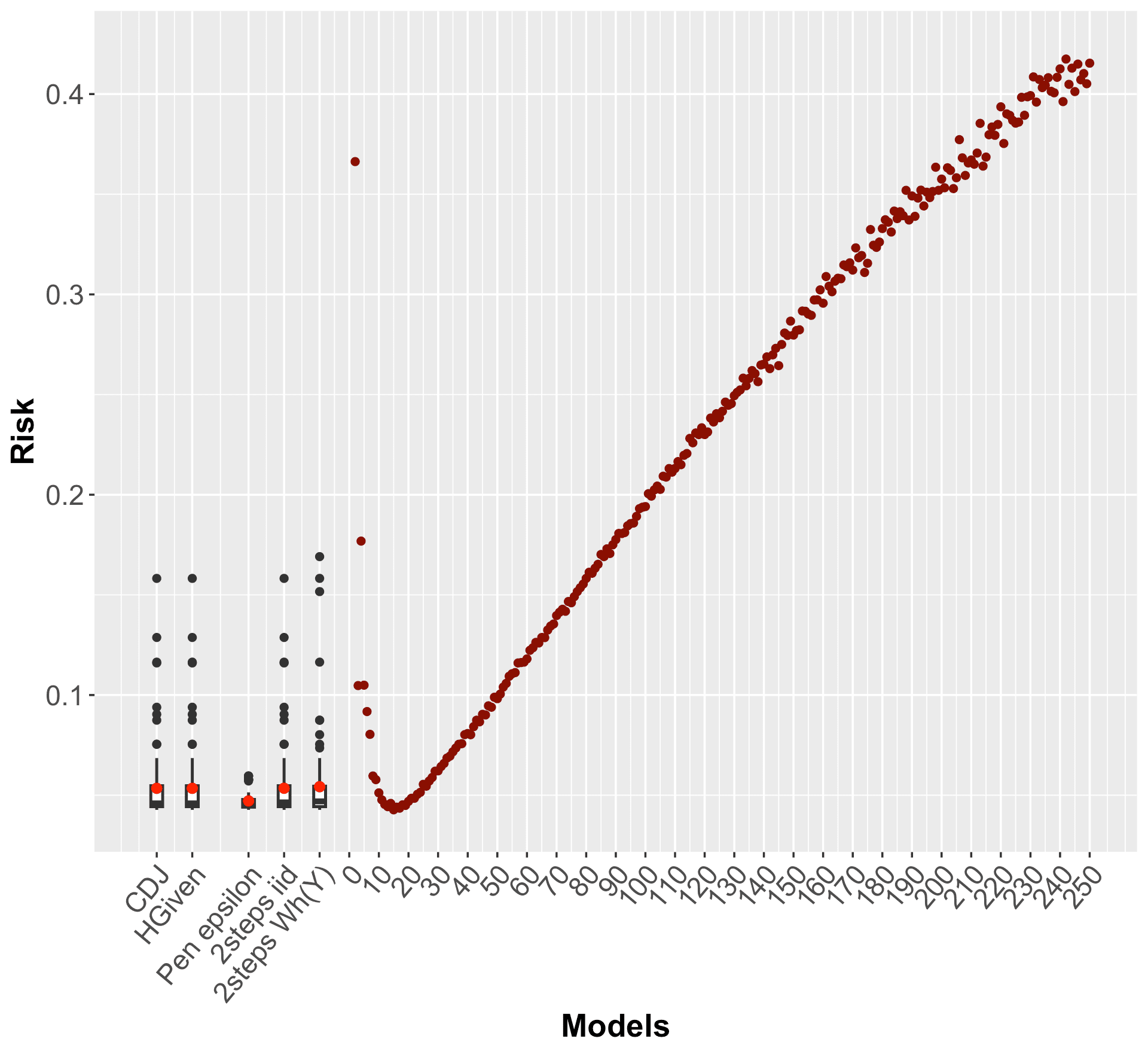}
\caption{$f_2$, $n=500$}
\end{subfigure}
\hfill
\begin{subfigure}{0.48\textwidth}
\centering
\includegraphics[width=\linewidth]{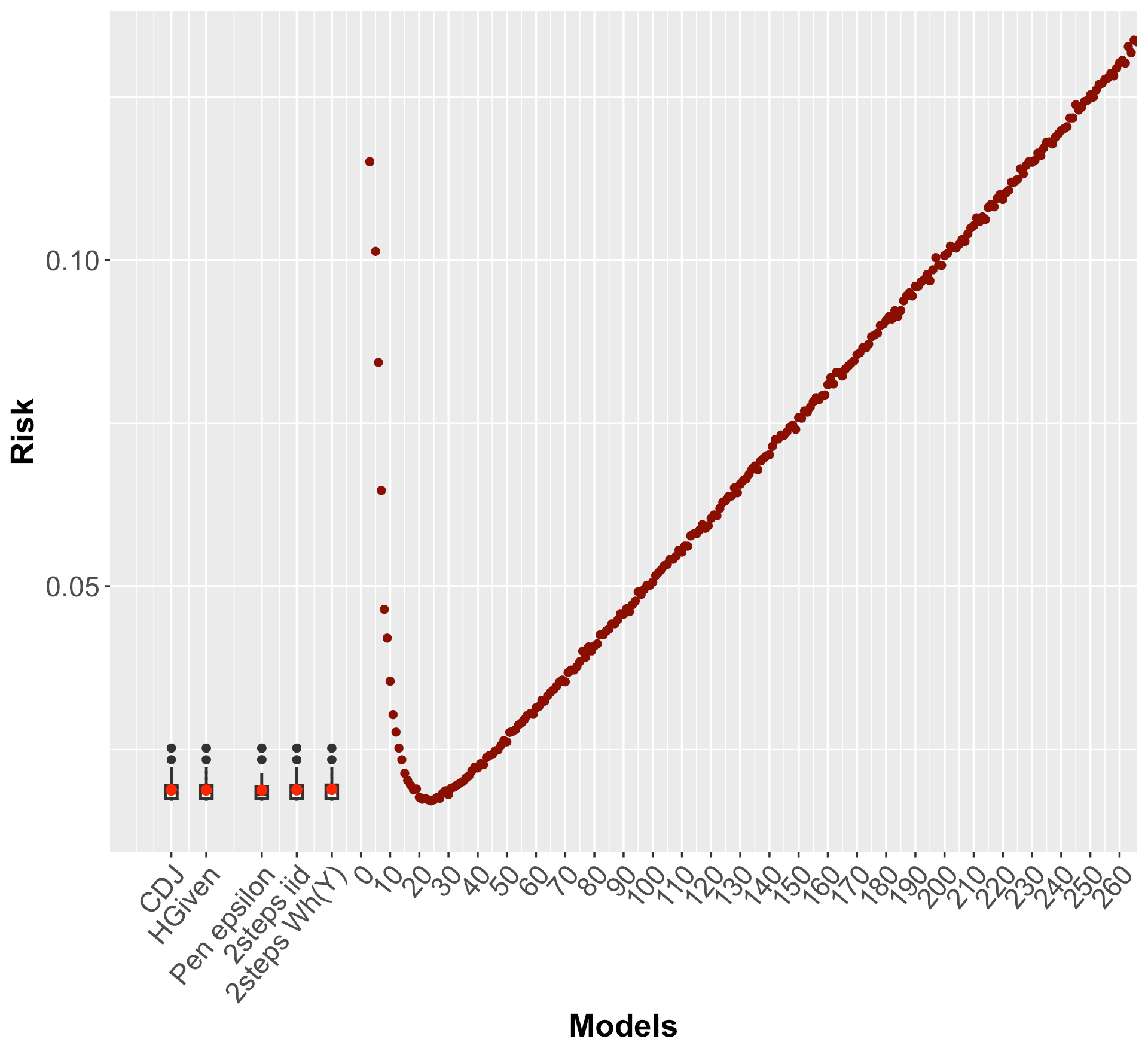}
\caption{$f_2$, $n=2000$}
\end{subfigure}

\caption{Risk performance of the  procedures for the two regression functions and two sample sizes in Experiment~1. The error process is  i.i.d. Gaussian, and the design  is i.i.d. with uniform distribution over $[0,1]$.}

\label{fig:short_range_dependence_error_rnorm_design_unif}
\end{figure}

$\bullet$ {\bf Experiment 1}: the error process is Gaussian i.i.d., and the design is i.i.d. with uniform distribution over $[0,1]$, see Fig.~\ref{fig:short_range_dependence_error_rnorm_design_unif}.

 \medskip 

 \medskip 

\begin{figure}[htbp]
\centering
\begin{subfigure}{0.48\textwidth}
\centering
\includegraphics[width=\linewidth]{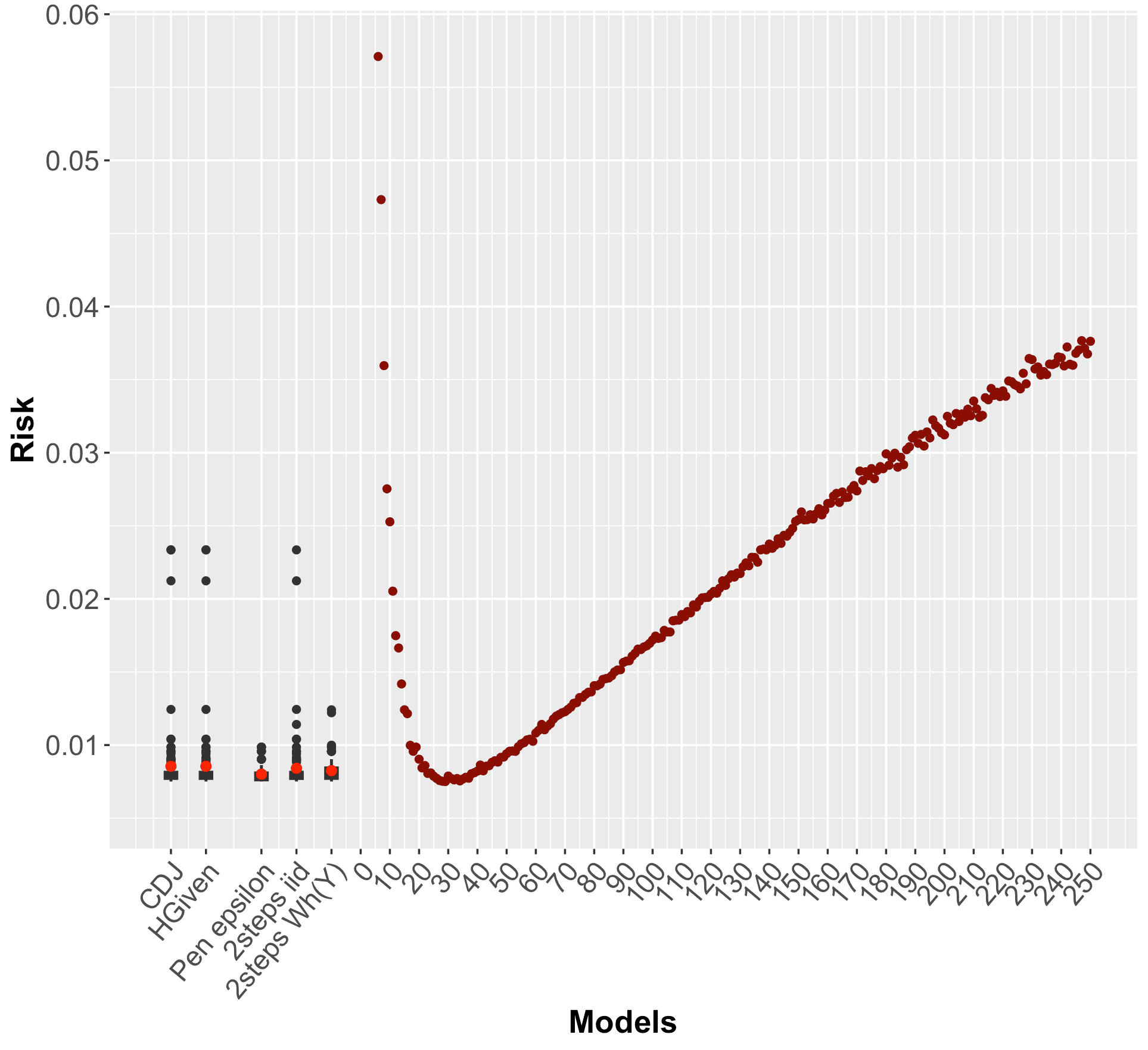}
\caption{$f_1$, $n=500$}
\end{subfigure}
\hfill
\begin{subfigure}{0.48\textwidth}
\centering
\includegraphics[width=\linewidth]{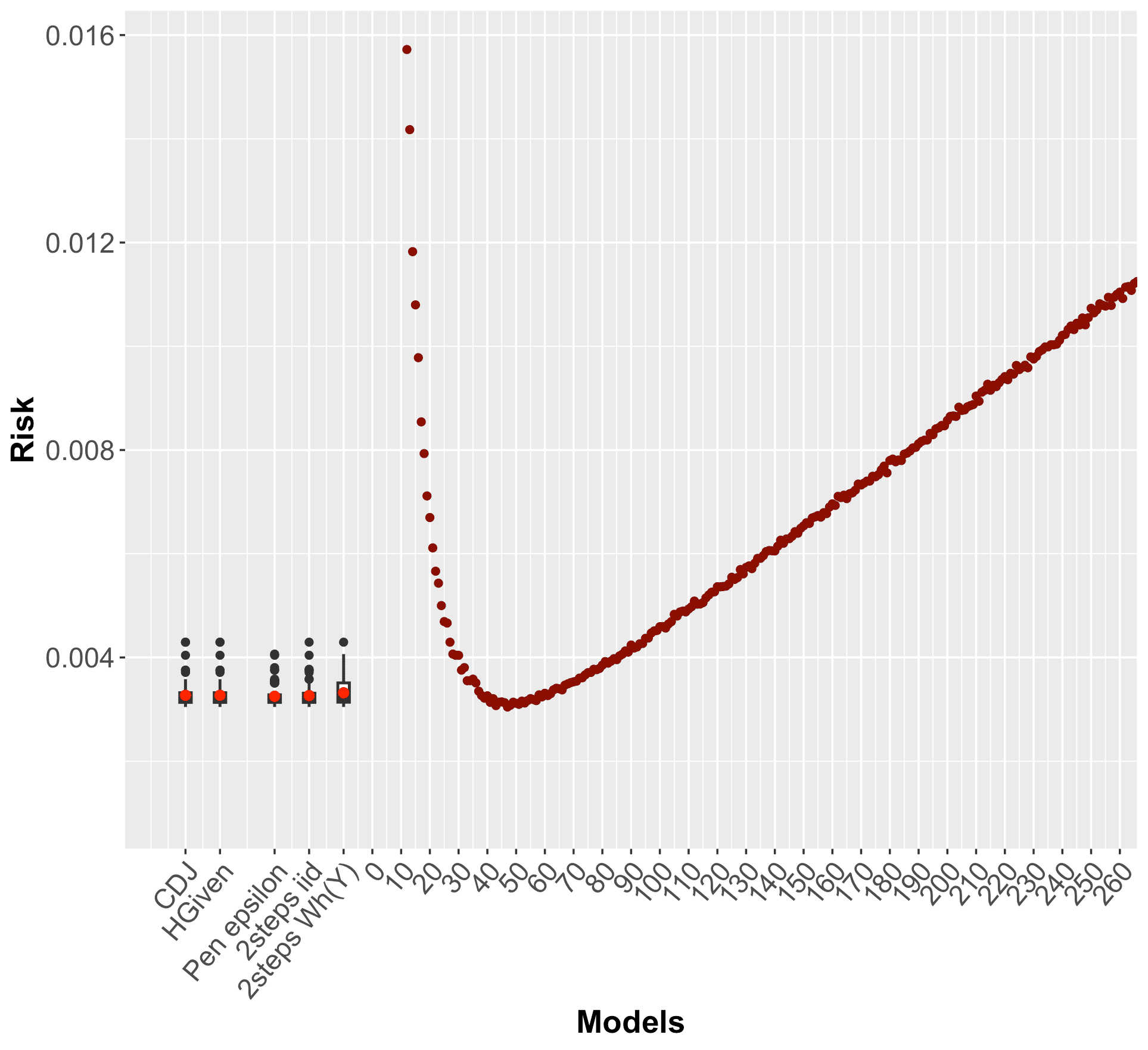}
\caption{$f_1$, $n=2000$}
\end{subfigure}

\vspace{0.3cm}

\begin{subfigure}{0.48\textwidth}
\centering
\includegraphics[width=\linewidth]{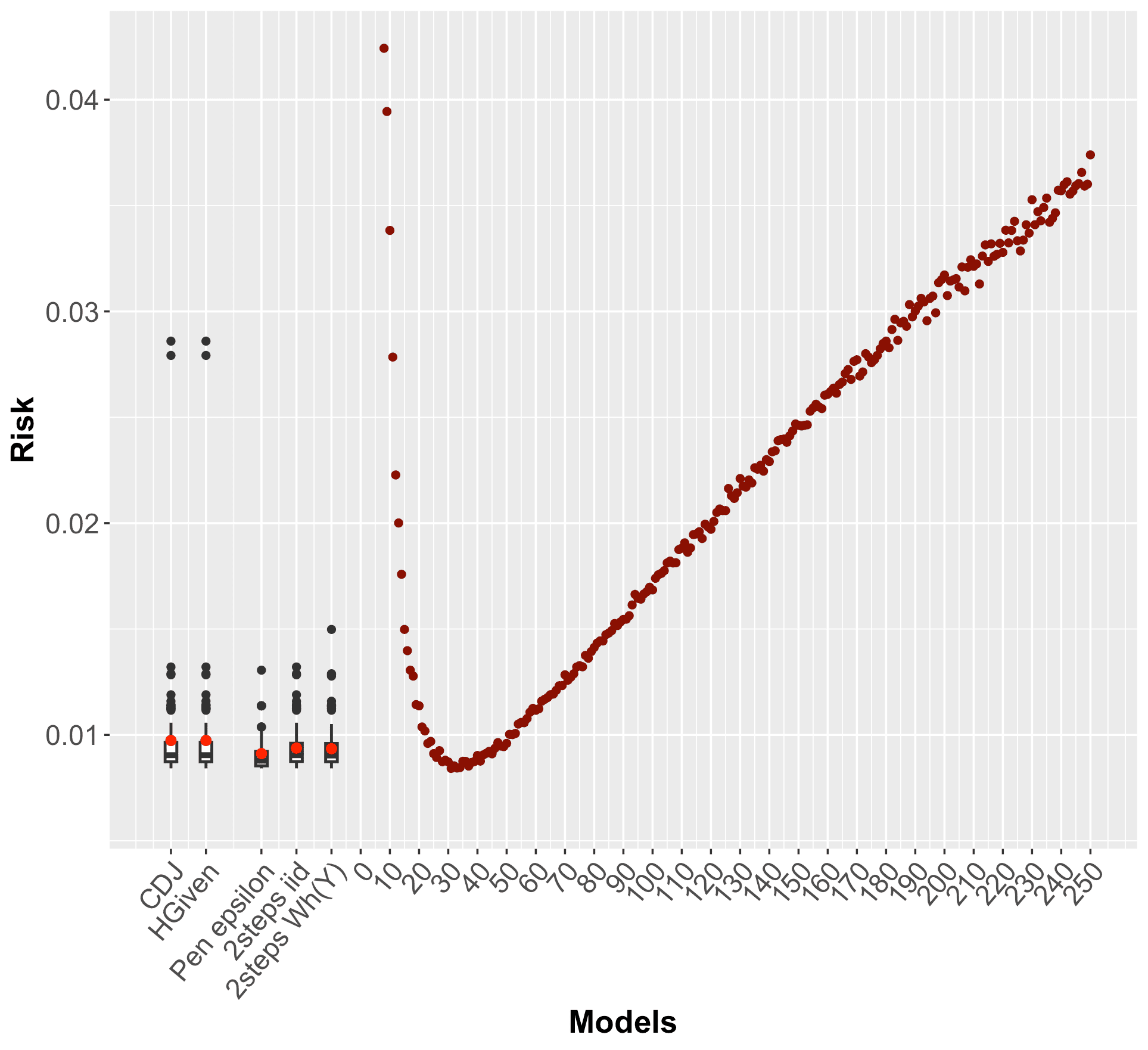}
\caption{$f_2$, $n=500$}
\end{subfigure}
\hfill
\begin{subfigure}{0.48\textwidth}
\centering
\includegraphics[width=\linewidth]{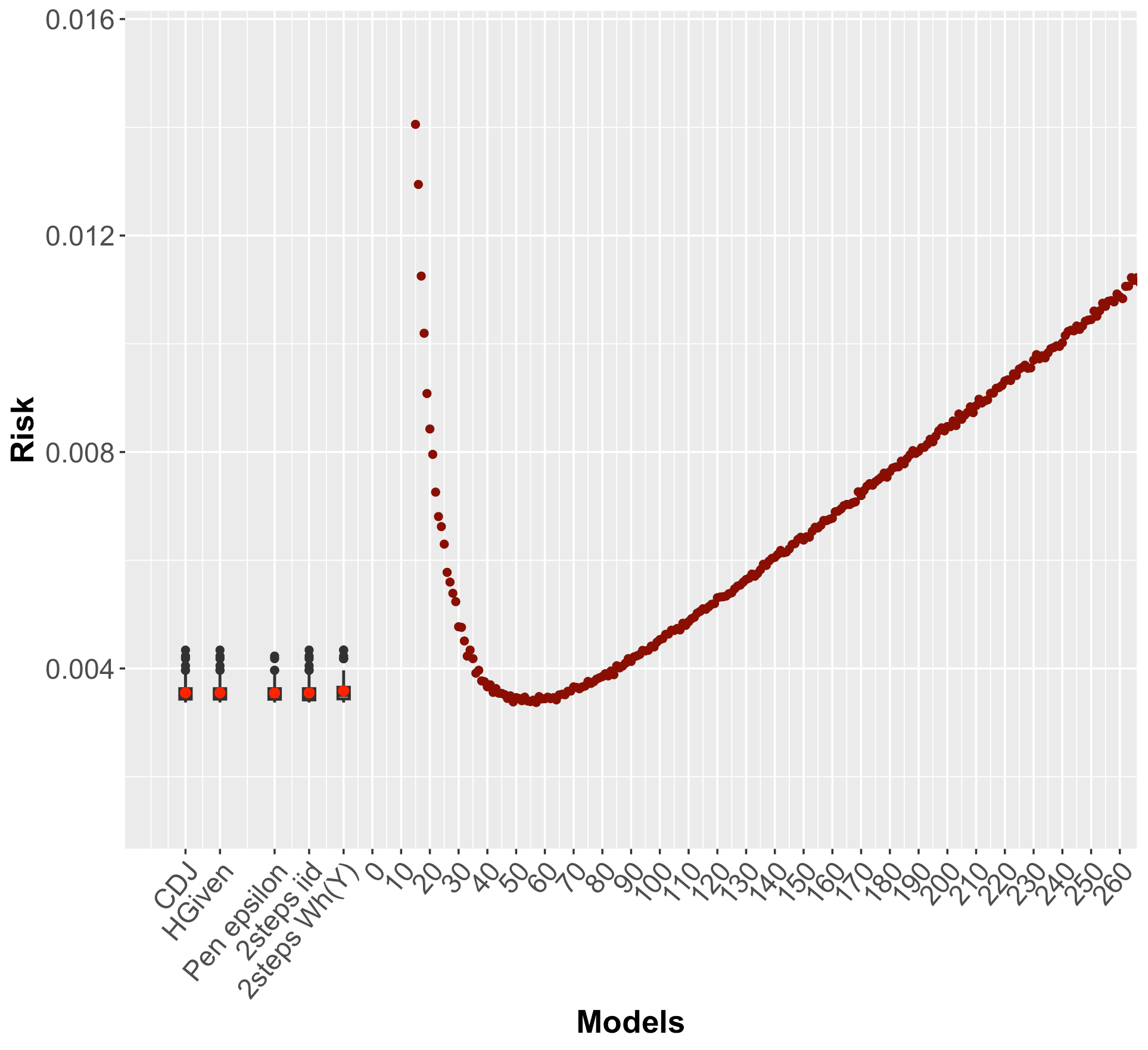}
\caption{$f_2$, $n=2000$}
\end{subfigure}
\caption{Risk performance of the penalization procedures for the two regression functions and two sample sizes in Experiment~2. The error process is the non mixing centered AR(1) process \eqref{Nonmix_ar1},  and the design  is i.i.d. with uniform distribution over $[0,1]$.}
\label{fig:short_range_dependence_error_nonmixar1_design_unif}
\end{figure}

 $\bullet$ {\bf Experiment 2}: the error process is the centered non mixing AR(1) process~\eqref{Nonmix_ar1}, and the design  is i.i.d. with uniform distribution over $[0,1]$, see Fig.~\ref{fig:short_range_dependence_error_nonmixar1_design_unif}.

 \medskip 

 \medskip  
 
\begin{figure}[htbp]
\centering
\begin{subfigure}{0.48\textwidth}
\centering
\includegraphics[width=\linewidth]{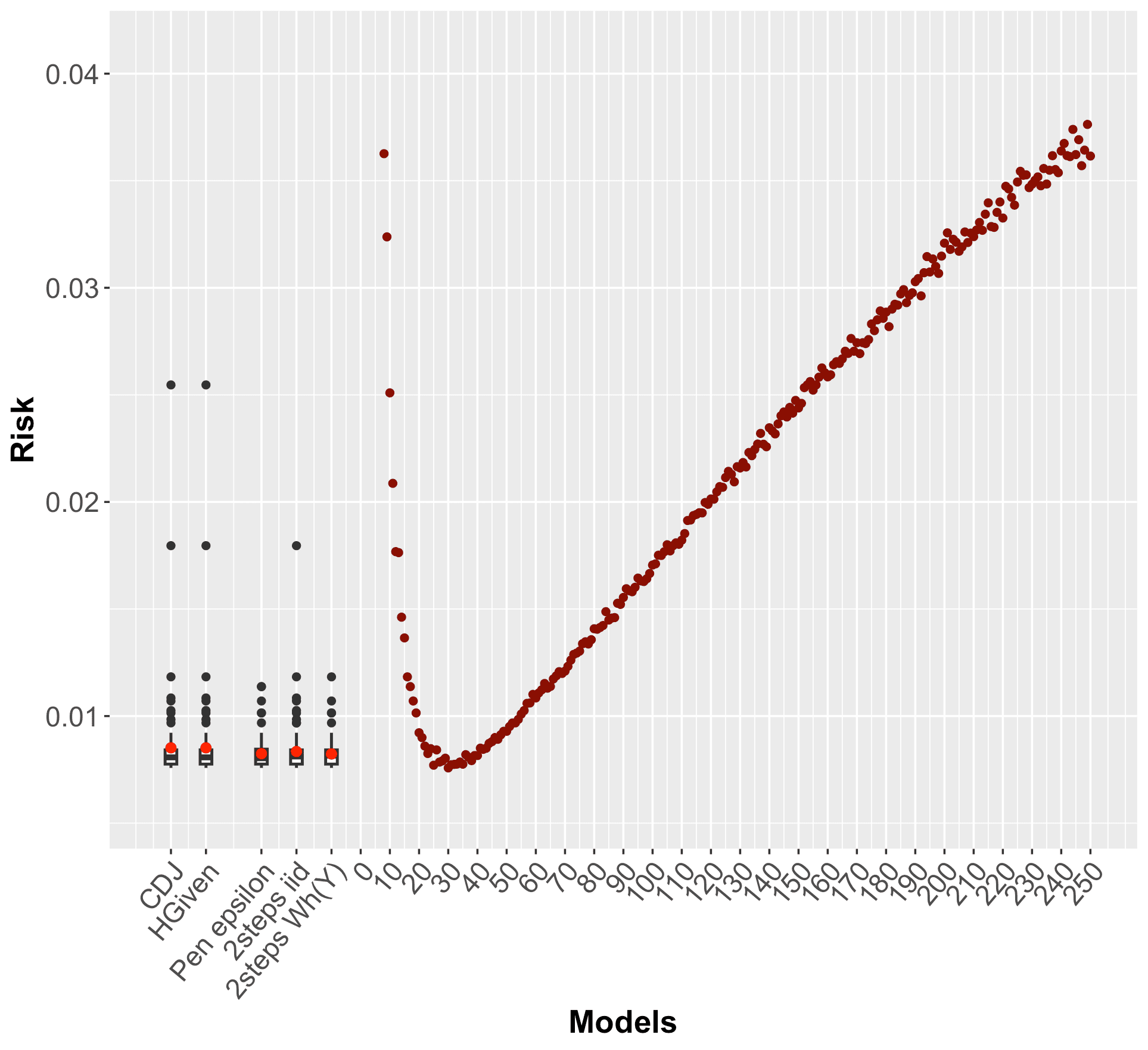}
\caption{$f_1$, $n=500$}
\end{subfigure}
\hfill
\begin{subfigure}{0.48\textwidth}
\centering
\includegraphics[width=\linewidth]{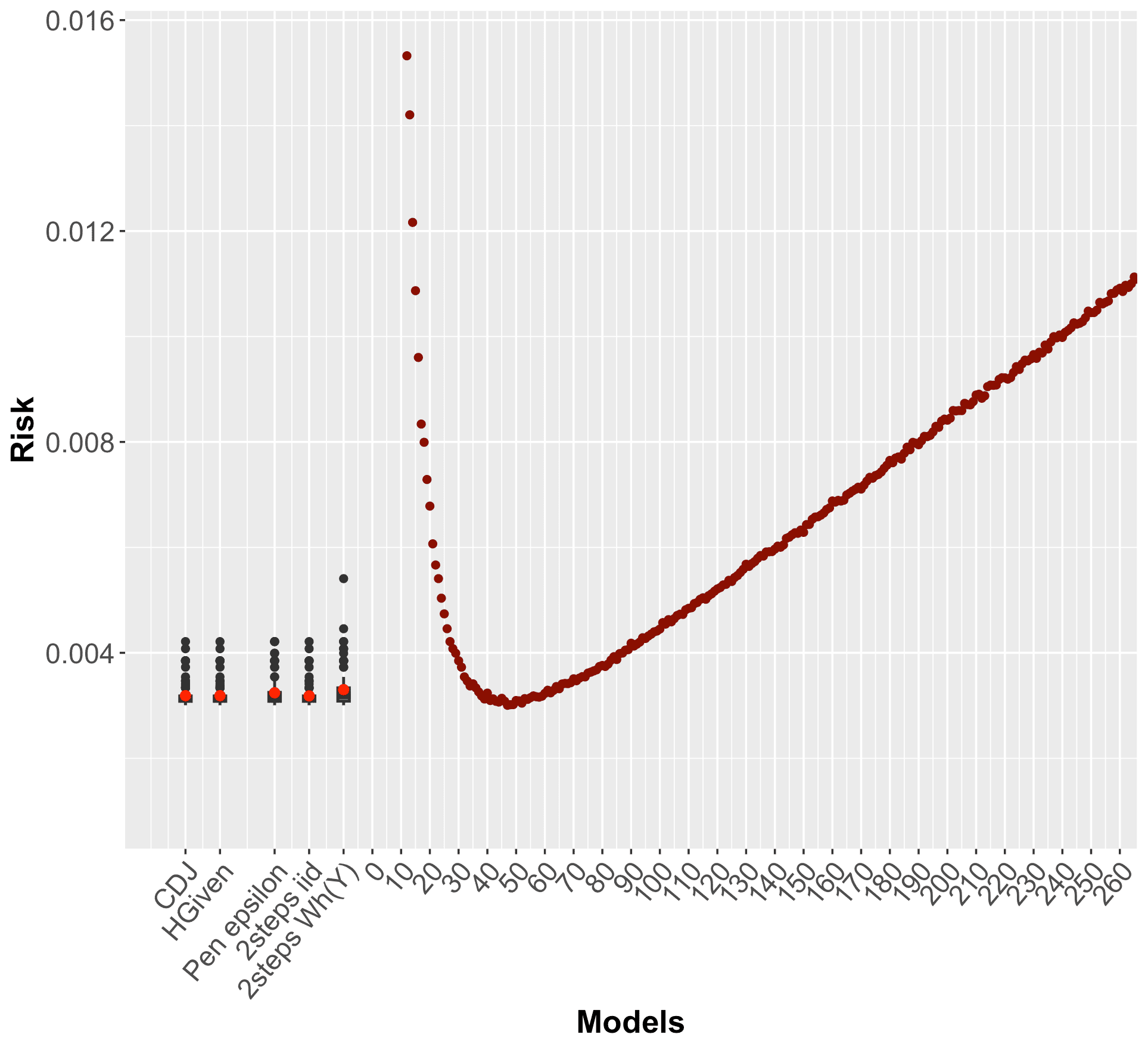}
\caption{$f_1$, $n=2000$}
\end{subfigure}

\vspace{0.3cm}

\begin{subfigure}{0.48\textwidth}
\centering
\includegraphics[width=\linewidth]{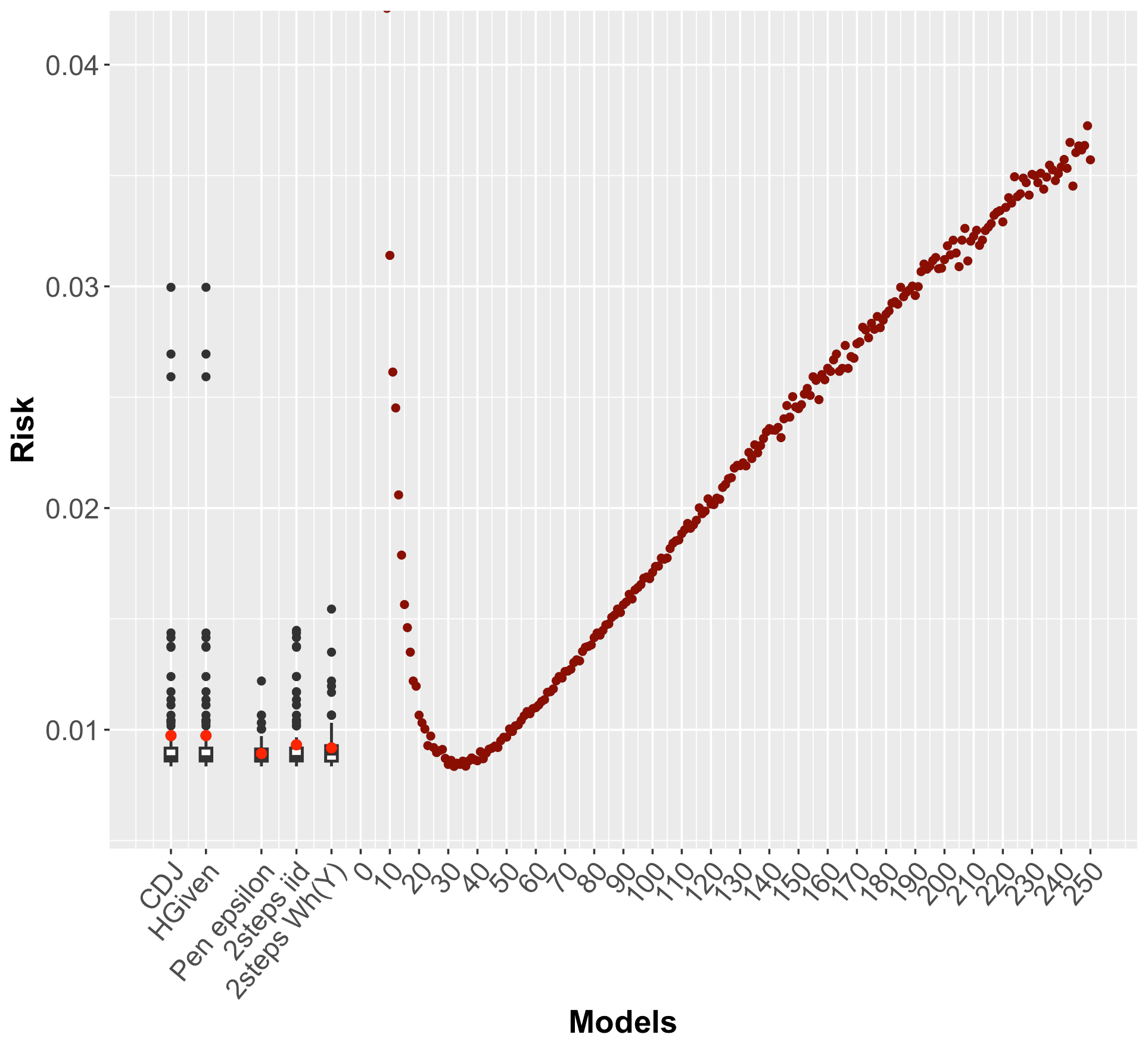}
\caption{$f_2$, $n=500$}
\end{subfigure}
\hfill
\begin{subfigure}{0.48\textwidth}
\centering
\includegraphics[width=\linewidth]{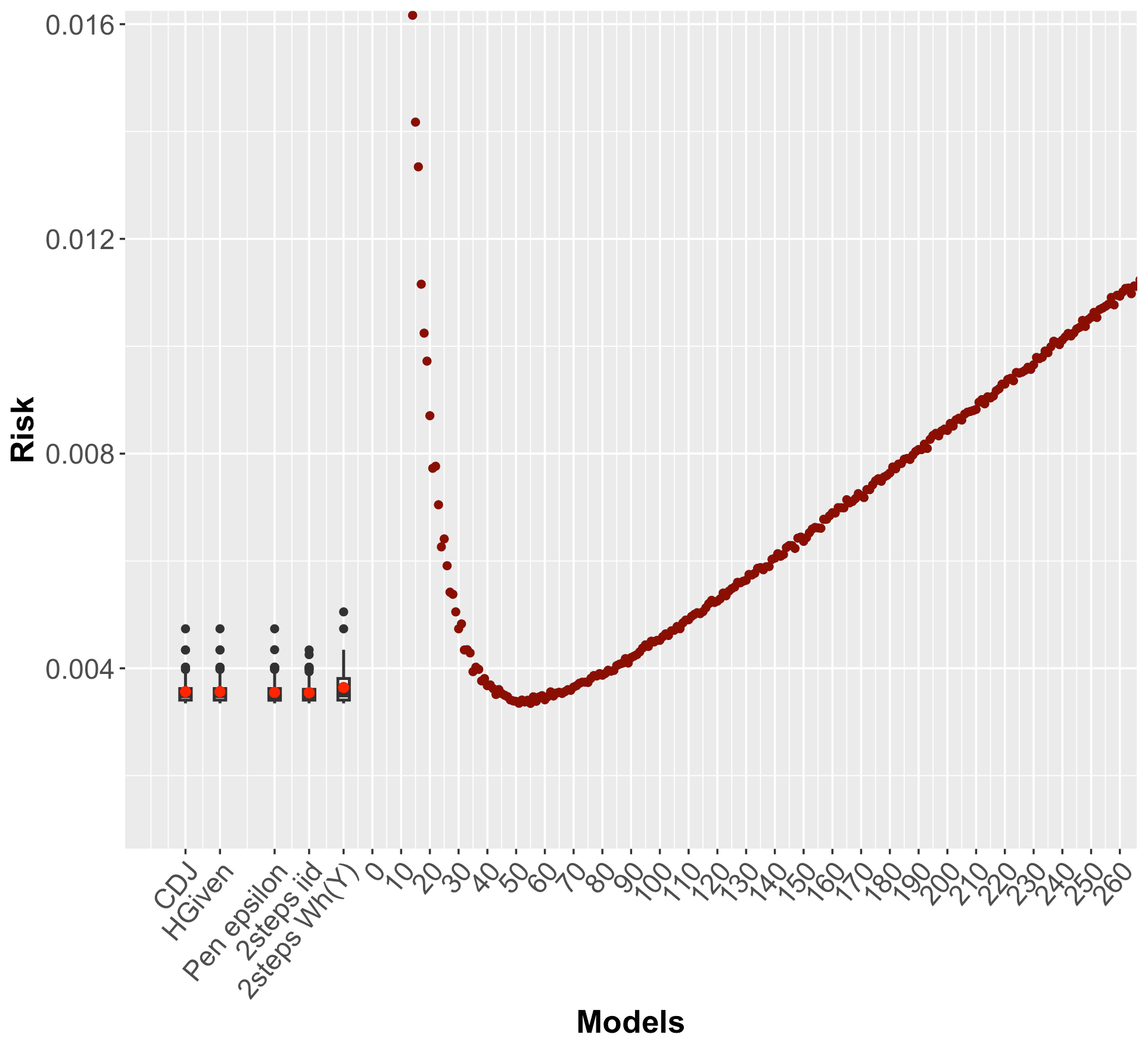}
\caption{$f_2$, $n=2000$}
\end{subfigure}
\caption{Risk performance of the penalization procedures for the two regression functions and two sample sizes in Experiment~3. The error process is  the non mixing centered AR(1) process \eqref{Nonmix_ar1},  and the design is obtained from  a FGN  with $H=0.7$ via the transformation \eqref{transFGN}.}
\label{fig:short_range_dependence_error_nonmixar1_design_FGN07}
\end{figure}

$\bullet$  {\bf Experiment 3}: the error process is the centered non mixing AR(1) process~\eqref{Nonmix_ar1}, and the design  is obtained from  a FGN  with $H=0.7$ via the transformation \eqref{transFGN}, see Fig.~\ref{fig:short_range_dependence_error_nonmixar1_design_FGN07}.

 \medskip 

 \medskip 

$\bullet$  {\bf Experiment 4}: the error process is a FGN  with $H=0.4$,  and the design  is obtained from  a FGN  with $H=0.7$ via the transformation \eqref{transFGN}, see Fig.~\ref{fig:short_range_dependence_error_nonmixar1_design_FGN07}.

 \medskip 

\begin{figure}[htbp]
\centering
\begin{subfigure}{0.48\textwidth}
\centering
\includegraphics[width=\linewidth]{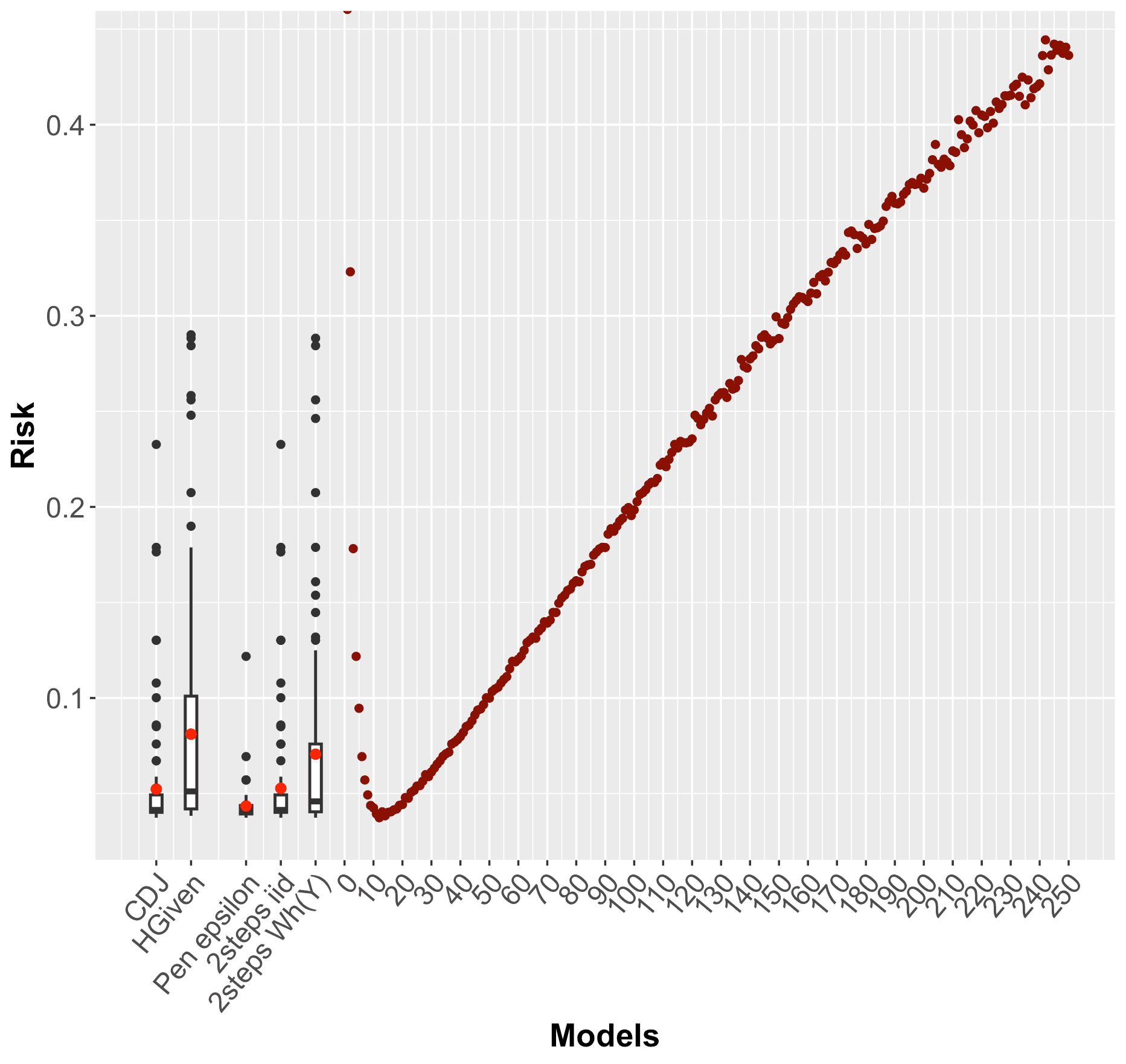}
\caption{$f_1$, $n=500$}
\end{subfigure}
\hfill
\begin{subfigure}{0.48\textwidth}
\centering
\includegraphics[width=\linewidth]{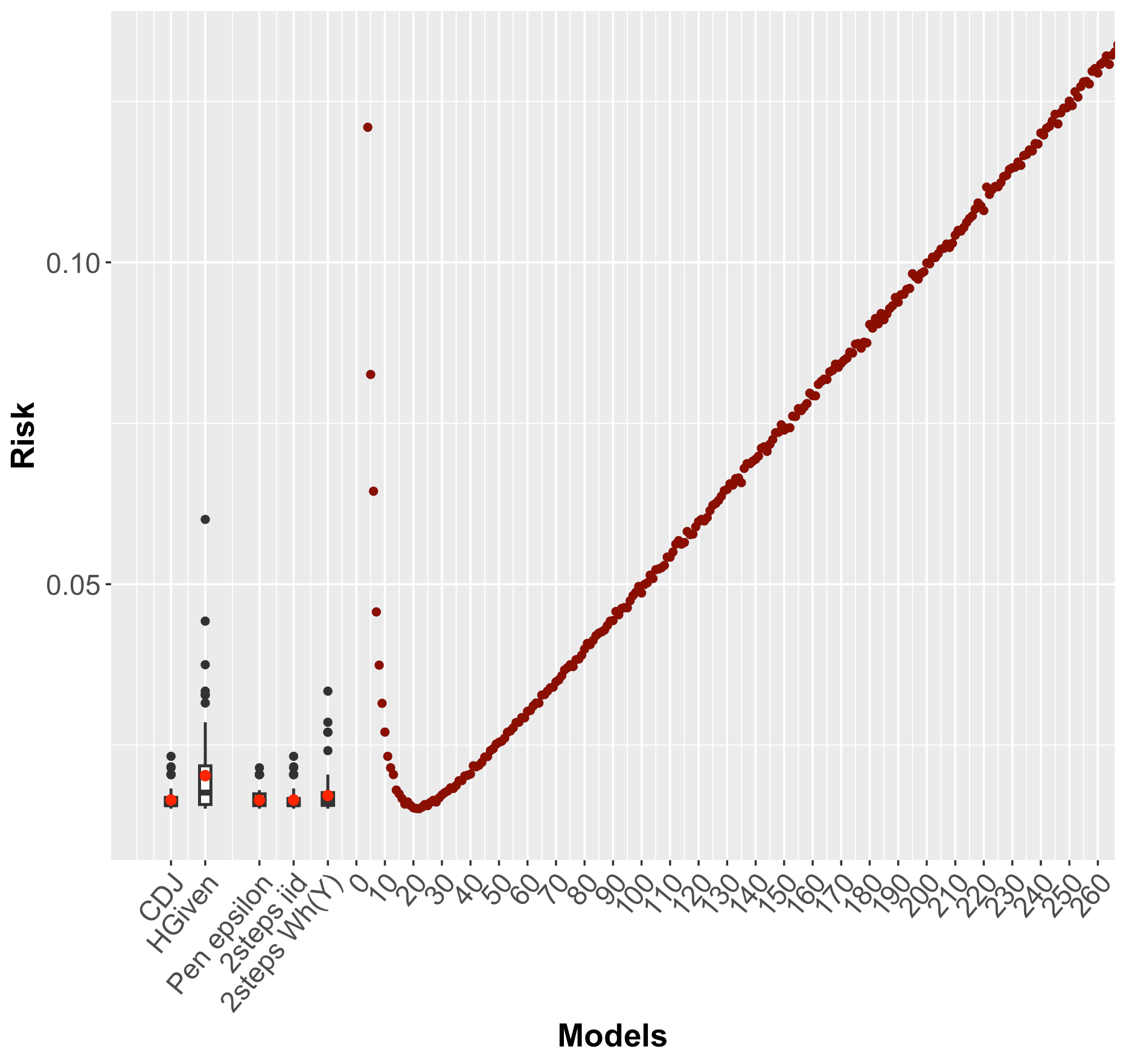}
\caption{$f_1$, $n=2000$}
\end{subfigure}

\vspace{0.3cm}

\begin{subfigure}{0.48\textwidth}
\centering
\includegraphics[width=\linewidth]{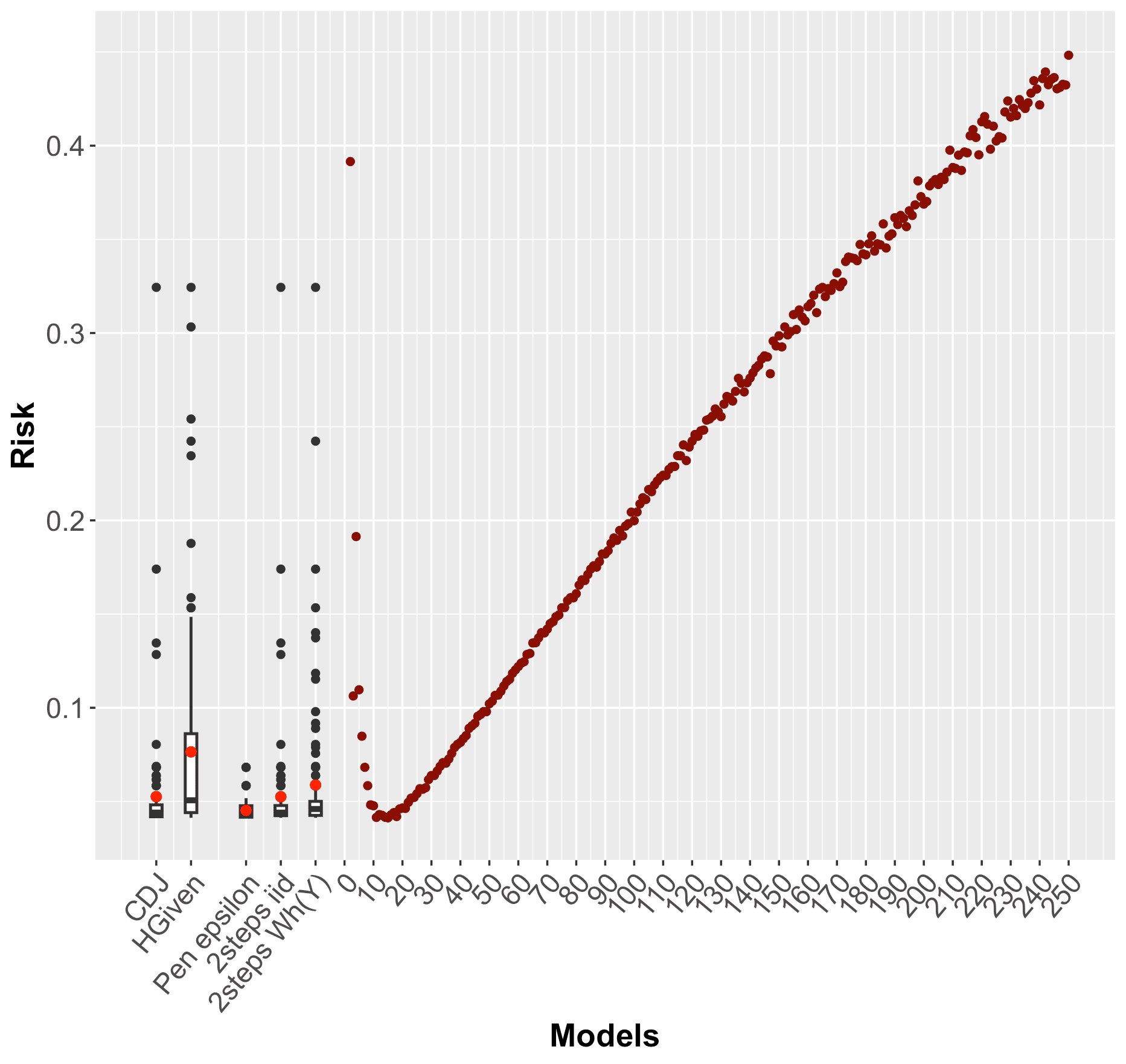}
\caption{$f_2$, $n=500$}
\end{subfigure}
\hfill
\begin{subfigure}{0.48\textwidth}
\centering
\includegraphics[width=\linewidth]{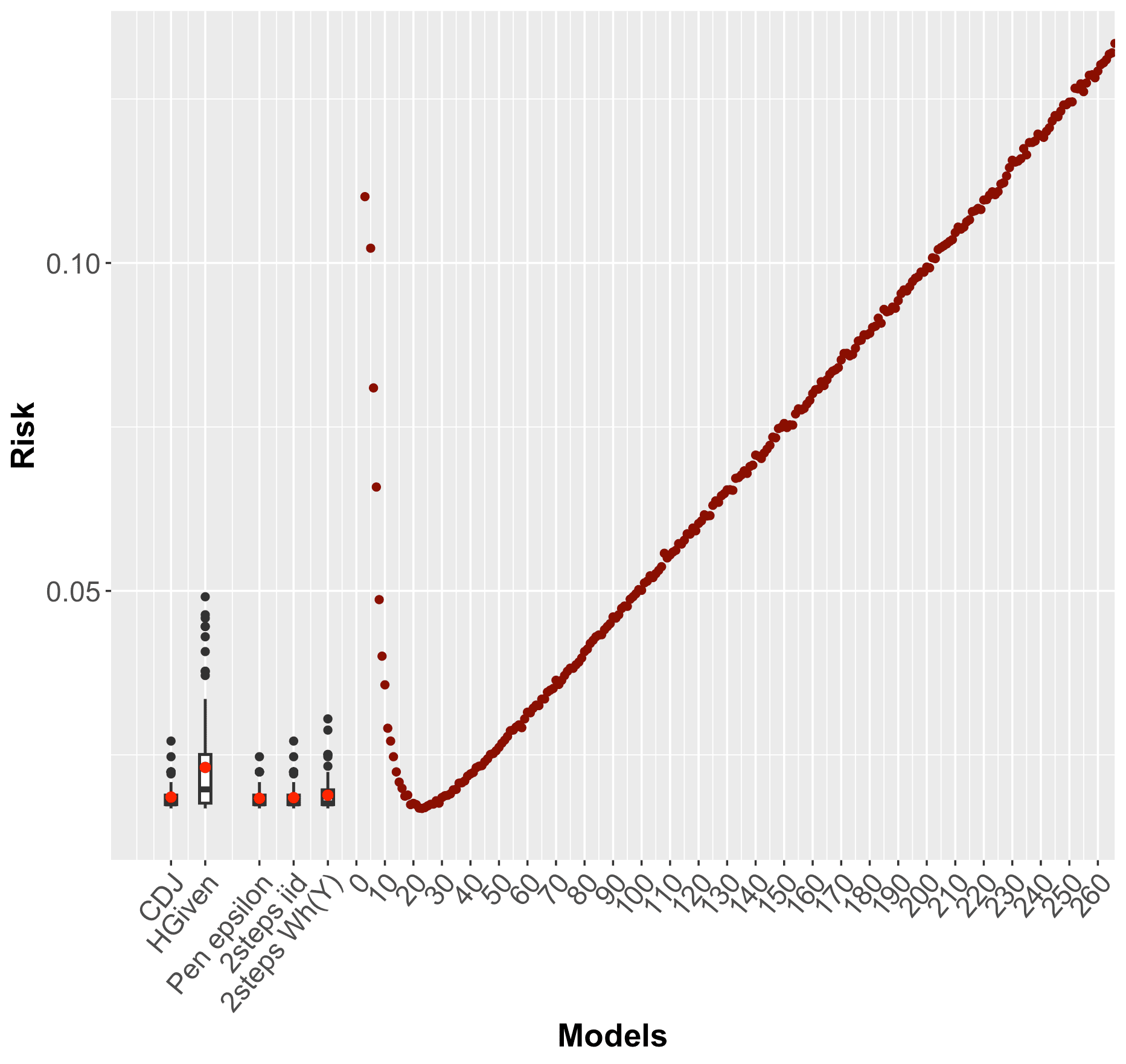}
\caption{$f_2$, $n=2000$}
\end{subfigure}
\caption{Risk performance of the penalization procedures for the two regression functions and two sample sizes in Experiment~4. The error process is a FGN  with $H=0.4$,  and the design  is obtained from  a FGN  with $H=0.7$ via the transformation \eqref{transFGN}.}
\label{fig:short_range_dependence_error_FGN04_design_FGN07}
\end{figure}

\medskip 

For these four error and design process configurations, the dimension-proportional penalty CDJ yields good results, in line  with the theoretical  results of Section~\ref{sec:main_shortrange}. Clearly, for the three first experiments, the penalty with $H = 0.5$ (HGiven procedure) gives exactly the same results. This is not the case for Experiment 4, for which the penalty with $H = 0.4$ performs poorly. This observation is  again consistent with the results of Section~\ref{sec:main_shortrange}.

The performances of the other two penalties will be discussed later in section~\ref{subsec:FullyDataDriven}. It should be noted that, so far, the error processes have been simulated independently of the design processes. 

\medskip

We now turn our attention to a heteroscedastic setting.

 \medskip 
 
\begin{figure}[htbp]
\centering
\begin{subfigure}{0.48\textwidth}
\centering
\includegraphics[width=\linewidth]{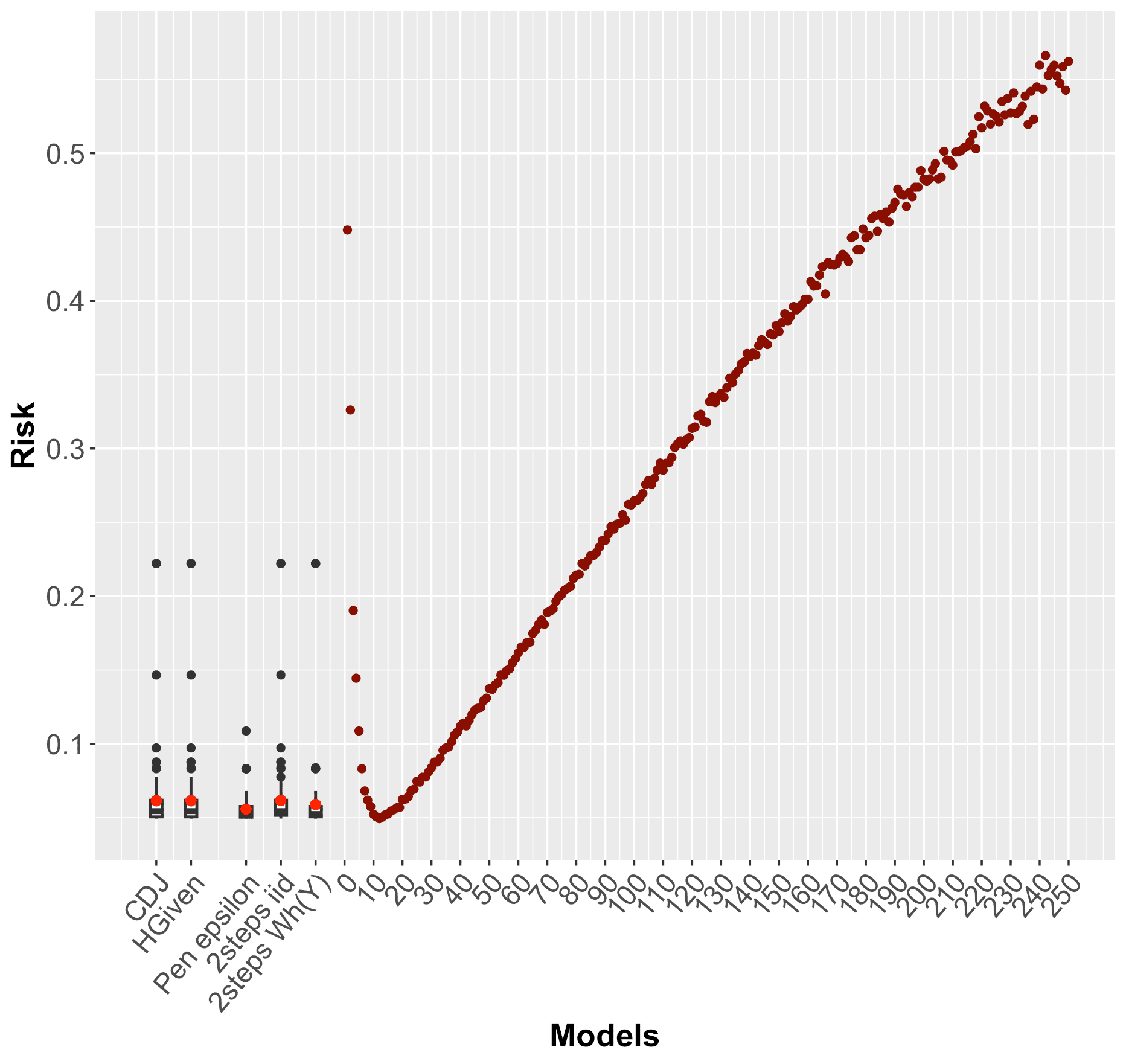}
\caption{$f_1$, $n=500$}
\end{subfigure}
\hfill
\begin{subfigure}{0.48\textwidth}
\centering
\includegraphics[width=\linewidth]{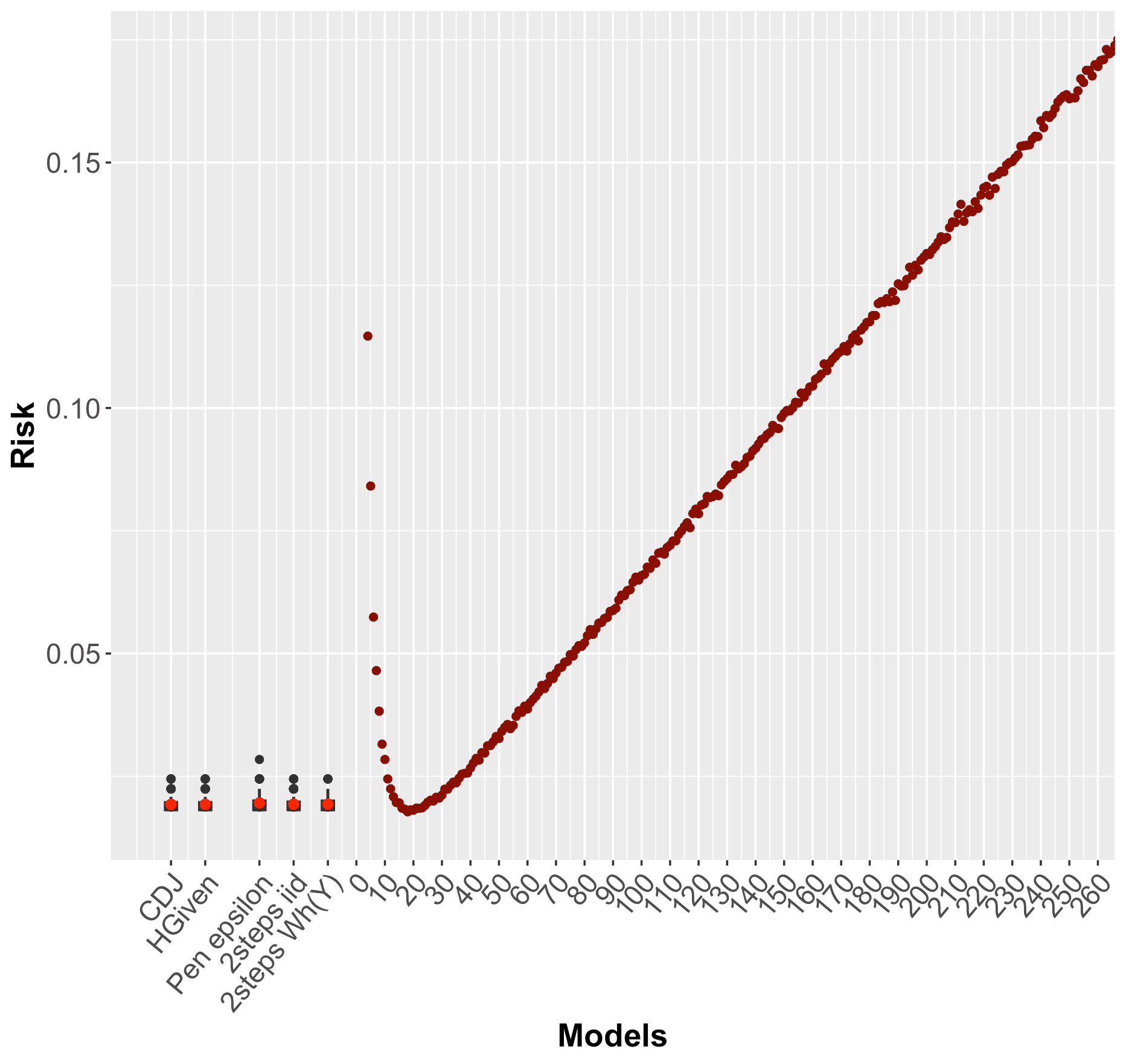}
\caption{$f_1$, $n=2000$}
\end{subfigure}

\vspace{0.3cm}

\begin{subfigure}{0.48\textwidth}
\centering
\includegraphics[width=\linewidth]{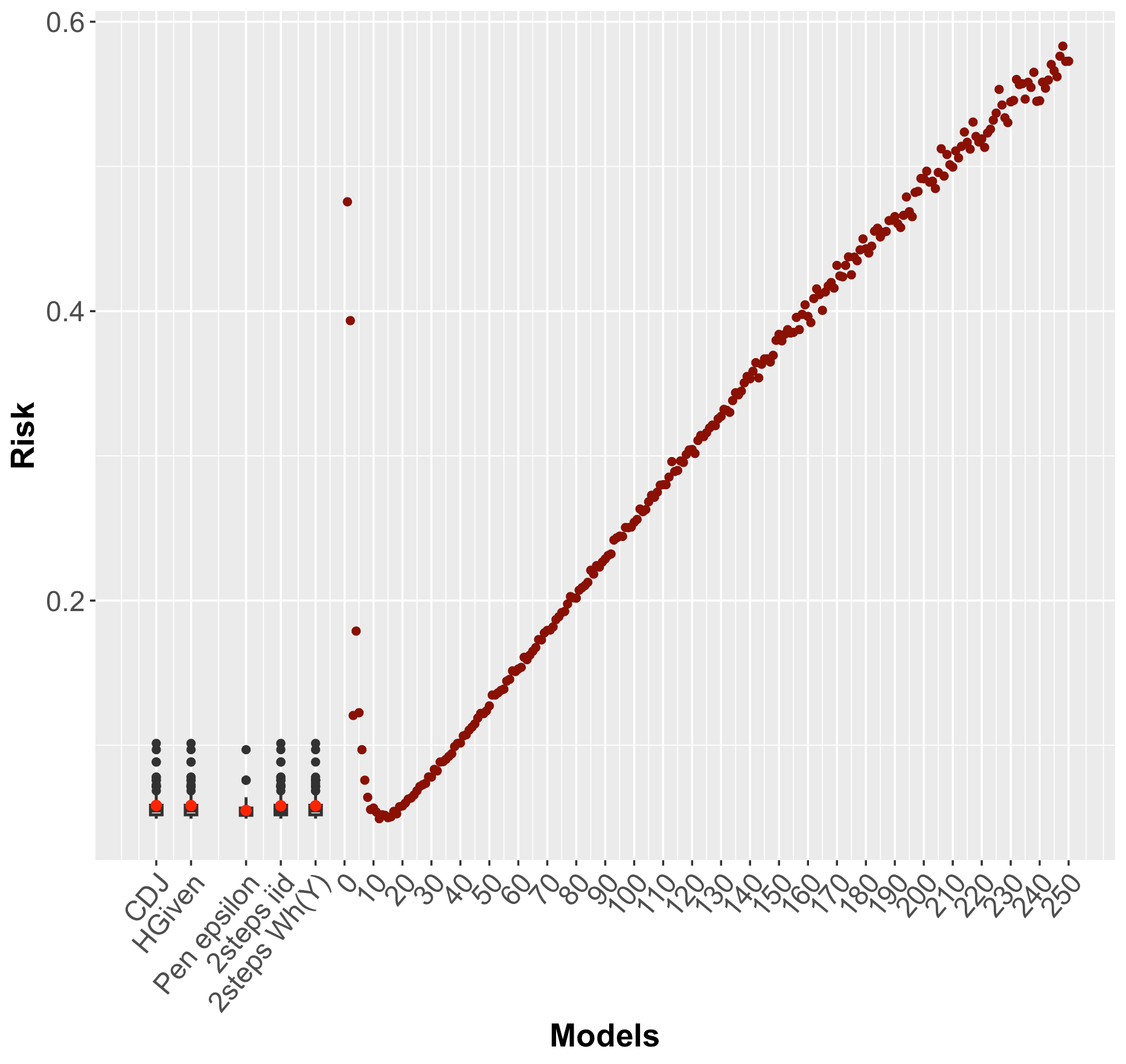}
\caption{$f_2$, $n=500$}
\end{subfigure}
\hfill
\begin{subfigure}{0.48\textwidth}
\centering
\includegraphics[width=\linewidth]{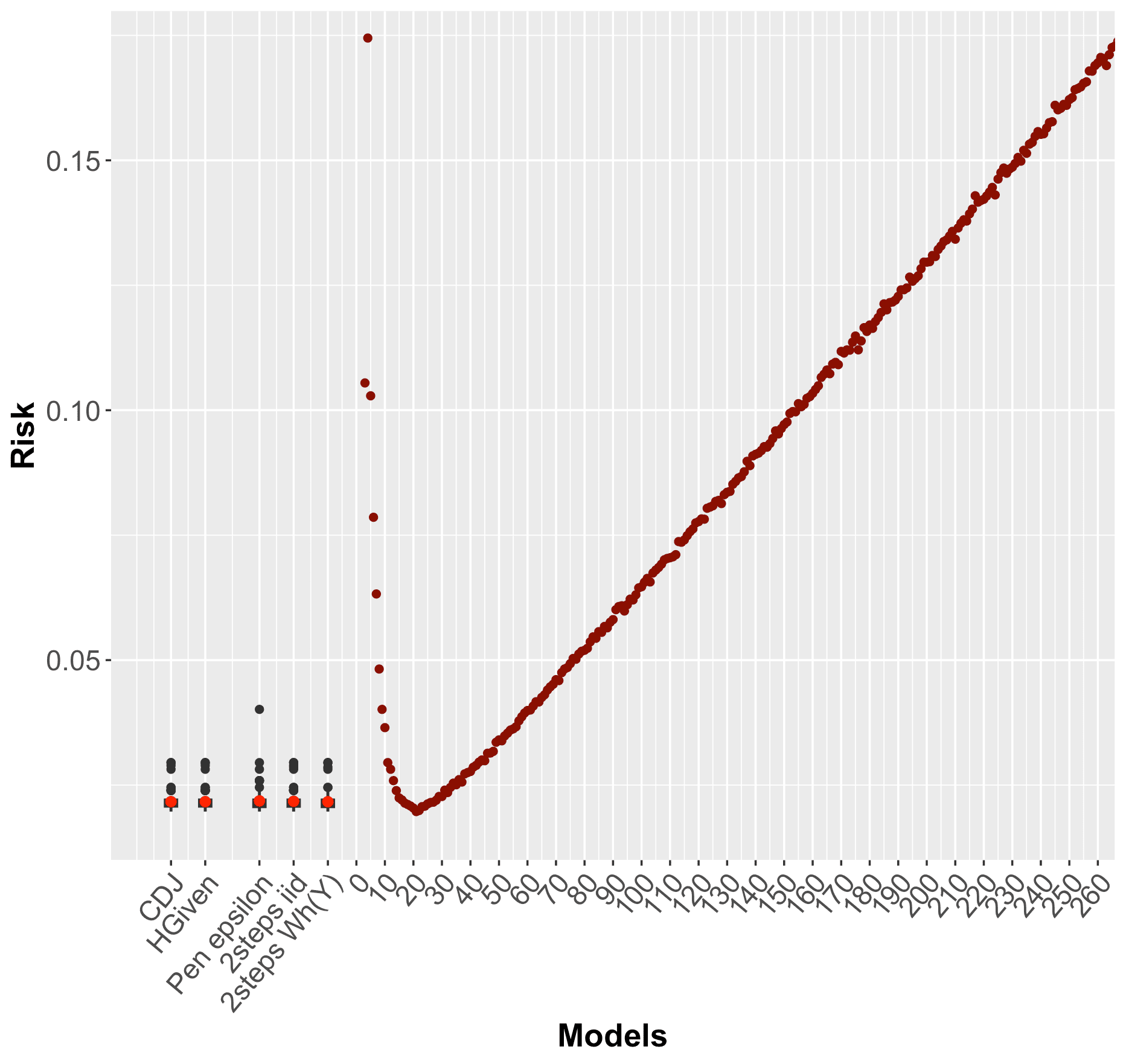}
\caption{$f_2$, $n=2000$}
\end{subfigure}
\caption{Risk performance of the penalization procedures for the two regression functions and two sample sizes in Experiment~5. The error process is an heteroscedastic FGN  with $H=0.5$ (see \eqref{Heteps}),  and the design is obtained from  a FGN  with $H=0.7$ via the transformation \eqref{transFGN}.}
\label{fig:short_range_dependence_error_FGN05hetero_design_FGN07}
\end{figure}

$\bullet$  {\bf Experiment 5}:
The design  is obtained from  a FGN  with $H=0.7$ via the transformation \eqref{transFGN}.
see Fig.~\ref{fig:short_range_dependence_error_FGN05hetero_design_FGN07}.
The heteroscedastic error process is built as follows. Let $(e_i)_{1\le i\le n}$ be a centered FGN process with Hurst parameter $H=0.5$, independent of the design process. We define
\begin{equation}\label{Heteps}
\varepsilon_i = g_i({\mathbf X}_n) e_i
\end{equation}
where, for $2 \le i \le n-1$,
$$
g_i({\mathbf X}_n) = \sqrt{1+0.2X_{i-1}^2+0.5X_i^2+0.2X_{i+1}^2},$$ 
with the usual adaptations at the boundaries. Consequently, we get
$$\operatorname{Var}(\varepsilon_i\mid {\mathbf X}_n) = g_i({\mathbf X}_n)^2.$$
The simulations  are presented in Fig.~\ref{fig:short_range_dependence_error_FGN05hetero_design_FGN07}.

 \medskip

The CDJ penalty (equivalently HGiven = 0.5) again yields good results, in line with the results given in Section~\ref{sec:main_shortrange}.

\subsubsection{Long range dependent errors}

We now consider several experiments with long range dependent errors, as in the context of Section~\ref{sec:main_longrange}. These experiments typically involve \(\beta\)-mixing (and thus \(\alpha\)-dependent) designs.
 
 \medskip 
 
\begin{figure}[htbp]
\centering
\begin{subfigure}{0.48\textwidth}
\centering
\includegraphics[width=\linewidth]{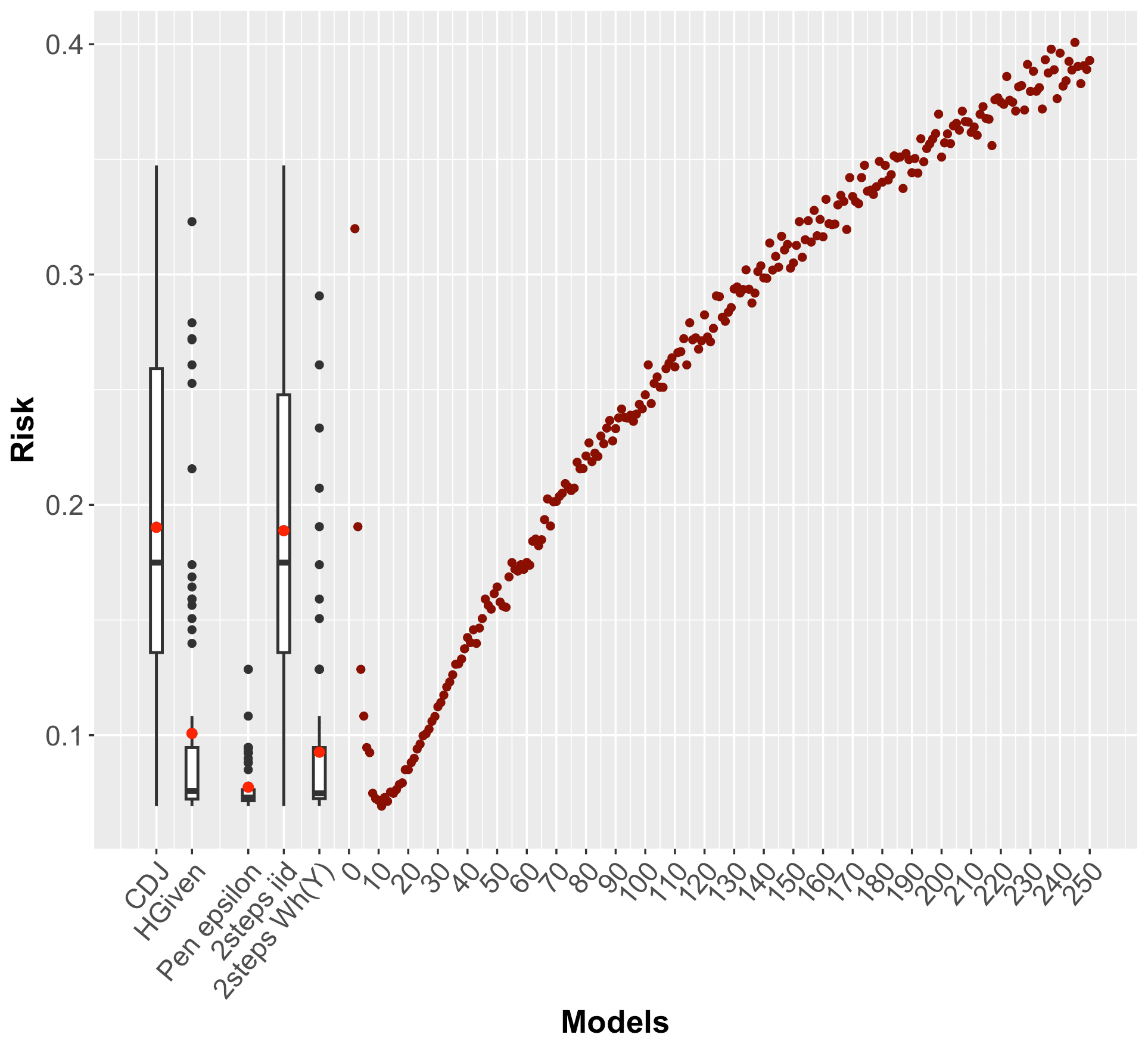}
\caption{$f_1$, $n=500$}
\end{subfigure}
\hfill
\begin{subfigure}{0.48\textwidth}
\centering
\includegraphics[width=\linewidth]{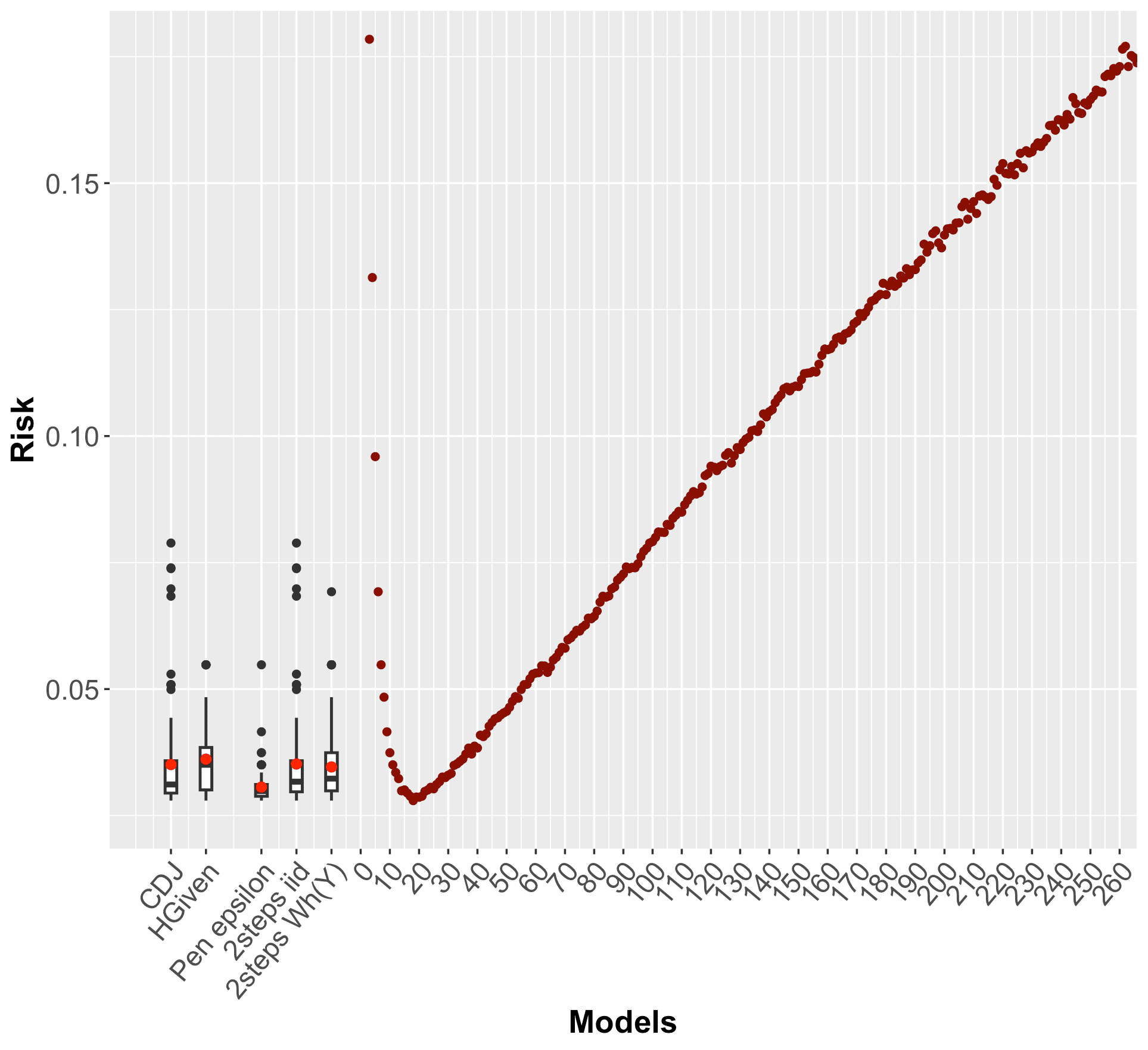}
\caption{$f_1$, $n=2000$}
\end{subfigure}

\vspace{0.3cm}

\begin{subfigure}{0.48\textwidth}
\centering
\includegraphics[width=\linewidth]{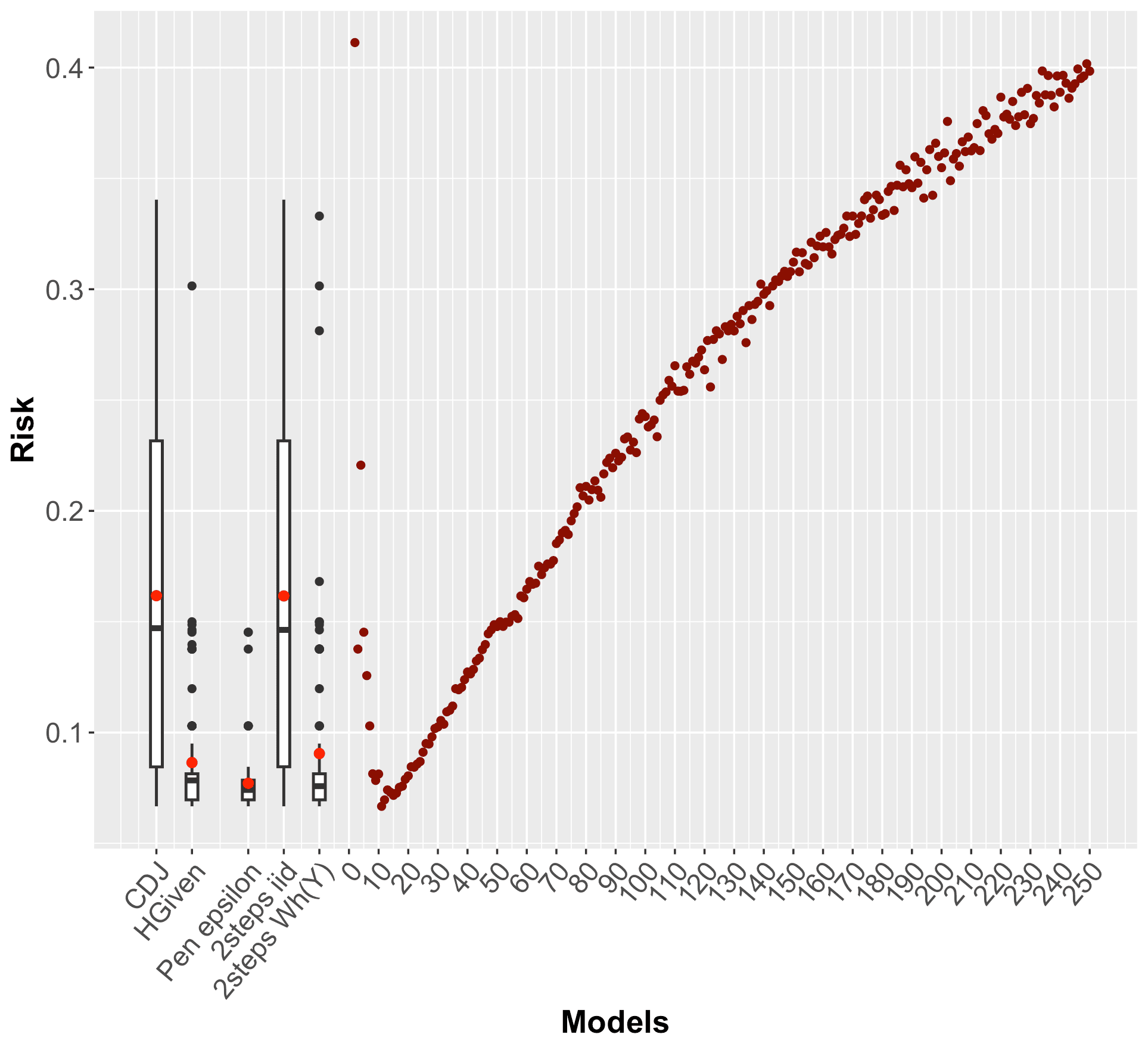}
\caption{$f_2$, $n=500$}
\end{subfigure}
\hfill
\begin{subfigure}{0.48\textwidth}
\centering
\includegraphics[width=\linewidth]{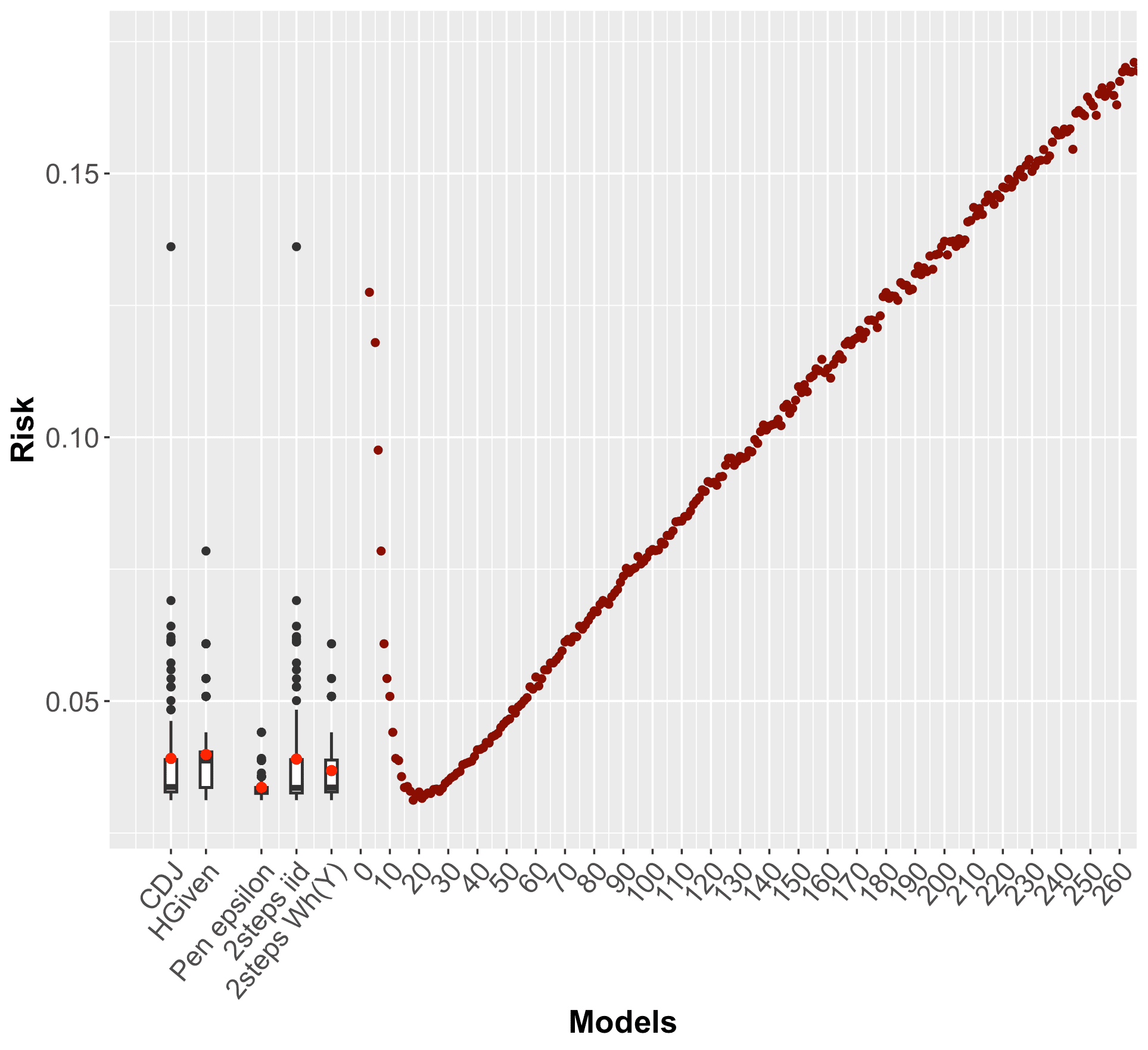}
\caption{$f_2$, $n=2000$}
\end{subfigure}
\caption{Risk performance of the penalization procedures for the two regression functions and two sample sizes in Experiment~6. The error process is a FGN  with $H=0.7$,  and the design is the non Gaussian Markov chain \eqref{markov_chain_X} with $a=1.5$.}

\label{fig:long_range_dependence_error_FGN07_design_betamix1_5}
\end{figure}

$\bullet$  {\bf Experiment 6}: the error process is a FGN process with $H=0.7$, and the design process is the non Gaussian Markov chain \eqref{markov_chain_X} with $a=1.5$, see Fig.~\ref{fig:long_range_dependence_error_FGN07_design_betamix1_5}.

 \medskip 

\begin{figure}[htbp]
\centering
\begin{subfigure}{0.48\textwidth}
\centering
\includegraphics[width=\linewidth]{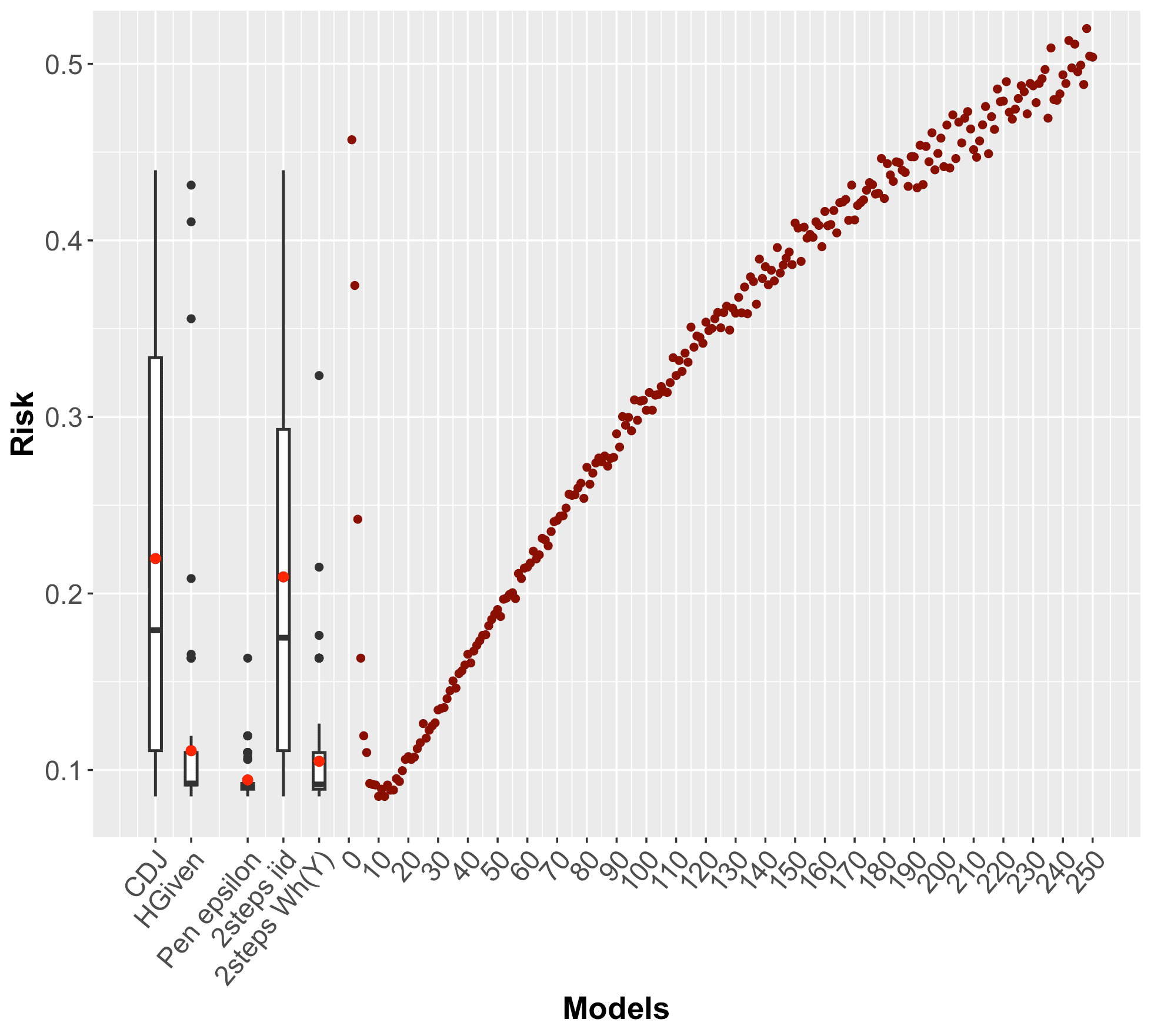}
\caption{$f_1$, $n=500$}
\end{subfigure}
\hfill
\begin{subfigure}{0.48\textwidth}
\centering
\includegraphics[width=\linewidth]{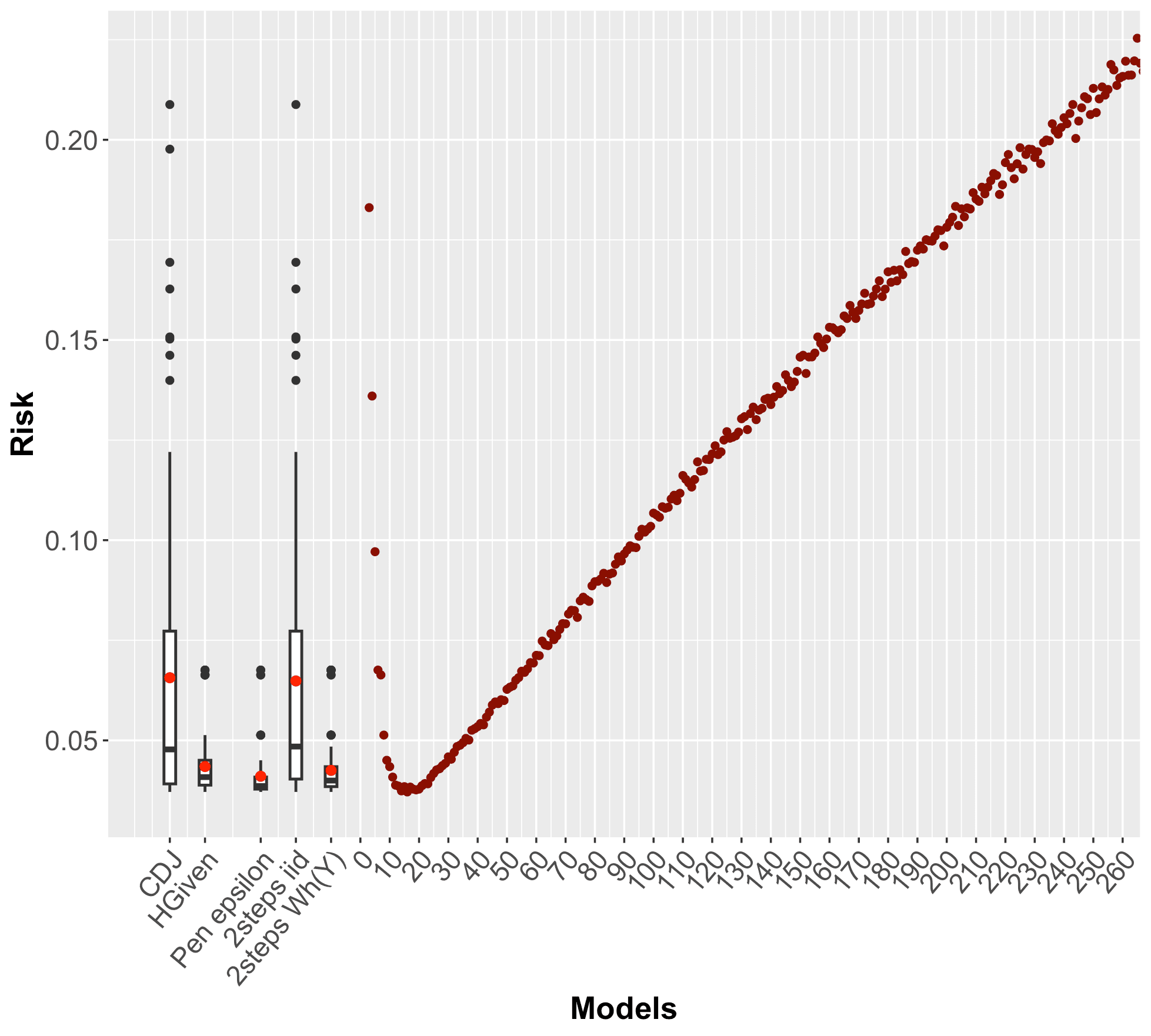}
\caption{$f_1$, $n=2000$}
\end{subfigure}

\vspace{0.3cm}

\begin{subfigure}{0.48\textwidth}
\centering
\includegraphics[width=\linewidth]{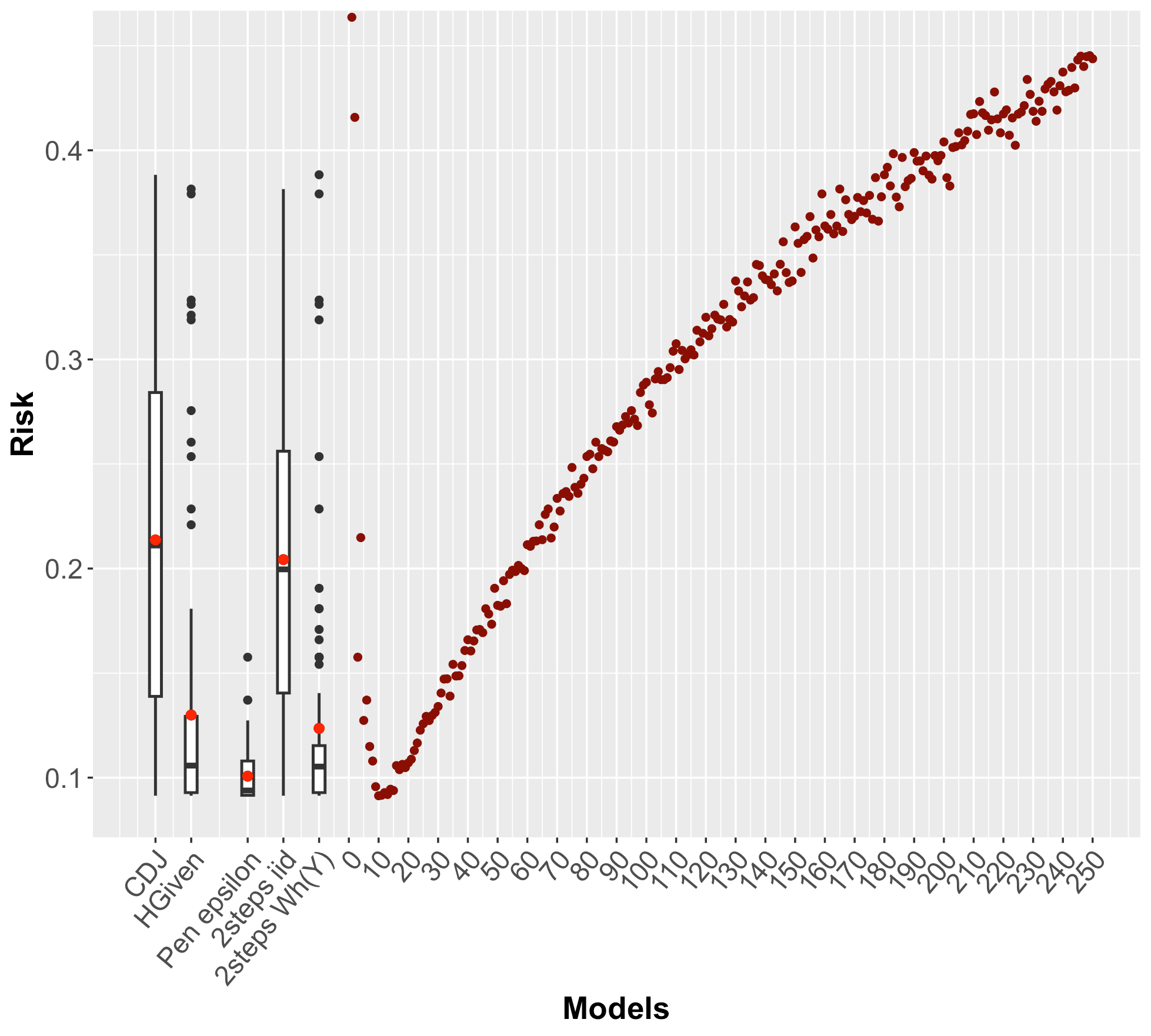}
\caption{$f_2$, $n=500$}
\end{subfigure}
\hfill
\begin{subfigure}{0.48\textwidth}
\centering
\includegraphics[width=\linewidth]{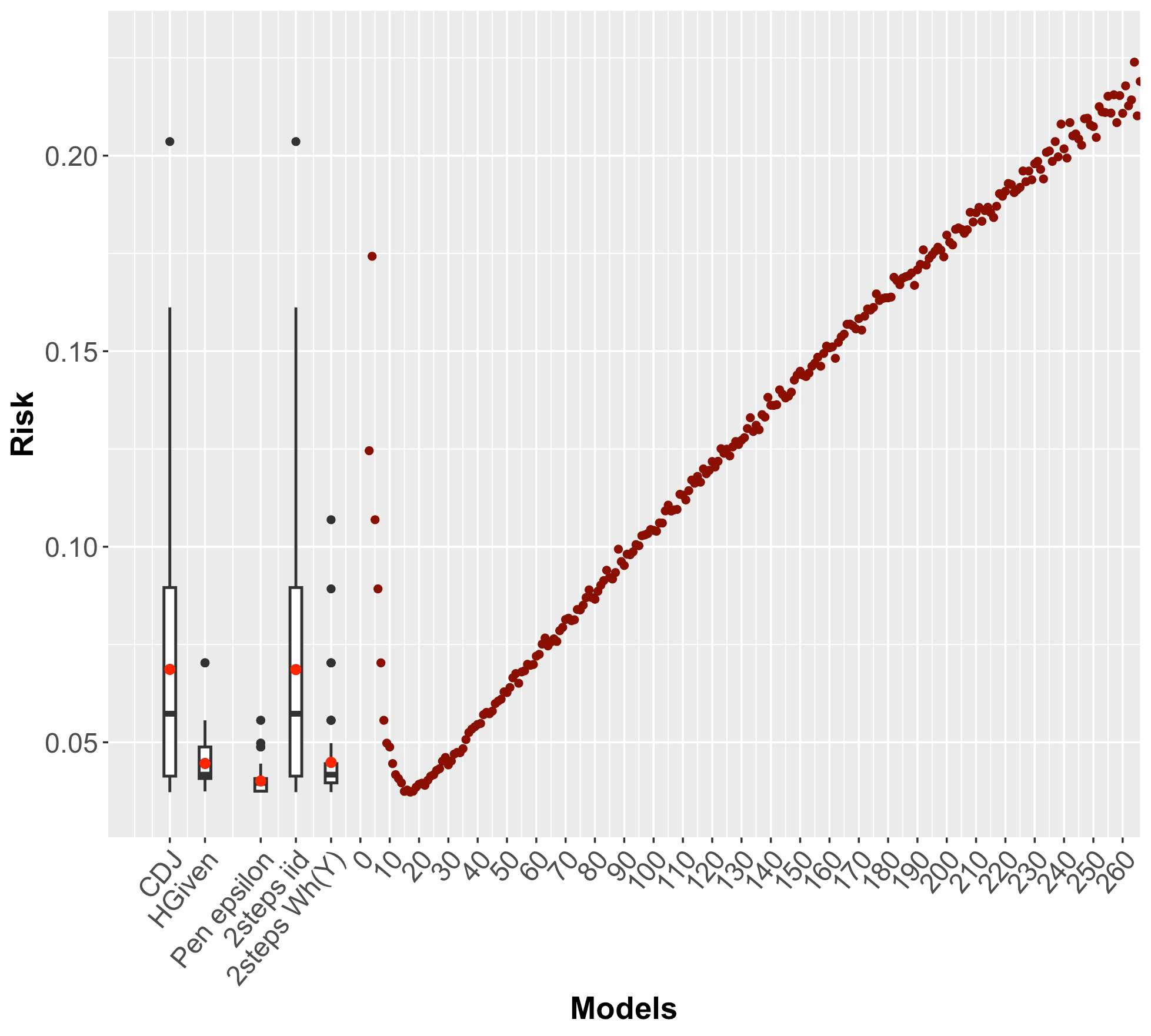}
\caption{$f_2$, $n=2000$}
\end{subfigure}
\caption{Risk performance of the penalization procedures for the two regression functions and two sample sizes in Experiment~7. The error process is an heteroscedastic FGN  with $H=0.7$ (as in Experience 5, see\eqref{Heteps}),  and the design  is the non Gaussian Markov chain \eqref{markov_chain_X} with $a=0.7$.}
\label{fig:long_range_dependence_error_FGNhetero07_design_betamix07}
\end{figure}

$\bullet$ {\bf Experiment 7}: the error process is the heteroscedastic process defined as in Experiment 5, except that the underlying FGN process has Hurst parameter $H=0.7$. The design process is the non Gaussian Markov chain~\eqref{markov_chain_X} with $a=0.7$, see Figure~\ref{fig:long_range_dependence_error_FGNhetero07_design_betamix07}.

\medskip

Experiments~6 and~7 correspond to long-memory settings, though not extreme long-memory cases. In this context, it is clearly evident that it is preferable to use the penalty corrected by the Hurst coefficient rather than the penalty proportional to the dimension, which performs poorly. These observations are consistent with the results presented in Section~\ref{sec:main_longrange}.

\medskip 

\begin{figure}[htbp]
\centering
\begin{subfigure}{0.48\textwidth}
\centering
\includegraphics[width=\linewidth]{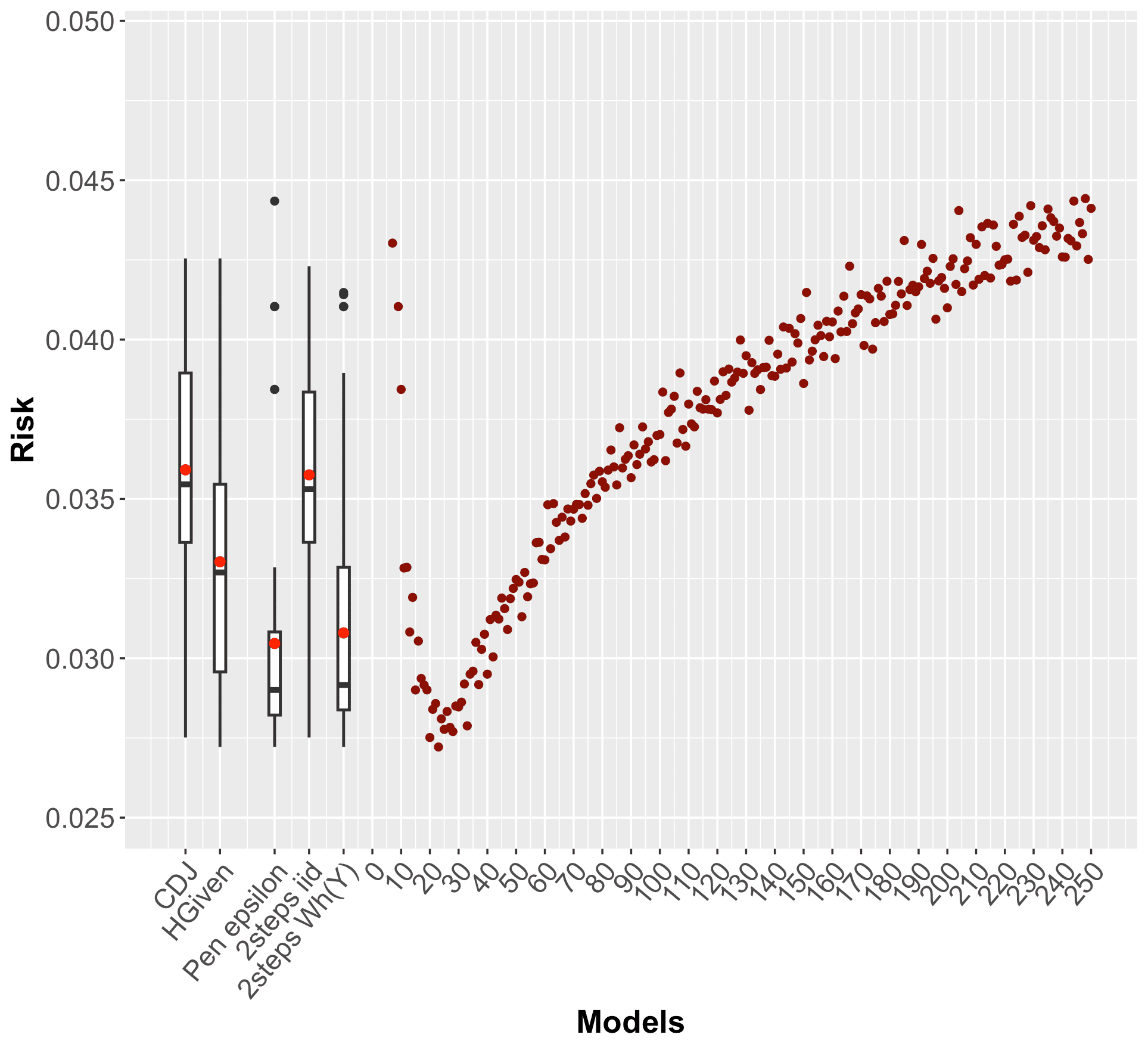}
\caption{$f_1$, $n=500$}
\end{subfigure}
\hfill
\begin{subfigure}{0.48\textwidth}
\centering
\includegraphics[width=\linewidth]{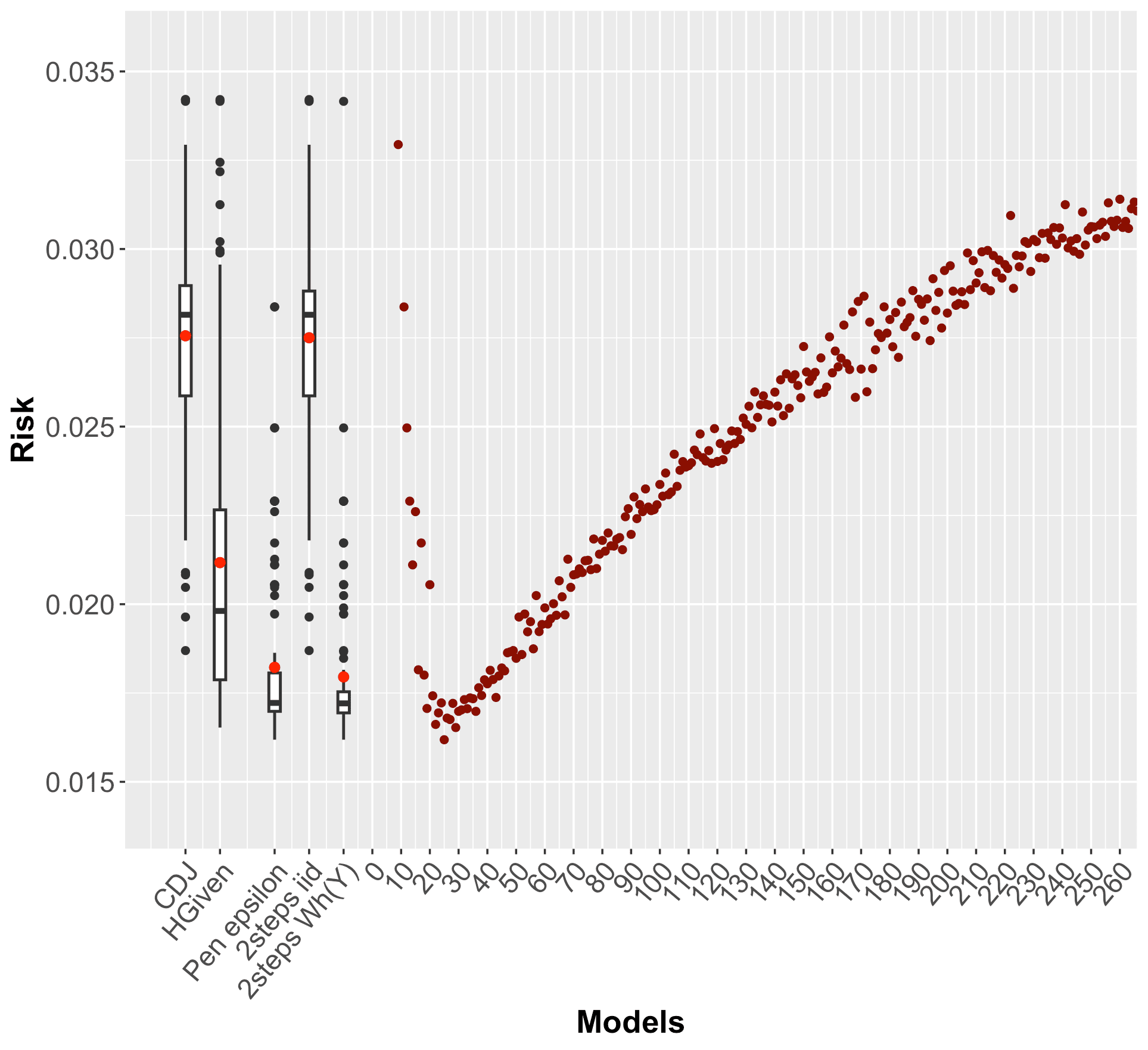}
\caption{$f_1$, $n=2000$}
\end{subfigure}

\vspace{0.3cm}

\begin{subfigure}{0.48\textwidth}
\centering
\includegraphics[width=\linewidth]{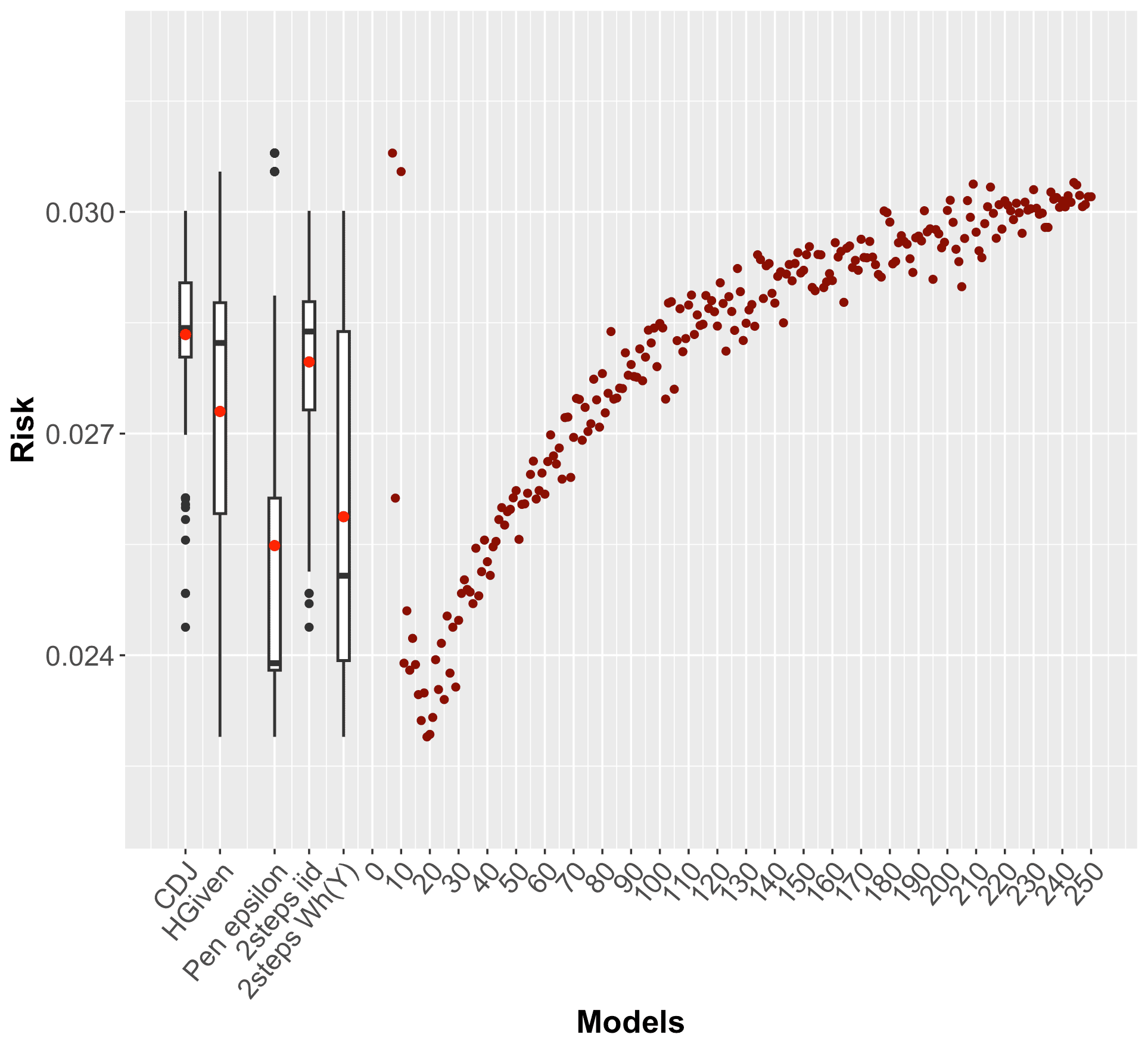}
\caption{$f_2$, $n=500$}
\end{subfigure}
\hfill
\begin{subfigure}{0.48\textwidth}
\centering
\includegraphics[width=\linewidth]{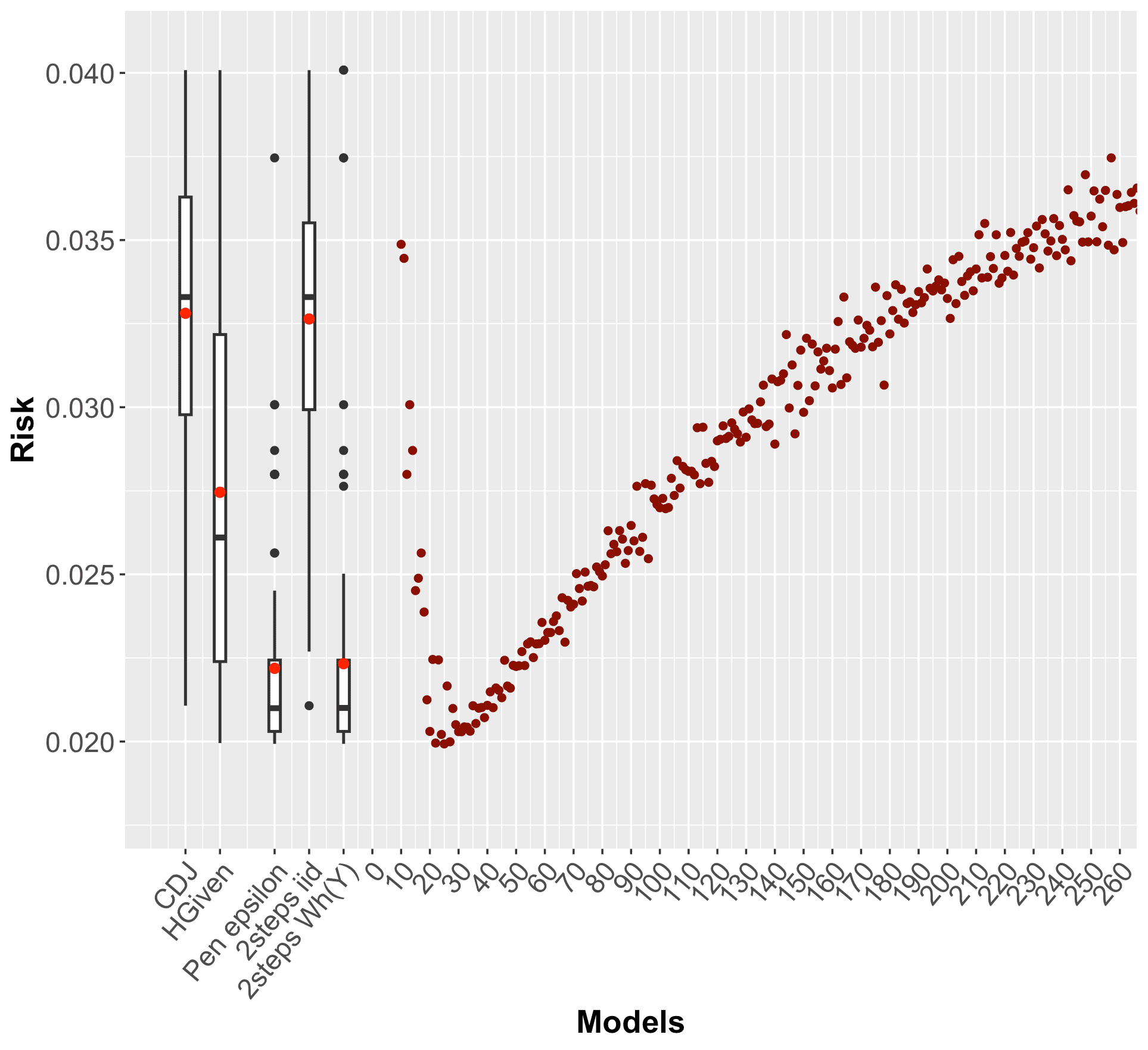}
\caption{$f_2$, $n=2000$}
\end{subfigure}

\caption{Risk performance of the penalization procedures for the two regression functions and two sample sizes in Experiment~8. The error process is the non Gaussian Markov chain \eqref{markov_chain_eps} with $a=0.5$,  and the design  is the non Gaussian Markov chain \eqref{markov_chain_X} with $a=0.3$}
\label{fig:long_range_dependence_error_betamix05_design_betamix03}
\end{figure}

$\bullet$ {\bf Experiment 8}: the error process is the centered non Gaussian Markov chain \eqref{markov_chain_eps} with $a=0.5$, and the design process is the non Gaussian Markov chain~\eqref{markov_chain_X} with $a=0.3$, see Figure~\ref{fig:long_range_dependence_error_betamix05_design_betamix03}.

\medskip 

\medskip

\begin{figure}[htbp]
\centering
\begin{subfigure}{0.48\textwidth}
\centering
\includegraphics[width=\linewidth]{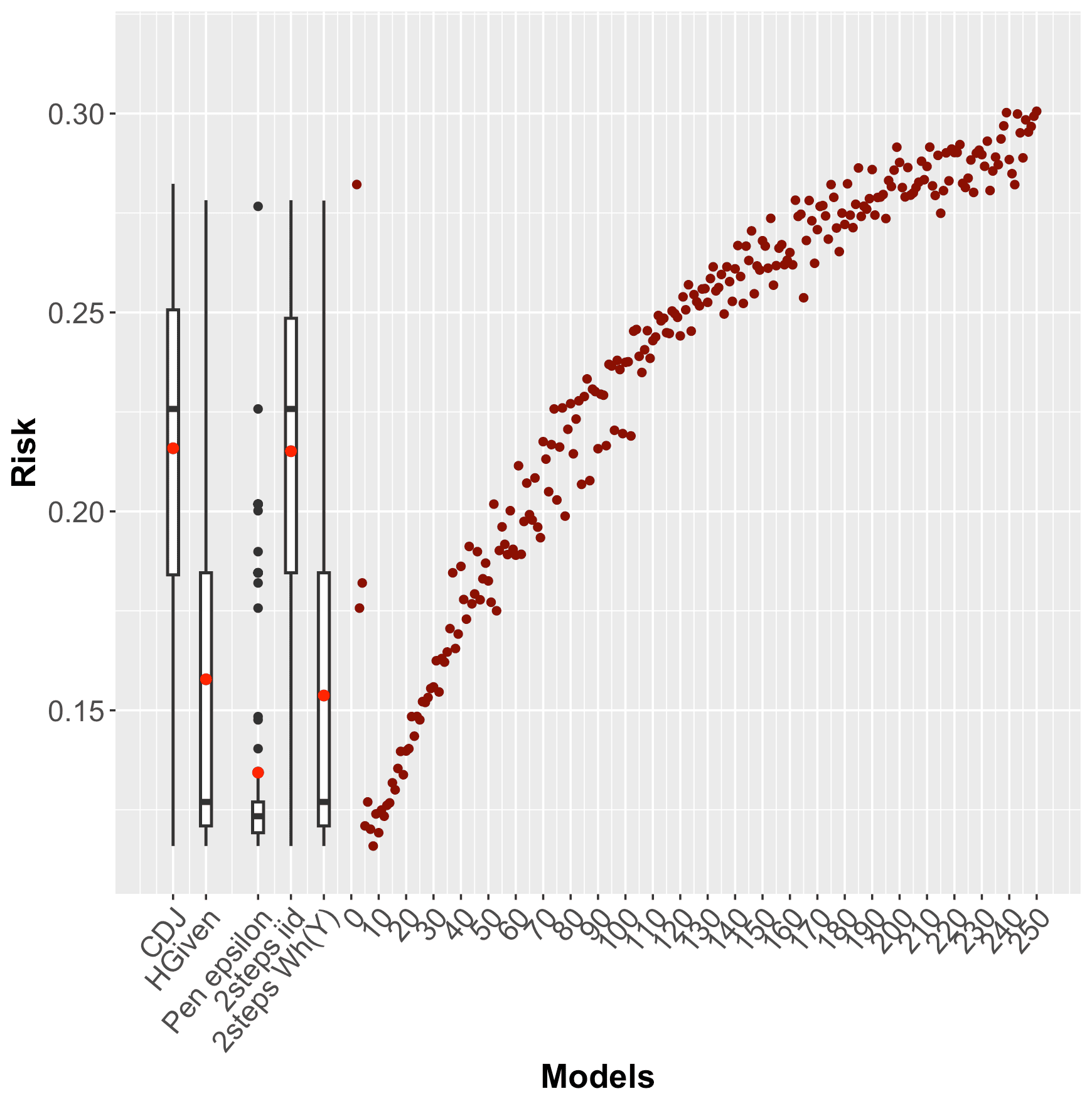}
\caption{$f_1$, $n=500$}
\end{subfigure}
\hfill
\begin{subfigure}{0.48\textwidth}
\centering
\includegraphics[width=\linewidth]{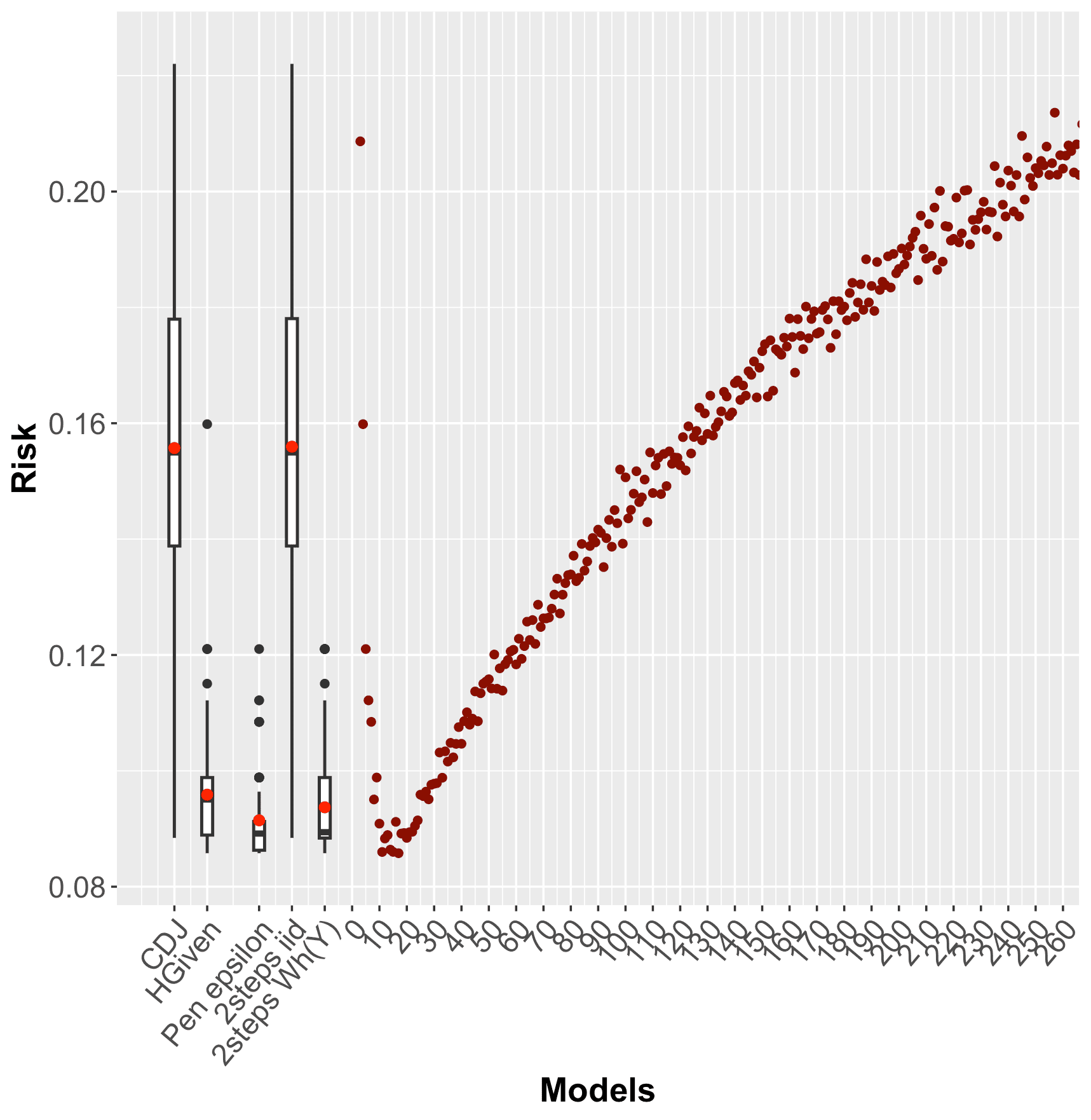}
\caption{$f_1$, $n=2000$}
\end{subfigure}

\vspace{0.3cm}

\begin{subfigure}{0.48\textwidth}
\centering
\includegraphics[width=\linewidth]{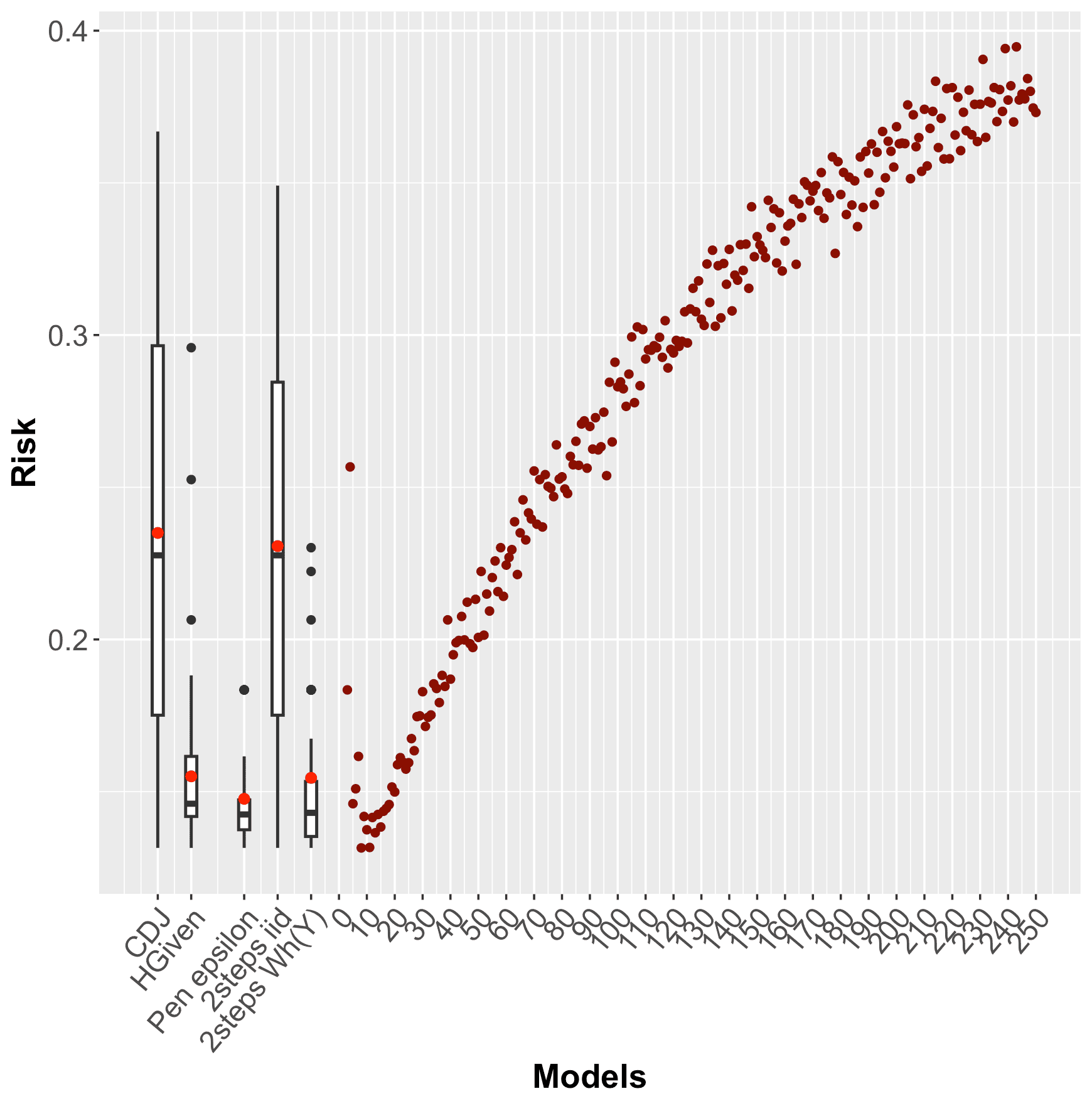}
\caption{$f_2$, $n=500$}
\end{subfigure}
\hfill
\begin{subfigure}{0.48\textwidth}
\centering
\includegraphics[width=\linewidth]{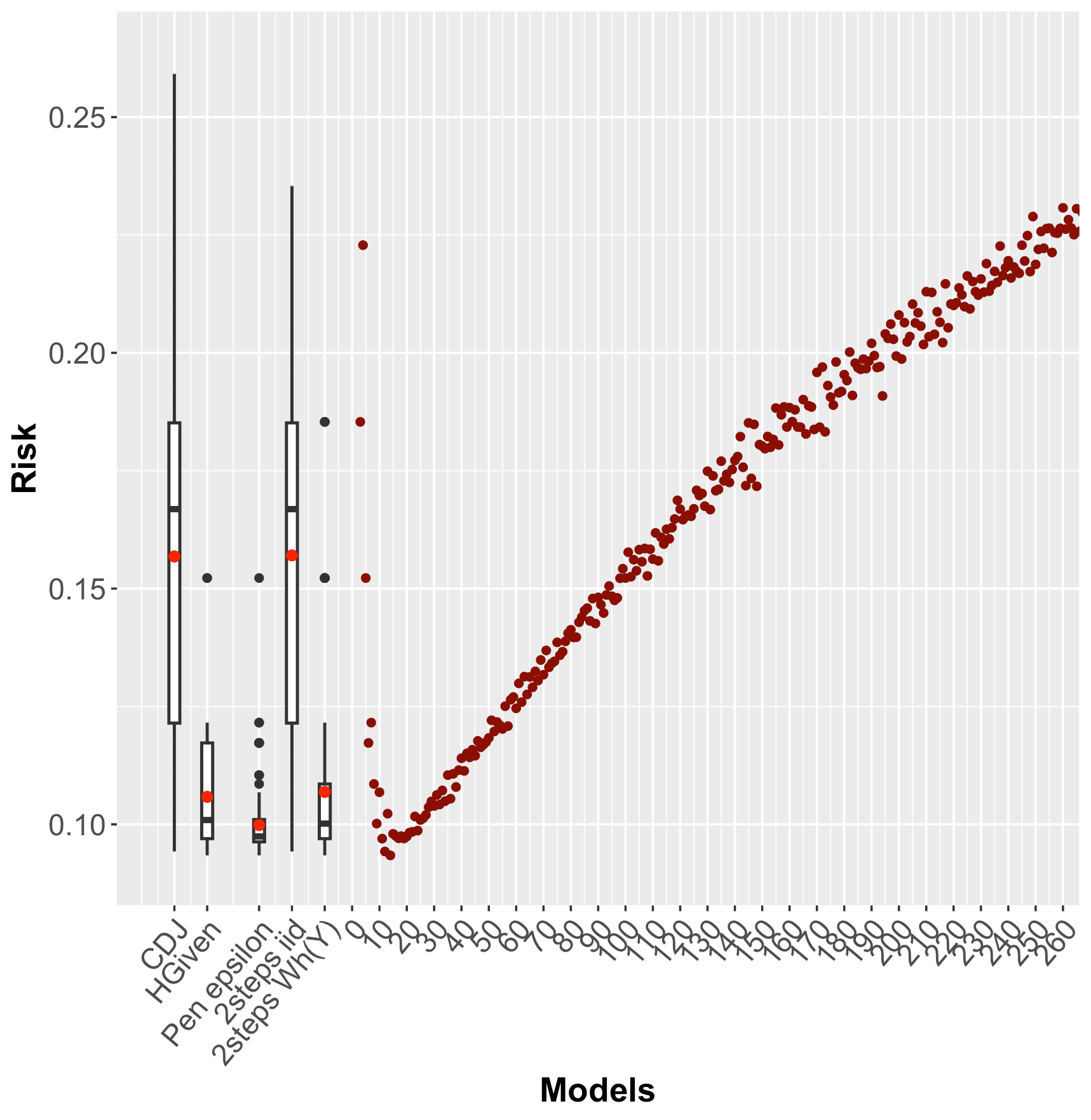}
\caption{$f_2$, $n=2000$}
\end{subfigure}

\caption{Risk performance of the penalization procedures for the two regression functions and two sample sizes in Experiment~9. The error process is a FGN  with $H=0.8$, and the design  is the non Gaussian Markov  chain \eqref{markov_chain_X} with $a=0.3$.}
\label{fig:long_range_dependence_error_FGN08_design_betamix03}
\end{figure}

$\bullet$ {\bf Experiment 9}: The error process is a FGN process with $H=0.8$,  and the design process is the non Gaussian Markov chain~\eqref{markov_chain_X} with $a=0.3$, see Figure~\ref{fig:long_range_dependence_error_FGN08_design_betamix03}.

\medskip

\medskip 

\begin{figure}[htbp]
\centering
\begin{subfigure}{0.48\textwidth}
\centering
\includegraphics[width=\linewidth]{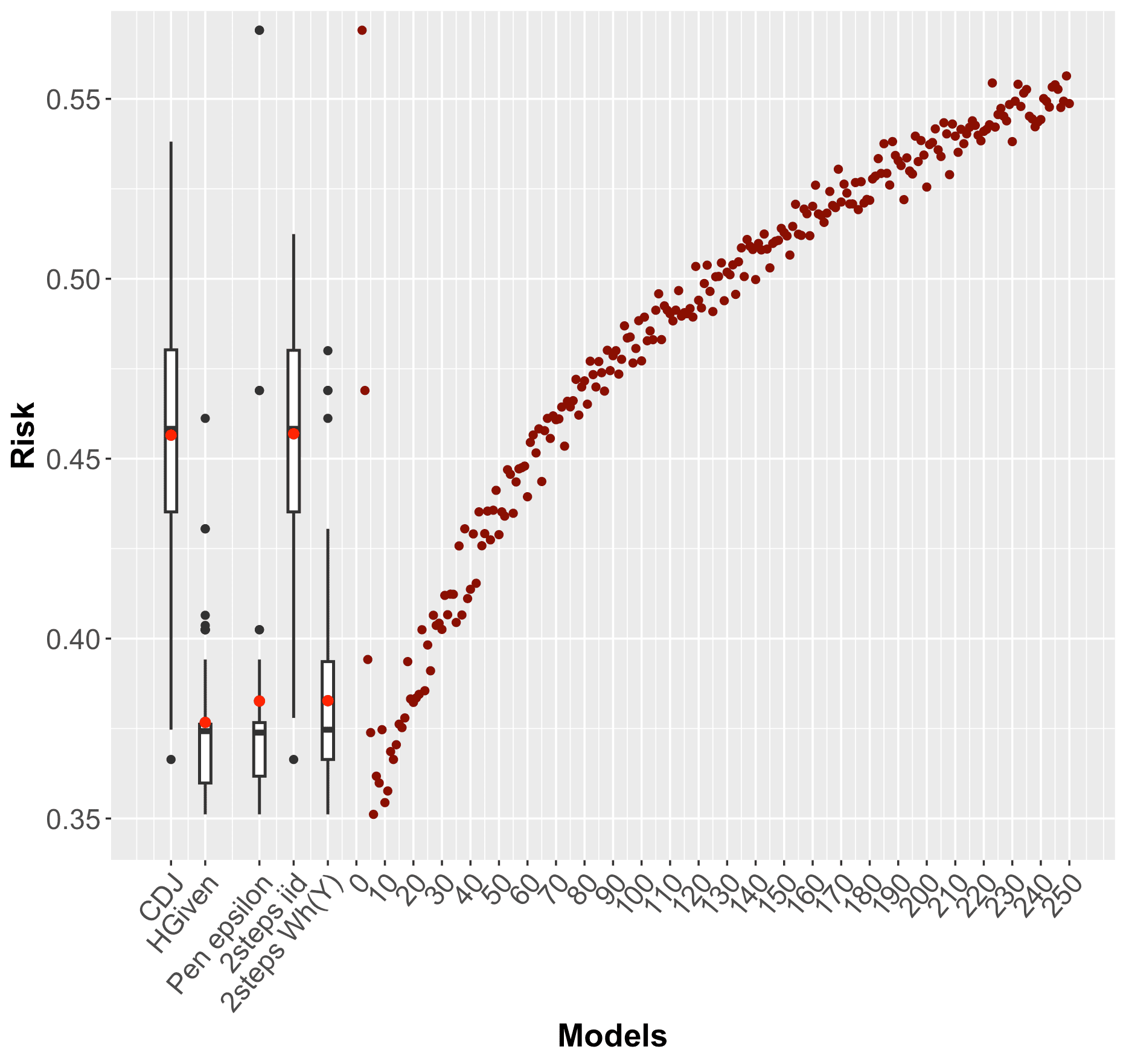}
\caption{$f_1$, $n=500$}
\end{subfigure}
\hfill
\begin{subfigure}{0.48\textwidth}
\centering
\includegraphics[width=\linewidth]{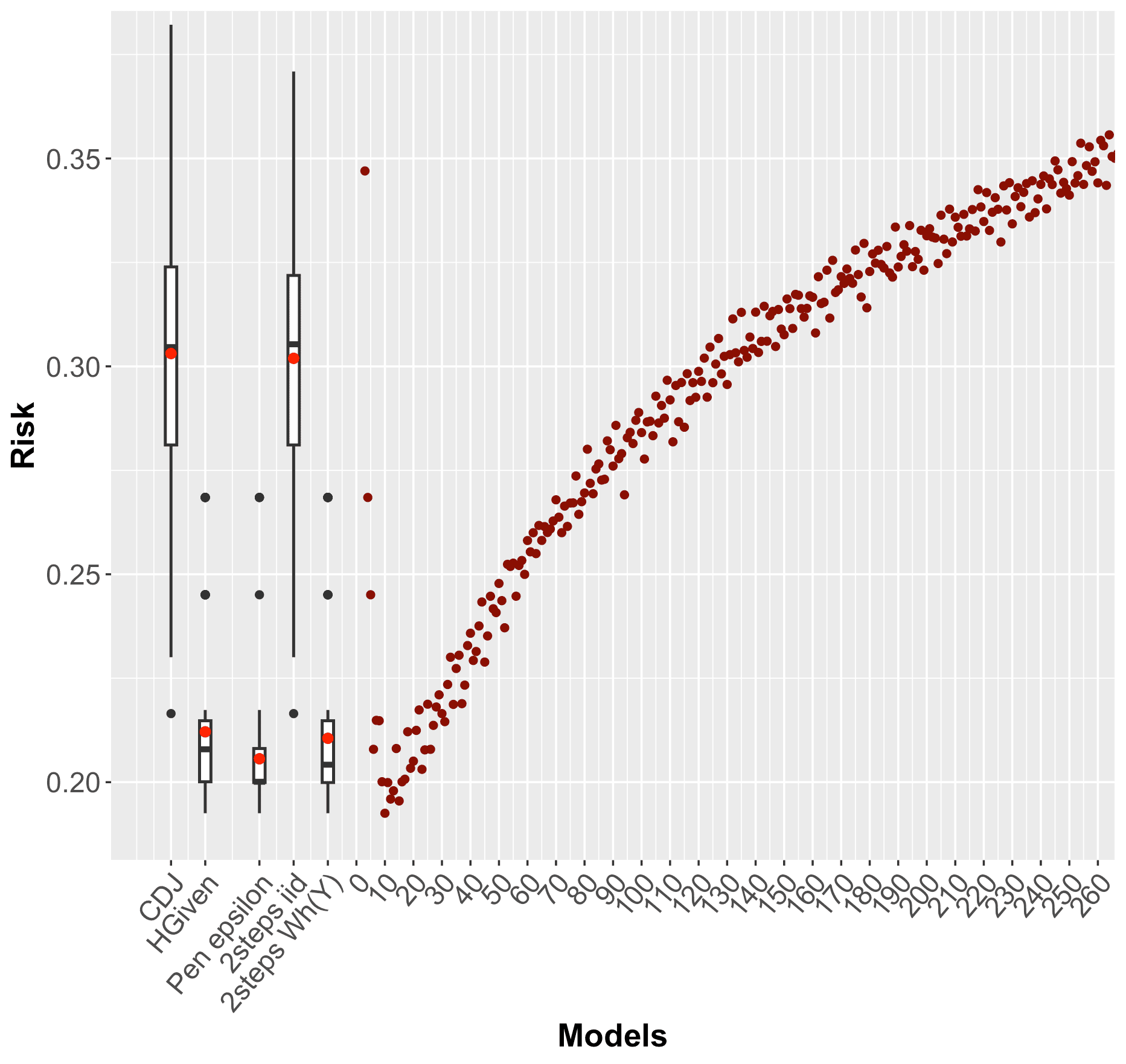}
\caption{$f_1$, $n=2000$}
\end{subfigure}

\vspace{0.3cm}

\begin{subfigure}{0.48\textwidth}
\centering
\includegraphics[width=\linewidth]{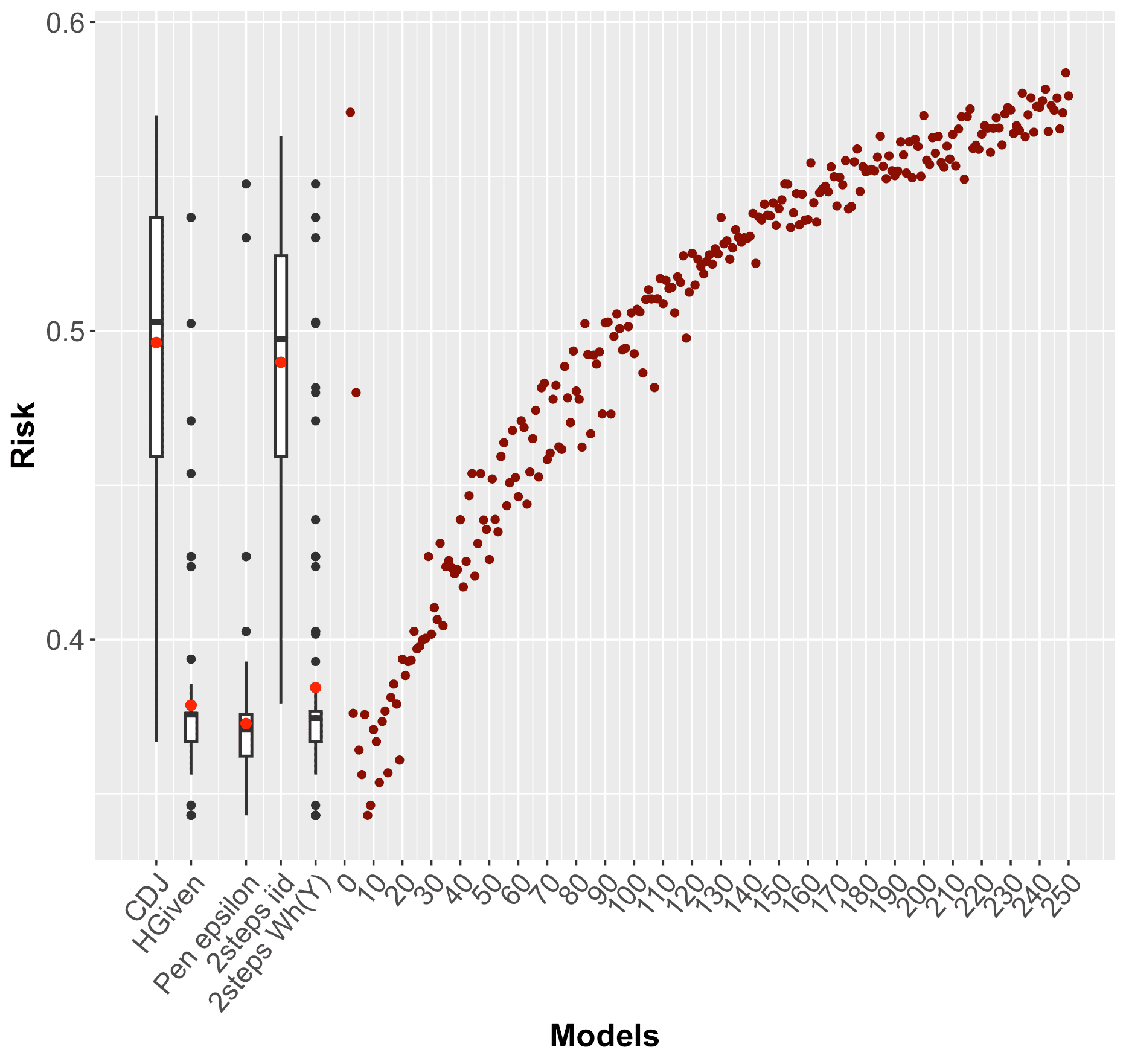}
\caption{$f_2$, $n=500$}
\end{subfigure}
\hfill
\begin{subfigure}{0.48\textwidth}
\centering
\includegraphics[width=\linewidth]{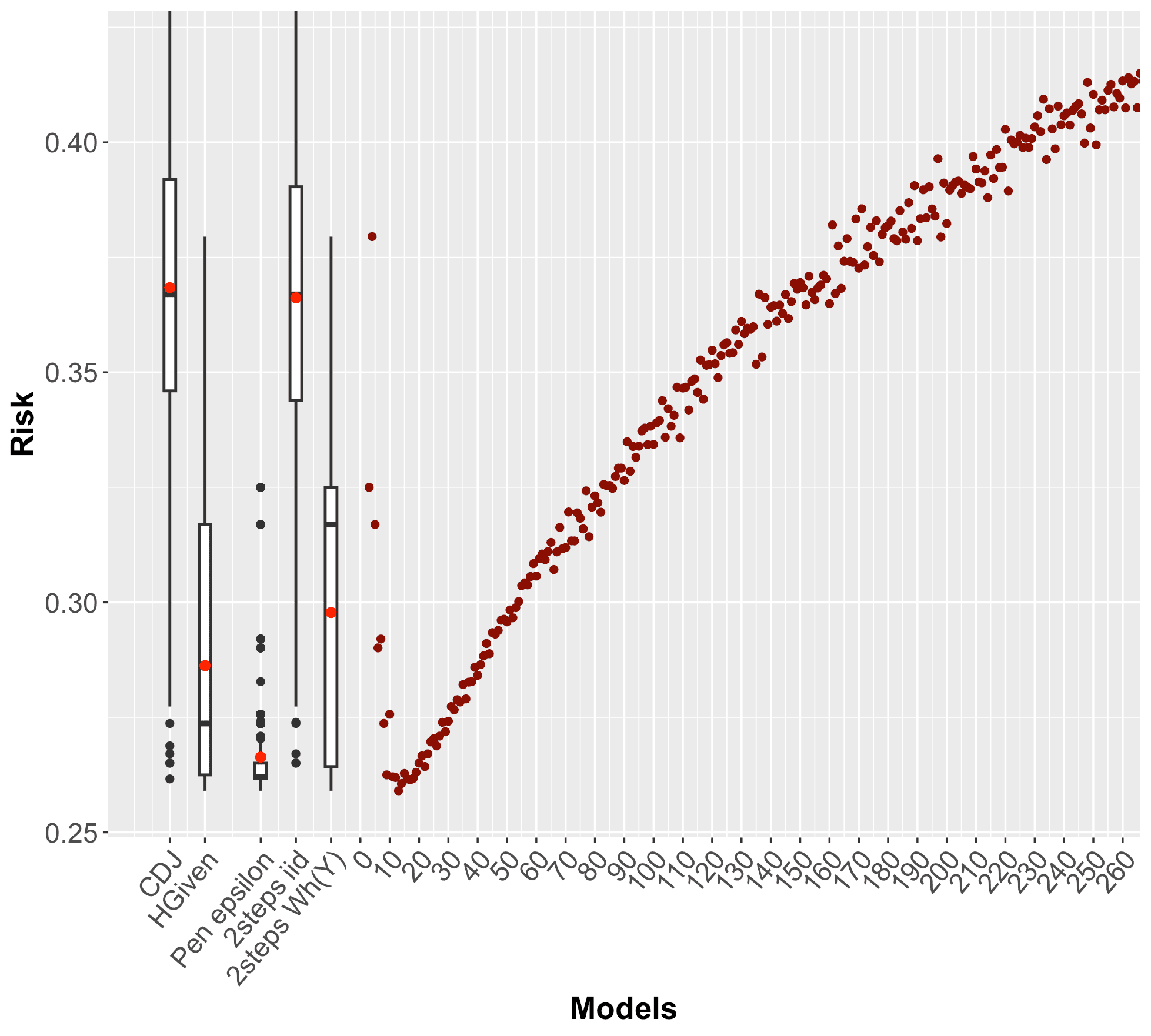}
\caption{$f_2$, $n=2000$}
\end{subfigure}

\caption{Risk performance of the penalization procedures for the two regression functions and two sample sizes in Experiment~10. The error process is a FGN  with $H=0.9$,  and the design  is the non Gaussian Markov chain \eqref{markov_chain_X} with $a=0.3$.}
\label{fig:long_range_dependence_error_FGN09_design_betamix03}
\end{figure}

$\bullet$ {\bf Experiment  10}: The error process is the FGN process with $H=0.9$, and the design process is the non Gaussian Markov chain~\eqref{markov_chain_X} with $a=0.3$, see Figure~\ref{fig:long_range_dependence_error_FGN09_design_betamix03}.

\medskip

Experiments 8, 9, and 10 focus on configurations with extremely long memory. Once again, we observe that the penalty proportional to the dimension leads to very poor performance. The penalty corrected by the Hurst coefficient yields significantly better results, even though, for small samples (\( n = 500 \)), the selected model may be far from the oracle.

\subsubsection{Long range dependence with $\rho$-dependent designs}
\label{subsec533}

We now consider the setting described in Section~\ref{sec:main_longrange2}, with a design process which is $\rho$-dependent.

\medskip 

\begin{figure}[htbp]
\centering
\begin{subfigure}{0.48\textwidth}
\centering
\includegraphics[width=\linewidth]{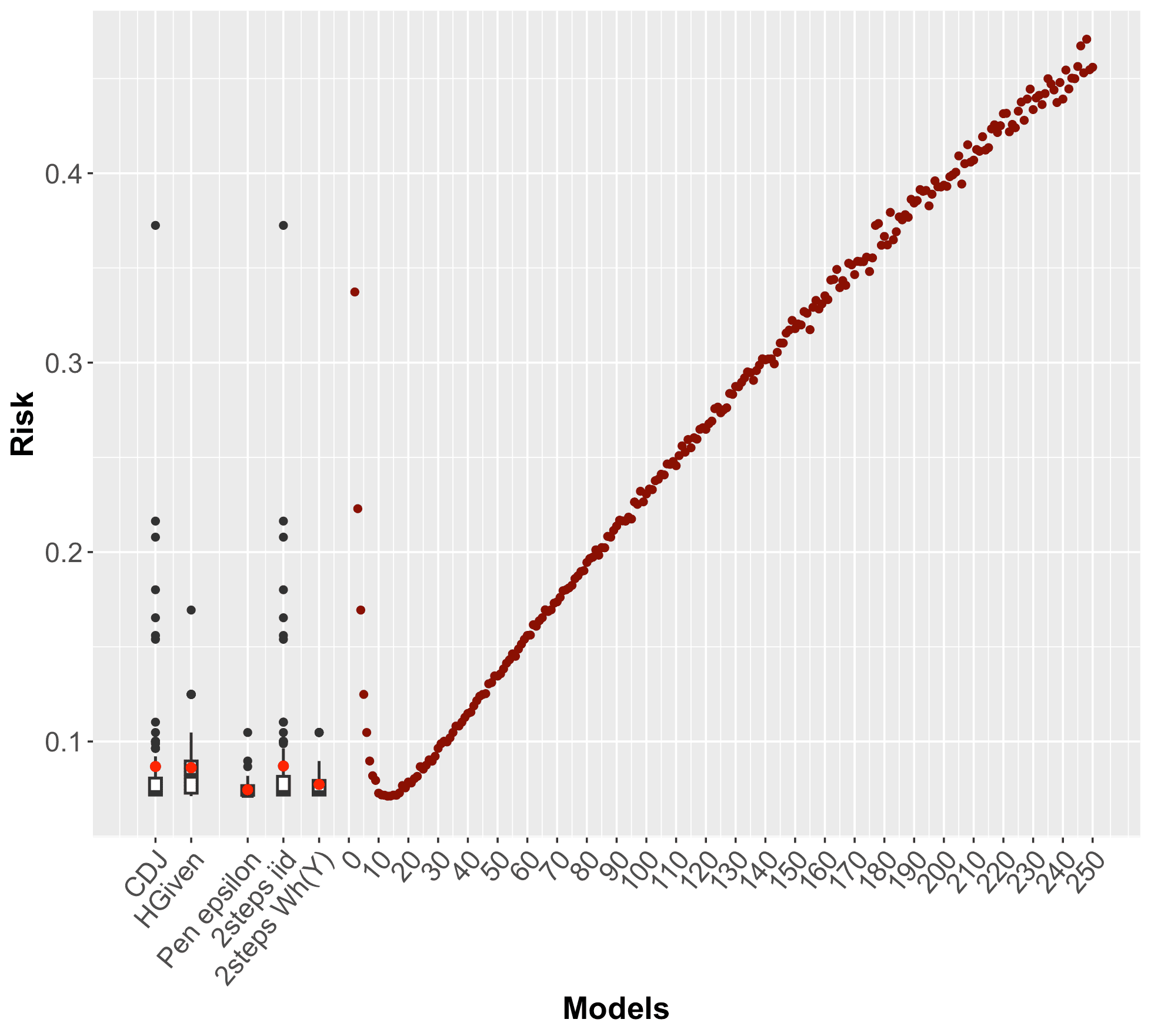}
\caption{$f_1$, $n=500$}
\end{subfigure}
\hfill
\begin{subfigure}{0.48\textwidth}
\centering
\includegraphics[width=\linewidth]{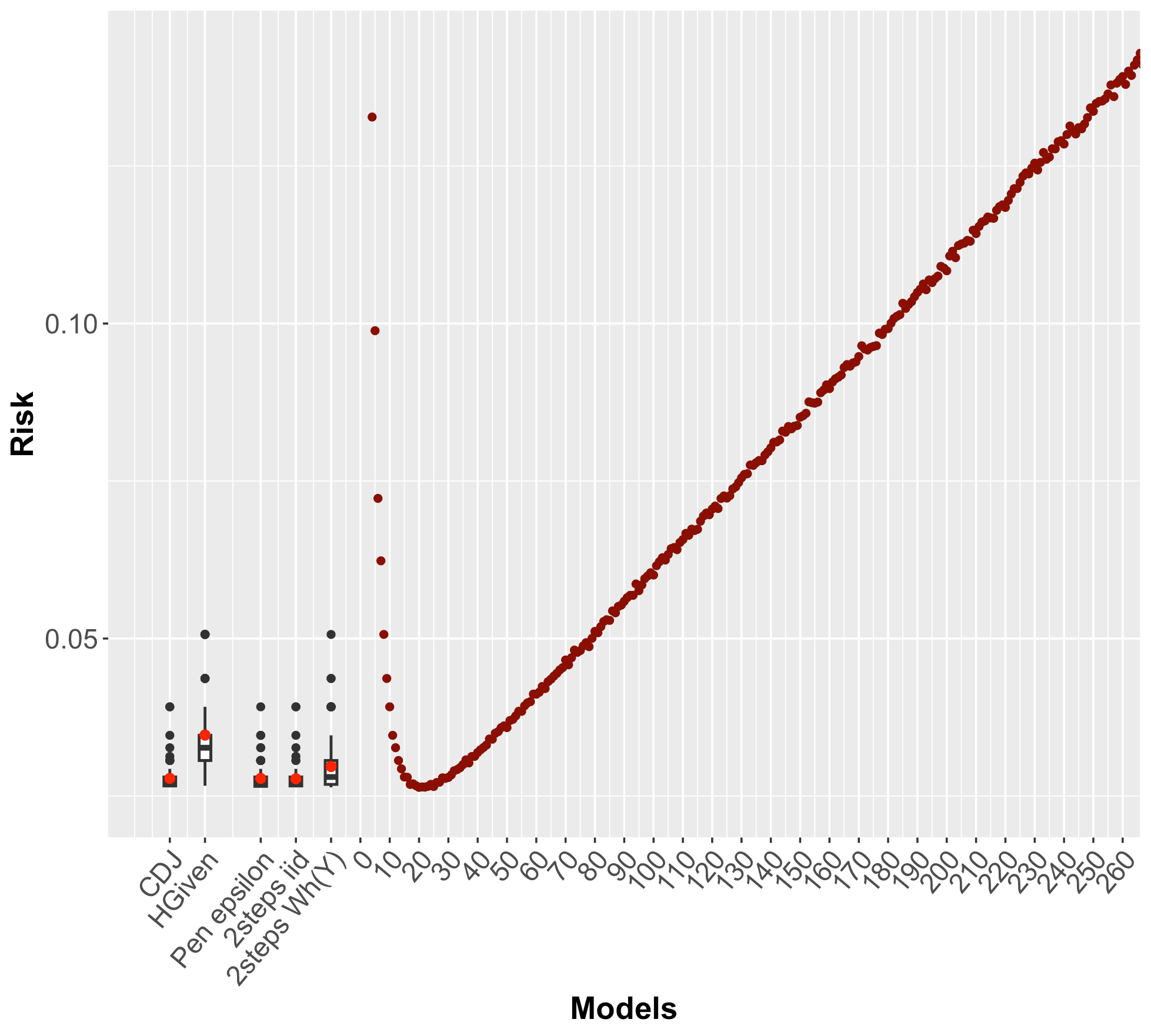}
\caption{$f_1$, $n=2000$}
\end{subfigure}

\vspace{0.3cm}

\begin{subfigure}{0.48\textwidth}
\centering
\includegraphics[width=\linewidth]{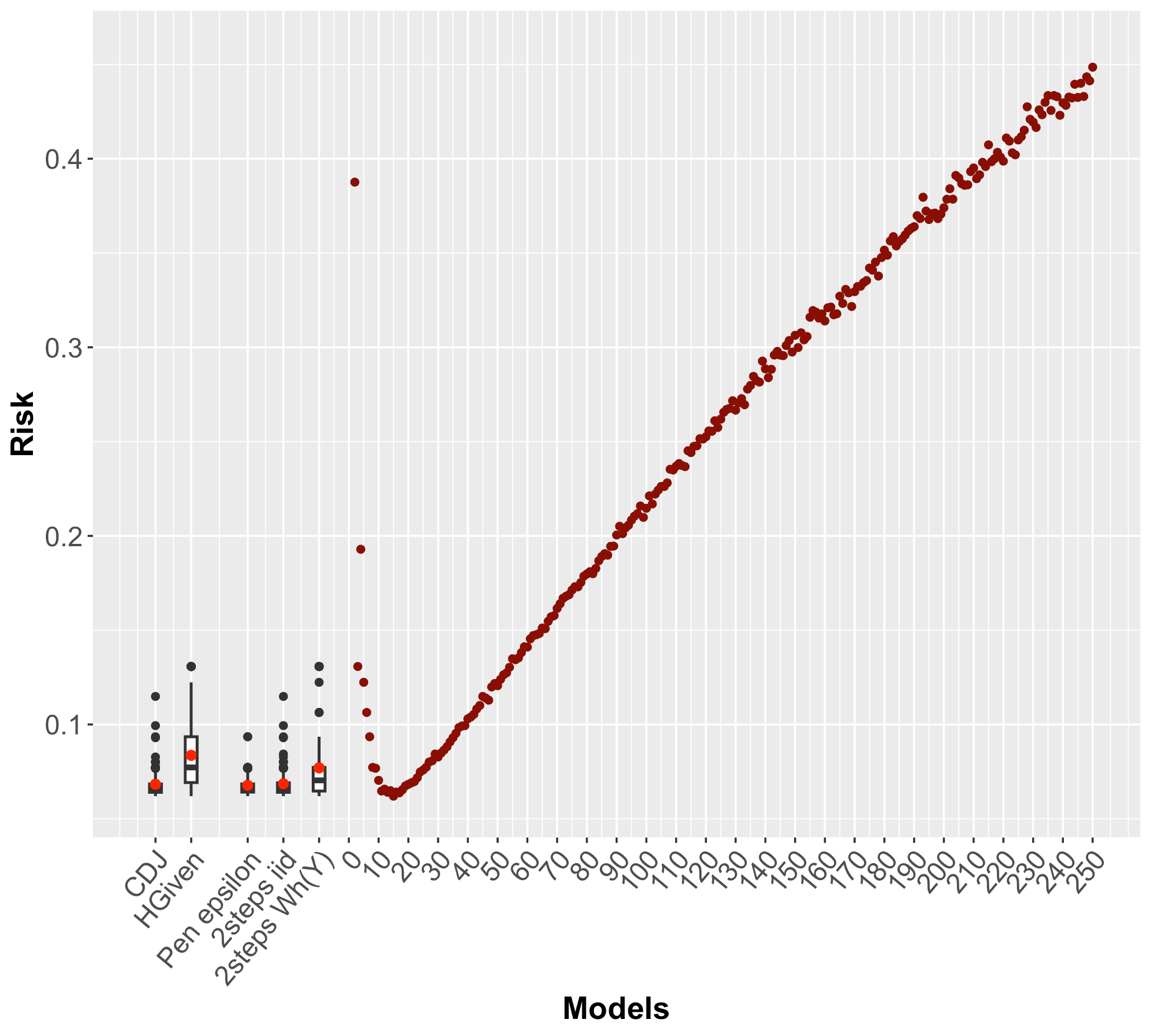}
\caption{$f_2$, $n=500$}
\end{subfigure}
\hfill
\begin{subfigure}{0.48\textwidth}
\centering
\includegraphics[width=\linewidth]{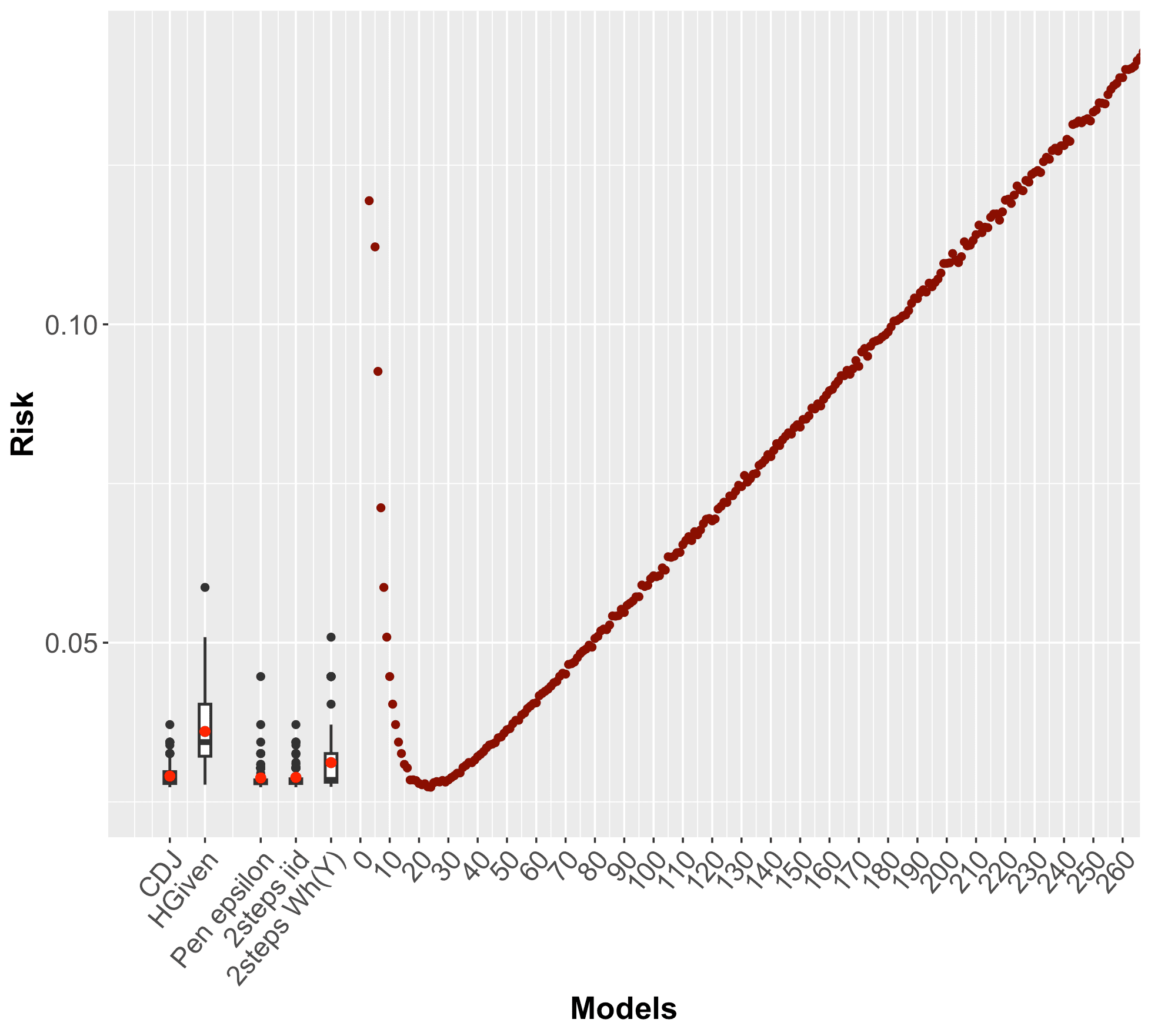}
\caption{$f_2$, $n=2000$}
\end{subfigure}
\caption{Risk performance of the penalization procedures for the two regression functions and two sample sizes in Experiment~11. The error process is a FGN with $H=0.7$, 
and the design  is obtained from  a FGN  with $H=0.7$ via the transformation \eqref{transFGN}.}
\label{fig:long_range_dependence_error_FGN07_design_FGN07}
\end{figure}

$\bullet$ {\bf Experiment 11} : the error process is a FGN with $H=0.7$, 
and the design  is obtained from  a FGN  with $H=0.7$ via the transformation \eqref{transFGN},
 see Figure~\ref{fig:long_range_dependence_error_FGN07_design_FGN07}. 
 

\medskip 

We verify here that the dimension-proportional penalty indeed yields satisfactory results, as stated in Proposition~\ref{prop_ld_improved} and the subsequent remarks. In this context, the Hurst coefficient-adjusted penalty proves to be suboptimal.

\medskip

\subsection{Fully data driven penalties}
\label{subsec:FullyDataDriven}

In practice, the Hurst coefficients are unknown, and more generally neither the nature of the error process nor that of the design process is available. Consequently, the strategies proposed above cannot be applied directly. We now introduce data-driven approaches for model selection that rely solely on the observed data. Although these methods do not come with the theoretical guarantees established in the previous sections, they nevertheless perform well empirically.

For estimating the Hurst index $H$ when it is unknown, we will use the Whittle MLE estimator~\cite{whittle1953estimation}
implemented in the {\ttfamily longmemo} package. 

To construct such data-driven penalties, let us return to the main term of the initial penalty \eqref{GenPen}, 
which is of the form 
$$
  \pen_{\varepsilon} =
\frac {\kappa}{(n({\mathbf X})\vee 1)}
 \text{tr}\left ( \text{Proj}_{S_m} \Sigma_I({\mathbf X}_n)\right )=
 \kappa \mathbb{E} \left[\left \| \text{Proj}_{S_m}(\varepsilon) \right \|_{n({\mathbf X})}^{2} \big | {\mathbf X}_n \right] \, .
$$
Of course, this quantity is not directly usable in practice either. We propose to use instead the naive quantity
$\left \| \text{Proj}_{S_m}(\varepsilon) \right \|_{n({\mathbf X})}^{2}$.
To obtain a penalty that is nondecreasing with the model dimension, we apply isotonic regression to the function
\begin{equation}
\label{eq:penvarpsilondirect}
m \mapsto      \left \| \text{Proj}_{S_m}(\varepsilon) \right \|_{n({\mathbf X})}^{2} \, .
\end{equation}
This procedure also has the effect of smoothing the fluctuations of $\| \text{Proj}_{S_m}(\varepsilon)  \|_{n({\mathbf X})}^{2}$. We still denote by $\pen_{\varepsilon}$ the resulting penalty. Although it requires access to the unobserved errors and is therefore not implementable in practice, we display this penalty in all risk-comparison figures for reference.   

\medskip 

A natural idea is to replace the error vector $\varepsilon$ in~\eqref{eq:penvarpsilondirect} by a vector of residuals. To do so, one first needs to select a "reasonable" preliminary model. We investigate and compare the following two-step procedures throughout the experiments of this section.

\begin{itemize}
    \item {\bf Two-step iid procedure}:
    \begin{enumerate}
        \item Pre-select a model using the dimension jump method with a penalty proportional to the model dimension.
        \item Compute the residuals $\hat{\varepsilon}$ in the pre-selected model.
        \item Estimate the penalty in~\eqref{eq:penvarpsilondirect} by
        \[
        \pen_{\mathrm{2step\text{-}iid}}(m)
        = \kappa \,  \text{Isot}\left (
        \left \| \text{Proj}_{S_m}(\hat \varepsilon) \right \|_{n({\mathbf X})}^{2}\right ),
        \]
        where the term Isot means that we  applied isotonic regresssion to $m \mapsto      \left \| \text{Proj}_{S_m}(\varepsilon) \right \|_{n({\mathbf X})}^{2}$.
        \item Select the final model using the dimension jump method with penalty $\pen_{\mathrm{2step\text{-}iid}}(m)$.
    \end{enumerate}

    \item {\bf Two-step Whittle($Y$) procedure}:
    \begin{enumerate}
        \setcounter{enumi}{-1}
        \item Estimate the Hurst exponent $\hat H(Y)$  directly from the observations $Y$ using the Whittle estimator.
        \item Pre-select a model using the dimension jump method with a penalty proportional to $m^{2-2\hat H(Y)}$.  
        \item Compute the residuals $\hat{\varepsilon}$ in the pre-selected model.
        \item Estimate the penalty in~\eqref{eq:penvarpsilondirect} by
        \[
        \pen_{\mathrm{2step\text{-}Whittle}(Y)}(m)
        =
       \kappa \,  \text{Isot}\left (
        \left \| \text{Proj}_{S_m}(\hat \varepsilon) \right \|_{n({\mathbf X})}^{2}\right ).
        \]
        \item Select the final model using the dimension jump method with penalty $\pen_{\mathrm{2step\text{-}Whittle}(Y)}(m)$.
    \end{enumerate}
\end{itemize}

For all the experiments presented above, the penalty $\pen_{\varepsilon} $ consistently achieves the best performance, often approaching the oracle's performance, even in cases involving very long memory. 

In the case of short memory, both proposed penalties $\pen_{\mathrm{2step\text{-}iid}}$ and $\pen_{\mathrm{2step\text{-}Whittle}(Y)}$ show good performances also. Specifically, with the {\it two-step i.i.d. procedure}, the  pre-model selected by the  penalty proportional to the dimension correctly estimates the regression function, so that the residuals are ``close'' to the unknown errors. This approach thus yields  performance similar to that of $\pen_{\varepsilon}$. Regarding the {\it two-step Whittle($Y$)} method, the Whittle estimator successfully estimates the Hurst index directly from $Y$ in this short memory context. The selected pre-model is thus close to that of the {\it two-step i.i.d. procedure} method, with comparable final performances.

We now discuss the results for long memory.  For Experiments 6 and 7, in which the error process is a Fractional Gaussian Noise and the design is a non Gaussian Markov chain, the {\it two-step i.i.d.} procedure exhibits very poor performance. Indeed, the penalty proportional to the dimension selects an inadequate model. However, even with a naive estimation of the Hurst index directly from $Y$, the {\it two-step Whittle($Y$)} method achieves reasonable model selection. Ultimately, the performance of the {\it two-step Whittle($Y$)} method closely approaches that of the epsilon penalty $\pen_{\varepsilon}$, particularly for sufficiently large sample sizes $n$.

The situation is similar for the more challenging experiment 8, 9 and 10 involving long memory or very long memory. Again, the sample size $n$ must be sufficiently large for the {\it two-step Whittle($Y$)} procedure to function effectively, whereas the {\it two-step i.i.d.} procedure remains unable to select a reasonable model in any case. In the context of very long memory of Experiment 10 in the case of the function $f^\star_2$, the {\it two-step Whittle($Y$)} procedure reaches its limits even for large $n$. 

For Experiment 11, where the error is FGN and the design is $\rho$-dependent, the {\it two-step i.i.d.} procedure outperforms the {\it two-step Whittle($Y$)} method. This is not surprising since it is expected that a penalty proportional to the dimension is effective in this context. However, the {\it two-step Whittle($Y$)} procedure still yields reasonable results.

\begin{figure}[htbp]
\centering
\begin{subfigure}{0.48\textwidth}
\centering
\includegraphics[width=\linewidth]{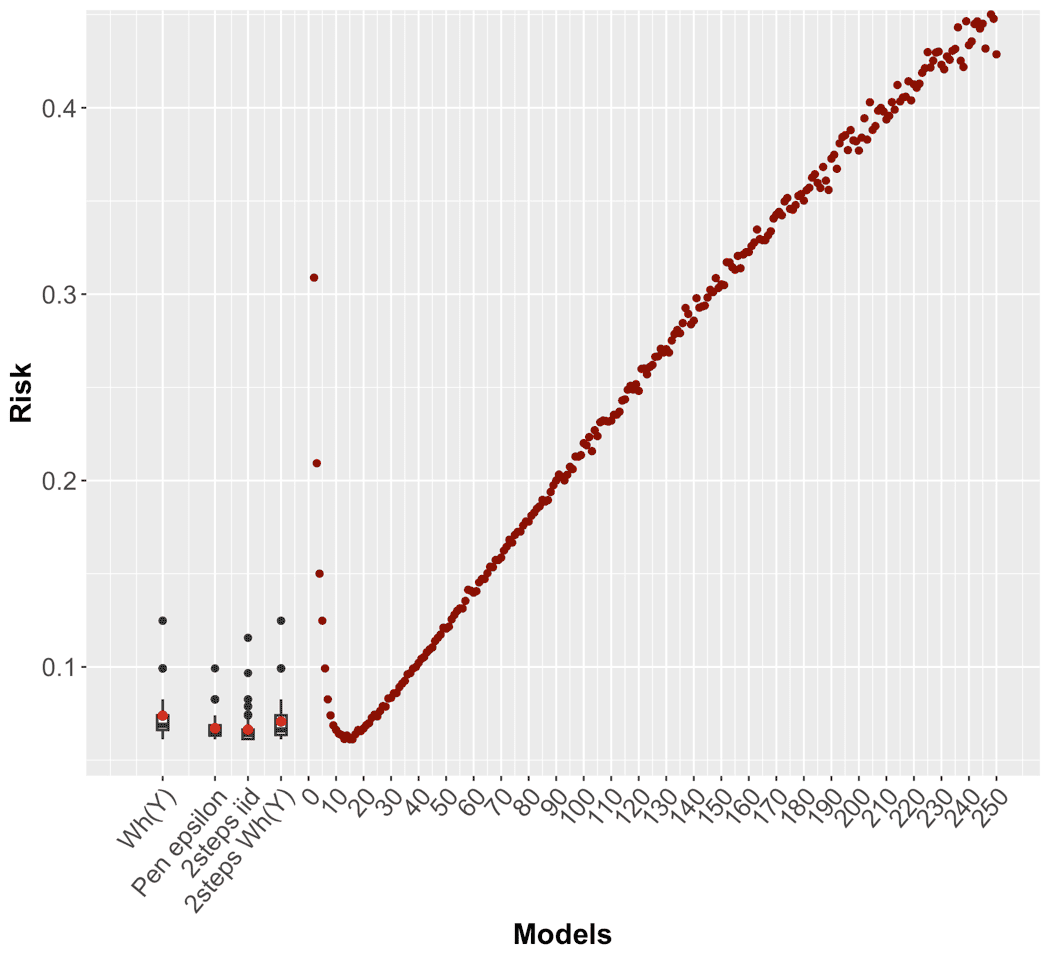}
\caption{$f_1$, $n=500$}
\end{subfigure}
\hfill
\begin{subfigure}{0.48\textwidth}
\centering
\includegraphics[width=\linewidth]{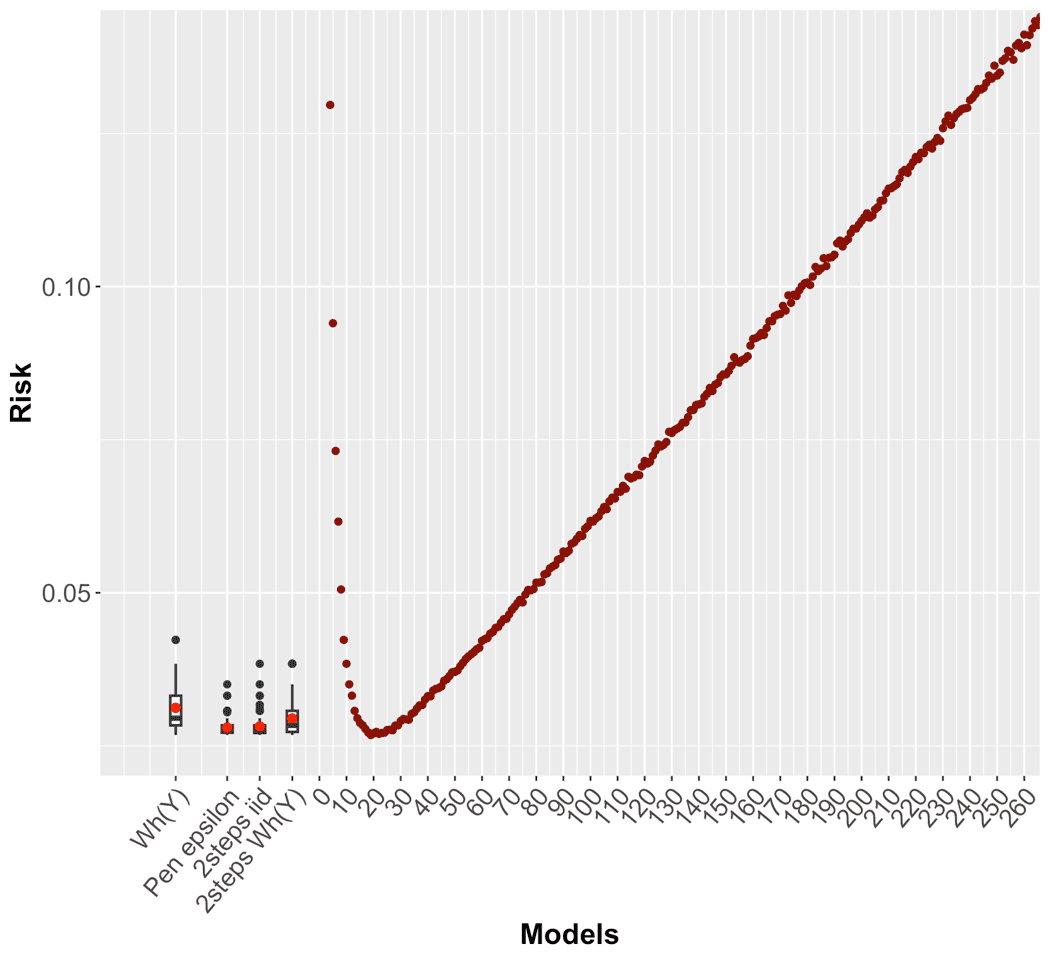}
\caption{$f_1$, $n=2000$}
\end{subfigure}

\vspace{0.3cm}

\begin{subfigure}{0.48\textwidth}
\centering
\includegraphics[width=\linewidth]{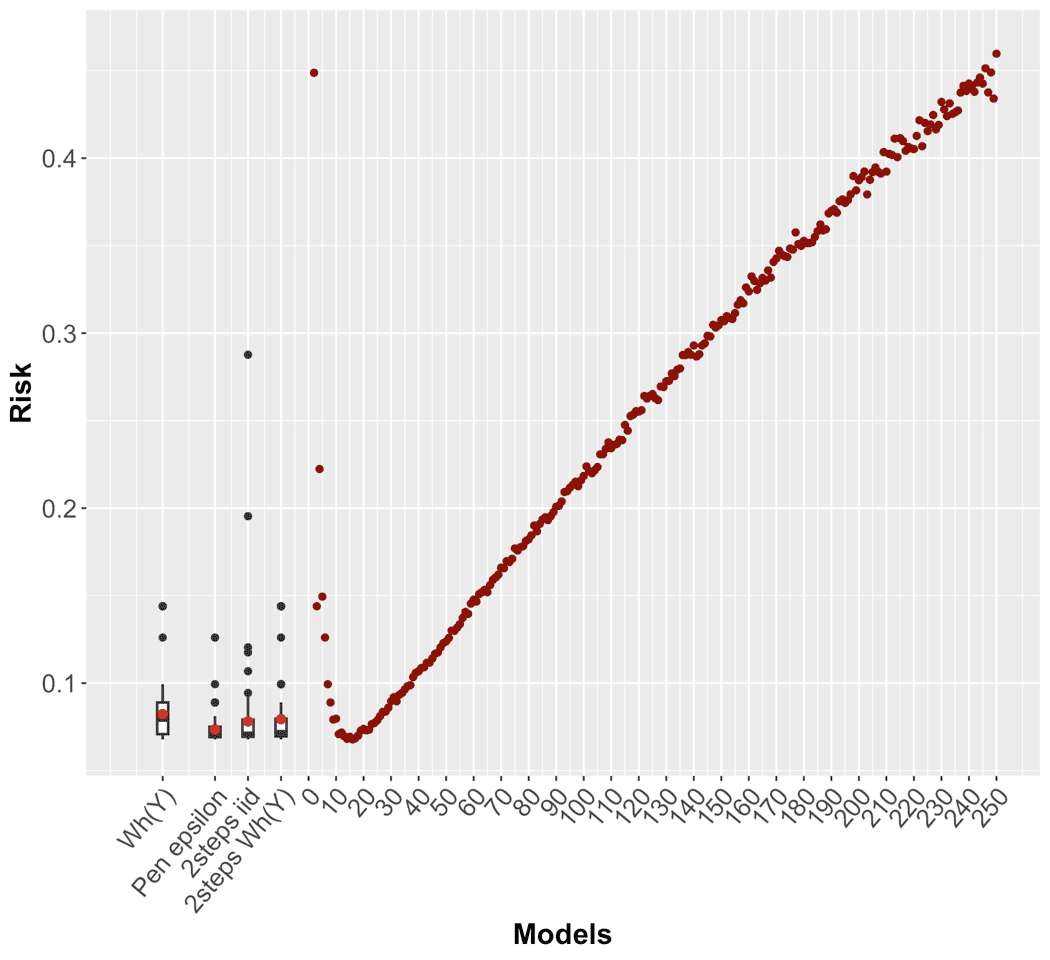}
\caption{$f_2$, $n=500$}
\end{subfigure}
\hfill
\begin{subfigure}{0.48\textwidth}
\centering
\includegraphics[width=\linewidth]{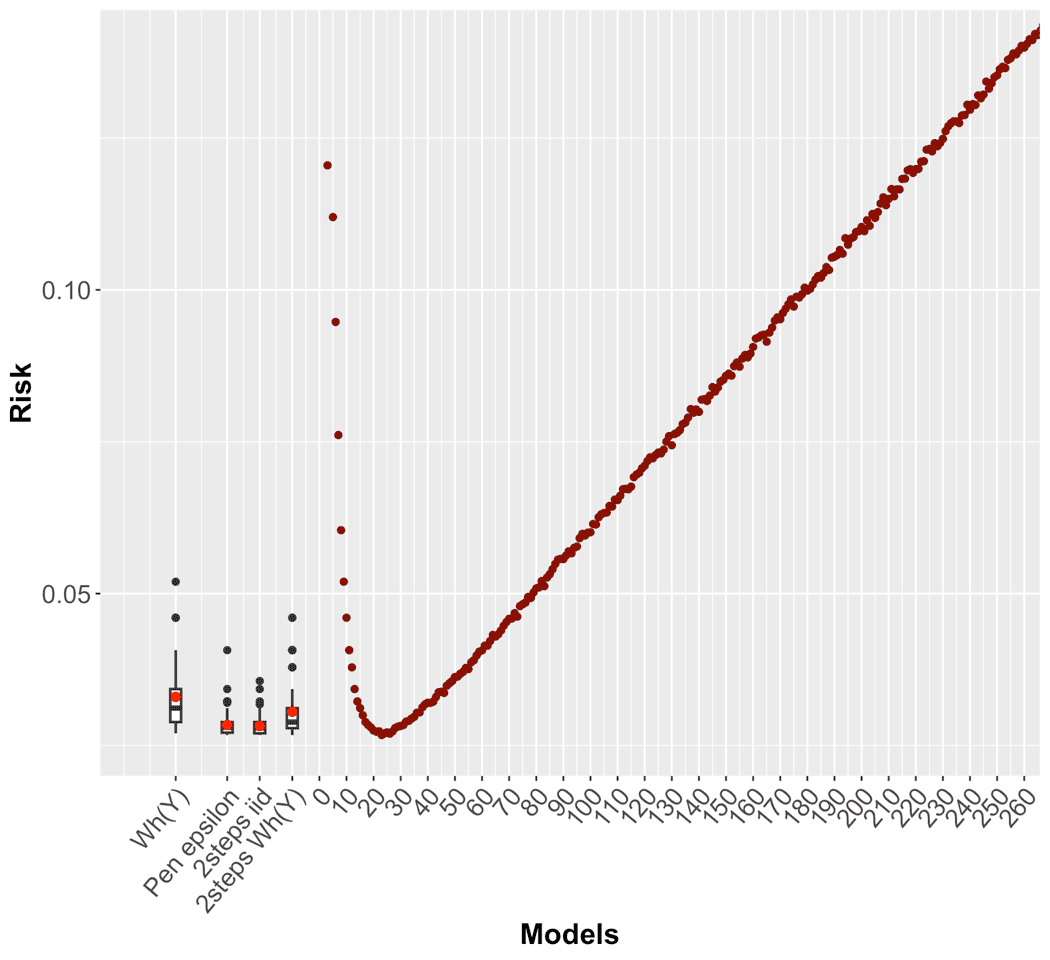}
\caption{$f_2$, $n=2000$}
\end{subfigure}
\caption{Risk performance of the penalization procedures for the two regression functions and two sample sizes in Experiment~11. The error process is a FGN with $H=0.7$, 
and the design  is obtained from  a FGN  with $H=0.7$ via the transformation \eqref{transFGN}. With the method $Wh(Y)$ the penalty is computed with a direct Whittle estimation of $H$ on $Y$ : it exhibits poorer performance compared the two steps approaches.}
\label{fig:long_range_dependence_error_FGN07_design_FGN07-YWhittleSimple}
\end{figure}

\medskip  
 
Overall, these experiments lead to two main conclusions.
First, the strong performance of \(\pen_{\epsilon}\) supports the validity of the proposed bounds on the variance.
This confirms that the theoretical guarantees derived for the variance are sufficiently sharp. Second, the \textit{two-step Whittle\((Y)\)} method proves sufficiently robust to enable mimic the behavior of \(\pen_{\epsilon}\), and then select a reasonable model, particularly in situations where the nature of the error processes and the design is unknown. This is particularly satisfying, as in data science applications, the underlying stochastic processes are rarely known in practice.
 
\begin{Rem}
A simpler strategy would consist in estimating the Hurst coefficient from $Y$ with the Whittle estimator, and next directly defining the penalty shape proportional to $ m^{2-2\hat H}$. However, this method proves less robust than the two-step approaches, at least for the set of experiments presented in this section. Figure~\ref{fig:long_range_dependence_error_FGN07_design_FGN07-YWhittleSimple} illustrates this observation for Experiment 11. 
 \end{Rem}

\bibliographystyle{apalike}
\bibliography{biblio}

\appendix

\section{Proofs}

\subsection{Proof of Corollary \ref{cor1_sd}}
We need to control the term 
$$
 \mathbb{E} \left[ \left \| f^{\ast} - \hat{f}_{m}  \right \|_{n({\mathbf X})}^{2} \Big | {\mathbf X}_n \right] \, .
$$
Recall that $\hat f_m= \text{Proj}_{S_m}(Y) $ (with a slight abuse of notation we use the notation $\text{Proj}_{S_m}$ for the Projection operator on $S_m$ and for its $n({\mathbf X})\times n({\mathbf X})$ matrix). Let $\varepsilon=(\varepsilon_1, \ldots, \varepsilon_n)^t$, and recall that we identify $f^*$ with the vector $(f^{\ast}(X_1), \ldots, f^{\ast}(X_n))^t$. Then, by orthogonality 
$$
\left \| f^{\ast} - \hat{f}_{m}  \right \|_{n({\mathbf X})}^{2}= \left \| f^{\ast} - \text{Proj}_{S_m}(f^*) \right \|_{n({\mathbf X})}^{2}+ \left \| \text{Proj}_{S_m}(\varepsilon) \right \|_{n({\mathbf X})}^{2}  \, , 
$$
and 
\begin{equation}
\label{ortdec}
\mathbb{E} \left[ \left \| f^{\ast} - \hat{f}_{m}  \right \|_{n({\mathbf X})}^{2} \Big | {\mathbf X}_n \right]= \left \| f^{\ast} - \text{Proj}_{S_m}(f^*) \right \|_{n({\mathbf X})}^{2}+ \mathbb{E} \left[\left \| \text{Proj}_{S_m}(\varepsilon) \right \|_{n({\mathbf X})}^{2} \big | {\mathbf X}_n \right]  \, .
\end{equation}
For the second term on right hand in \eqref{ortdec}, we write
\begin{equation}
\label{varterm_dec}
\mathbb{E} \left[\left \| \text{Proj}_{S_m}(\varepsilon) \right \|_{n({\mathbf X})}^{2} \big | {\mathbf X}_n \right]= \frac{1}{n({\mathbf X})\vee 1} \text{tr}\left ( \text{Proj}_{S_m} \Sigma_I({\mathbf X}_n)\right ) \leq (r+1)\frac{m}{n({\mathbf X})\vee 1} \rho_{{\mathbf e} |{\mathbf X}, I}  
\end{equation}
(see the upper bound (2.4) in \cite{caron2021gaussian}). We now consider the first term on right hand in \eqref{ortdec}. 
Since $\text{Proj}_{S_m}(f^*)(X_i)=\text{Proj}_{V_k}(f^*)(X_i)$ when $X_i \in Q_k$, we get
\begin{align}
\label{bias1}
\left \| f^{\ast} - \text{Proj}_{S_m}(f^*) \right \|_{n({\mathbf X})}^{2}&=\frac{1}{n({\mathbf X})\vee 1}\sum_{i=1}^n \sum_{k=1}^m \left ( f^{\ast}(X_i) - \text{Proj}_{S_m}(f^*)(X_i) \right )^2{\bf 1}_{X_i \in Q_k}\nonumber \\
&=\frac{1}{n({\mathbf X})\vee 1}\sum_{i=1}^n \sum_{k=1}^m \left ( f^{\ast}(X_i) - \text{Proj}_{V_k}(f^*)(X_i) \right )^2{\bf 1}_{X_i \in Q_k} \, .
\end{align}
Next, we use an approximation theorem of Jackson (see inequality (22) page 261 in the book by \cite{Timan1963}). Since $\text{Proj}_{V_k}(f^*)$ is the best approximation of 
$f^*{\bf 1}_{Q_k}$ of the form $a_0+ a_1x + \cdots +a_r x^r $ for the square norm $ g \mapsto \sum_{i=1}^n (g(X_i))^2 {\bf 1}_{X_i \in Q_k} $, it follows from this approximation theorem that, if $f^*{\bf 1}_I$ belongs to the class ${\mathcal H}_s(L)$ for some $s \in (0, r+1]$, then there exists a positive constant $M_r$  such that
\begin{equation}
\label{bias2}
\sum_{i=1}^n  \left ( f^{\ast}(X_i) - \text{Proj}_{V_k}(f^*)(X_i) \right )^2{\bf 1}_{X_i \in Q_k} \leq M_r\frac{L^2}{m^{2s}} \sum_{i=1}^n {\bf 1}_{X_i \in Q_k} \, .
\end{equation}
Combining  \eqref{bias1} and \eqref{bias2}, we get that, if $f^*$ belongs to the class ${\mathcal H}_s(L)$ for some $s \in (0, r+1]$, then 
\begin{equation}
\label{bias3}
\left \| f^{\ast} - \text{Proj}_{S_m}(f^*) \right \|_{n({\mathbf X})}^{2} \leq M_r\frac{L^2}{m^{2s}} \, .
\end{equation}
Combining the bounds \eqref{ortdec}, \eqref{varterm_dec}, \eqref{bias3} and Theorem \ref{genth}, we see that: if $f^*{\bf 1}_I$ belongs to the class ${\mathcal H}_s(L)$ for some $s \in (0, r+1]$, then there exist a positive constant $A_r$ such that, almost surely
\begin{equation}
\label{proof1:laststep}
\mathbb{E} \left[ \left \| f^{\ast} - \hat{f}_{\hat m} \right \|_{n({\mathbf X})}^{2} \Big |{\mathbf X}_n\right] \leq A_r \left( \inf_{m} \left\{ \frac{L^2}{m^{2s}} + \rho_{{\mathbf e} | {\mathbf X},I}\frac{m}{n({\mathbf X}) \vee 1} \right\} + \frac {\rho_{{\mathbf e} | {\mathbf X},I}}{n({\mathbf X})\vee 1} \right) \, .
\end{equation}
The proof of  Corollary \ref{cor1_sd} is complete by computing the minimum in $m$ in \eqref{proof1:laststep}.

\subsection{Proof of Corollary \ref{cor2_sd_alpha}}

Since $\rho_{{\mathbf e} | {\mathbf X}} \leq  \rho
$ almost surely, it follows from Corollary \ref{cor1_sd} that there exists a positive constant $K$ such that 
$$
\mathbb{E} \left[ \left \| f^{\ast} - \hat{f}_{\hat m} \right \|_{n({\mathbf X})}^{2} \right] \leq  K {\mathbb E}\left ( \frac {1}{(n({\mathbf X})\vee 1)^{\frac{2 s}{2s+1}}} \right ) \, .
$$
Now, let $c={\mathbb P}(X_1 \in I)/2$, and recall that $c>0$ by assumption.
Introducing the indicator functions of the  sets $\{n({\bf X})>nc\}$ and $\{n({\bf X})\leq nc\}$ in the right-hand expectation,  it follows that there exists a positive constant $A$ such that
\begin{equation}
\label{proof2:dec1}
\mathbb{E} \left[ \left \| f^{\ast} - \hat{f}_{\hat m} \right \|_{n({\mathbf X})}^{2} \right] \leq  \frac {A}{n^{\frac{2 s}{2s+1}}} +
K {\mathbb P}\left ( n({\mathbf X}) \leq n c \right )  
\, .
\end{equation}
Now, 
$$
{\mathbb P}\left ( n({\mathbf X}) \leq n c \right )= 
{\mathbb P}\left ( \frac{n({\mathbf X})}{n}- {\mathbb P}(X_1 \in I) \leq - c \right ) \leq {\mathbb P}\left (\left | \frac{n({\mathbf X})}{n}- {\mathbb P}(X_1 \in I) \right | \geq  c \right ) \, .
$$
Applying Markov Inequality at order 2 and using the stationarity of the sequence $(X_i)_{i \geq 1}$, we get
\begin{equation}\label{proof2:sumcov}
{\mathbb P}\left ( n({\mathbf X}) \leq n c \right )\leq \frac{2}{n c^2}\sum_{k=1}^n \left |\text{Cov}({\bf 1}_{X_1 \in I}, {\bf 1}_{X_{k} \in I}) \right | \, .
\end{equation}
Now, by \eqref{defalpha}
\begin{equation}
\label{proof2:covbound}
\left |\text{Cov}({\bf 1}_{X_1 \in I}, {\bf 1}_{X_{k} \in I}) \right |  \leq \alpha(k-1) \, .
\end{equation}
From \eqref{proof2:dec1}, \eqref{proof2:sumcov} and \eqref{proof2:covbound}, we get that there exists a positive constant $B$ such that 
\begin{equation}
\mathbb{E} \left[ \left \| f^{\ast} - \hat{f}_{\hat m} \right \|_{n({\mathbf X})}^{2} \right] \leq  \frac {A}{n^{\frac{2 s}{2s+1}}} +
\frac B n  \sum_{k=0}^{n-1} \alpha(k) \, .
\end{equation}
This completes the proof of Corollary \ref{cor2_sd_alpha}.

\subsection{Proof of Corollary \ref{cor_heteroiid}}
Recall that we assume here that $(X_i)_{1 \leq i \leq n}$ is i.i.d.
Let $Z_i=\sigma^2(X_i) {\mathbf 1}_{X_i \in I}$ and 
$Z_n^*= \max_{1 \leq i \leq n} Z_i$. Applying Corollary \ref{cor1_sd} to this particular situation, there exists a positive constant $K$ such that 
$$
\mathbb{E} \left[ \left \| f^{\ast} - \hat{f}_{\hat m} \right \|_{n({\mathbf X})}^{2} \right] \leq K {\mathbb E} \left (\frac {(Z^*_n)^{\frac{2 s}{2s+1}}}{(n({\mathbf X})\vee 1)^{\frac{2s}{2s+1}}} \right ) + K {\mathbb E} \left (
\frac {Z^*_n}{n({\mathbf X})\vee 1} \right ) \,.
$$
Keeping the same notations as in the proof of Corollary \ref{cor2_sd_alpha}, and introducing 
the indicator functions of the  sets $\{n({\bf X})>nc\}$ and $\{n({\bf X})\leq nc\}$ in the right-hand expectations,  it follows that there exists a positive constant $A$ such that
$$
\mathbb{E} \left[ \left \| f^{\ast} - \hat{f}_{\hat m} \right \|_{n({\mathbf X})}^{2} \right] \leq  A\frac {{\mathbb E}\big((Z^*_n)^{\frac{2 s}{2s+1} }\big )}{n^{\frac{2 s}{2s+1}}} + A\frac {{\mathbb E}(Z^*_n)}{n} +
A {\mathbb E}\left ( ((Z^*_n)^{\frac{2 s}{2s+1}}+Z^*_n){\bf 1}_{n({\mathbf X}) \leq n c} \right )  
\, .
$$
For the last term on right hand, since we always assume that $Z_i$ has a moment of order $p>1$, we write
$$
{\mathbb E}\left ( ((Z^*_n)^{\frac{2 s}{2s+1}}+Z^*_n){\bf 1}_{n({\mathbf X}) \leq n c} \right ) \leq \|(Z^*_n)^{\frac{2 s}{2s+1}}+Z^*_n \|_p \left ({\mathbb P}\left ( n({\mathbf X}) \leq n c \right ) \right )^{(p-1)/p} \, .
$$
Now, by Hoeffding Inequality applied to ${\mathbb P}\left ( n({\mathbf X}) \leq n c \right )$, we see that that  this last upper bound tends to zero at an exponential rate, and is therefore negligible. It remains to evaluate the two expectations ${\mathbb E}\big((Z^*_n)^{\frac{2 s}{2s+1} }\big )$ and ${\mathbb E}(Z^*_n)$. We apply here the general results of \cite{CorreaRomero2021}. If $Z_i$ has a moment of order $p>1$, then ${\mathbb E}(Z^*_n)=o(n^{1/p})$, from which we immediately get that
${\mathbb E}\big((Z^*_n)^{\frac{2 s}{2s+1} }\big )=o(n^{\frac{2 s}{(2s+1)p}})$. This proves Item 1 of Corollary \ref{cor_heteroiid}. From the same paper, we also know that if $Z_i$ has an exponential moment then 
${\mathbb E}(Z^*_n)=o(\log (n))$, from which we immediately get that ${\mathbb E}\big((Z^*_n)^{\frac{2 s}{2s+1} }\big )=o\big( (\log (n))^{\frac{2s}{2s+1}}\big )$. This proves Item 2 of Corollary \ref{cor_heteroiid}. Item 3 of Corollary \ref{cor_heteroiid} is immediate, since in that case $\|Z_n^*\|_\infty = \|Z_1\|_\infty < \infty$. 

\subsection{Proof of Theorem \ref{ThLD}}
\label{PrThLD}
According to the general form of the penalty given in \eqref{GenPen}, we have to bound up the term
$$
\text{tr}\left ( \text{Proj}_{S_m} \Sigma_I({\mathbf X}_n)\right )=(n({\mathbf X})\vee 1) \mathbb{E} \left[\left \| \text{Proj}_{S_m}(\varepsilon) \right \|_{n({\mathbf X})}^{2} \big | {\mathbf X}_n \right] \, .
$$
Note that the spaces  $V_{j}$ introduced in \eqref{Vj} are mutually orthogonal, in such a way that :
\begin{equation}
\label{ort1_ld}
\left  \|\text{Proj}_{S_{m}}(\varepsilon)\right \|_{n({\mathbf X})}^{2} = \sum_{j=1}^{m}  
\left\|\text{Proj}_{V_{j}}(\varepsilon) \right \|_{n({\mathbf X})}^{2} \, .
\end{equation}
Let $\langle \cdot, \cdot \rangle$ be the usual scalar product on ${\mathbb R}^n$, and let $(e_{(j, k)})_{ 1 \leq k \leq \kappa_{j}}$ be an orthonormal basis of $V_{j}$ with respect to $\langle \cdot, \cdot \rangle$. Note that $\kappa_{j} \leq r+1$, and that the coordinates  $e_{(j, k), p}$ of the vector  $e_{(j, k)}$ are 0 if  $p$ does not belong to  $K_{j}= \{ k \in \{1, \ldots, n \}: X_{k} \in Q_{j} \}$. Let $b_{j} : \{1, \ldots  ,\ell_{j}\} \mapsto K_{j}$ be the  increasing map  describing the set $K_{j}$.
We have 
\begin{equation}
\label{ort2_ld}
(n({\mathbf X})\vee 1)  {\mathbb E} \left [ \left \|\text{Proj}_{V_{j}}(\varepsilon) \right \|_{n({\mathbf X})}^{2} \big | {\mathbf X}_n \right ]= \sum_{k=1}^{\kappa_{j} }{\mathbb E}\left(\langle \varepsilon, e_{(j, k)}\rangle^2 \big | {\mathbf X}_n \right)\, .
\end{equation}
Now
\begin{multline*}
{\mathbb E}\left(\langle \varepsilon, e_{(j, k)}\rangle^2 \big | {\mathbf X}_n \right)= {\mathbb E}\left( \left(\sum_{p=1}^{\ell_{j}}
\varepsilon_{b_{j}(p)} e_{(j, k), b_{j}(p)}\Big | {\mathbf X}_n \right )^2\right)\\
=\sum_{p=1}^{\ell_{j}}\sum_{q=1}^{\ell_{j}}
e_{(j, k), b_{j}(p)}  e_{(j, k), b_{j}(q)}
\mathrm{Cov}(\varepsilon_{b_{j}(p)}, \varepsilon_{b_{j}(q)} | {\mathbf X}_n)\\
\leq c_I({\mathbf X}_n)+ 2 \sum_{p=1}^{\ell_{j}-1} \sum_{q=p+1}^{\ell_{j}} | e_{(j, k), b_{j}(p)} e_{(j, k), b_{j}(q)} | \frac{c_I({\mathbf X}_n)}{ (b_{j}(q)-b_{j}(p)+1)^{\gamma}}  \, .
\end{multline*}
Since $b_j$ is increasing, we have 
$
 (b_{j}(q) - b_{j}(p)+1)^{- \gamma} \leq  (q -p +1)^{- \gamma} 
$. It follows that 
\begin{multline*}
{\mathbb E}\left(\langle \varepsilon, e_{(j, k)}\rangle^2 \big | {\mathbf X}_n \right)\leq c_I({\mathbf X}_n)+ 2 c_I({\mathbf X}_n) \sum_{p=1}^{\ell_{j}-1} \sum_{k=1}^{\ell_{j}-p}  |e_{(j, k), b_{j}(p)} e_{(j, k), b_{j}(k+p)}| (k+1)^{-\gamma}\\
\leq c_I({\mathbf X}_n)+ 2 c_I({\mathbf X}_n)\sum_{k=1}^{\ell_{j}-1} (k+1)^{-\gamma} \sum_{p=1}^{\ell_{j}-k}  |e_{(j, k), b_{j}(p)} e_{(j, k), b_{j}(k+p)}| \, . 
\end{multline*}
Since $\langle e_{(j, k)}, e_{(j, k)} \rangle =1$, by Cauchy-Schwarz's inequality
$$
\sum_{p=1}^{\ell_{j}-k}  |e_{(j, k), b_{j}(p)} e_{(j, k), b_{j}(k+p)}| \leq 1 \, ,
$$
and consequently
\begin{equation}
\label{bound_fin}
{\mathbb E}\left(\langle \varepsilon, e_{(j, k)}\rangle^2 \big | {\mathbf X}_n \right)\leq c_I({\mathbf X}_n) \kappa_\gamma \ell_{j}^{1-\gamma} 
\end{equation}
for some positive constant $\kappa_\gamma$.

Combining \eqref{ort1_ld}, \eqref{ort2_ld} and \eqref{bound_fin} we see that there exists a positive constant  $\kappa_{\gamma, r}$ such that
$$
 \text{tr}\left ( \text{Proj}_{S_m} \Sigma_I({\mathbf X}_n)\right )=(n({\mathbf X})\vee 1) \mathbb{E} \left[\left \| \text{Proj}_{S_m}(\varepsilon) \right \|_{n({\mathbf X})}^{2} \big | {\mathbf X}_n \right] \leq \kappa_{\gamma, r} c_I({\mathbf X}_n) \sum_{j=1}^{m}   \ell_{j}^{1-\gamma} \, .
$$
Applying Jensen's inequality (since $x \mapsto x^{1-\gamma}$ is a concave function), 
$$
   \sum_{j=1}^{m} \ell_{j}^{1-\gamma}= m   \sum_{j=1}^{m}  \frac {\ell_{j}^{1-\gamma}} {m}  \leq m \left (  \sum_{j=1}^{m}  \frac  {\ell_{j}} {m} \right )^{1- \gamma} \, .
 $$
 Since $  \sum_{j=1}^{m}  \ell_{j}=n({\mathbf X})$, we get that
 $$
  \sum_{j=1}^{m} \ell_{j}^{1-\gamma} \leq m^\gamma n({\mathbf X})^{1- \gamma}, 
 $$
 and consequently 
\begin{equation}
\label{mainTerm}
\frac{ \text{tr}\left ( \text{Proj}_{S_m} \Sigma_I({\mathbf X}_n)\right )}{(n({\mathbf X})\vee 1)}
  \leq \kappa_{\gamma, r} c_I({\mathbf X}_n) \left (\frac{m}{n({\mathbf X})\vee 1}\right )^\gamma \, .
\end{equation}
The upper  bound \eqref{mainTerm} gives a control of the main term of the penalty \eqref{GenPen}. However, we still have to check that the terms involving the spectral radius $\rho_{{\mathbf e} | {\mathbf X}, I}$ are of a smaller order with respect to this main term. 

Let $a : \{1, \ldots  , n({\mathbf X}) \} \mapsto A_I$ be the  increasing map  describing the set $A_I$ (defined in \eqref{defAI}). Let also
$$
B_1=\left \{ x \in {\mathbb R}^{n({\mathbf X})} : \sum_{i=1}^{n({\mathbf X})} x_i^2 =1\right \} \, .
$$
By definition of the spectral radius 
$$
\rho_{{\mathbf e} | {\mathbf X}, I}= \sup_{x \in B_1} x^t \Sigma_I({\mathbf X}_n) x= \sup_{x \in B_1} \text{Var}\left ( \sum_{i=1}^{n({\mathbf X})}x_i \varepsilon_{a(i)}\right ) = 
\sup_{x \in B_1}
\sum_{i=1}^{n({\mathbf X})} \sum_{j=1}^{n({\mathbf X})} x_ix_j \text{Cov}(\varepsilon_{a(i)}, \varepsilon_{a(j)} ) \, .
$$
Then, with the same arguments as to get \eqref{bound_fin},  there exists a positive constant $\kappa_\gamma$ such that
\begin{equation}
\label{Boundrho}
\rho_{{\mathbf e} | {\mathbf X}, I} \leq  \kappa_\gamma c_I({\mathbf X}_n) n({\mathbf X})^{1-\gamma}  \, .
\end{equation}
This together with \eqref{mainTerm} imply that a penalty of the form \eqref{pen2_ld} will indeed satisfy the inequality \eqref{GenPen}. Applying Theorem 2.1 of~\cite{caron2021gaussian}, Theorem \ref{ThLD} follows.

\subsection{Proof of Corollary \ref{cor1_ld}}
Following the proof of Corollary \ref{cor1_sd}, we obtain that: if $f^*{\bf 1}_I$ belongs to the class ${\mathcal H}_s(L)$ for some $s \in (0, r+1]$, then there exist a positive constant $A_{r, \gamma}$ such that, almost surely
\begin{equation}
\label{laststep2}
\mathbb{E} \left[ \left \| f^{\ast} - \hat{f}_{\hat m} \right \|_{n({\mathbf X})}^{2} \Big |{\mathbf X}_n\right] \leq A_{r, \gamma} \left( \inf_{m} \left\{ \frac{L^2}{m^{2s}} + c_I({\mathbf X}_n) \left (\frac{m}{n({\mathbf X})\vee 1}\right )^\gamma \right\} + \frac {\rho_{{\mathbf e} | {\mathbf X},I}}{n({\mathbf X})\vee 1} \right) \, .
\end{equation}
The proof of  Corollary \ref{cor1_ld} is complete by computing the minimum in $m$ in \eqref{laststep2}.

\subsection{Proof of Corollary \ref{cor2_ld_alpha}}
Since $c_I({\mathbf X}_n) \leq  c
$ almost surely and since \eqref{Boundrho} is satisfied, it follows from Corollary \ref{cor1_ld} that there exists a positive constant $K$ such that 
$$
\mathbb{E} \left[ \left \| f^{\ast} - \hat{f}_{\hat m} \right \|_{n({\mathbf X})}^{2} \right] \leq  K {\mathbb E}\left ( \frac {1}{(n({\mathbf X})\vee 1)^{\frac{2\gamma s}{2s+\gamma}}} \right ) \, .
$$
The end of the proof is similar to that of Corollary \ref{cor2_sd_alpha}.

\subsection{Proof of Proposition \ref{prop_ld_improved}}

Recall that, since $r=0$, we have
$$
\|\text{Proj}_{S_{m}}(\varepsilon)\|_n^2=  \frac 1 n\sum_{i=1}^{m}   \ell_{i}(\bar \varepsilon_{i})^2 \, ,
$$
where
$$
   \bar \varepsilon_{i} = \frac{1}{\ell_{i}}  \sum_{k=1}^{n}  \varepsilon_{k} 
   {\bf 1}_{X_{k} \in Q_{i}}\, , \quad  \ell_{i}=  \sum_{k=1}^{n}  
   {\bf 1}_{X_{k} \in Q_{i}}
$$
Let $p_i={\mathbb P}(X_1 \in Q_i)$ (without loss of generality, one can assume that $p_i>0$, otherwise there is no observations in $Q_i$). We start from the inequality
\begin{equation}
\label{dec1_improved}
\|\text{Proj}_{S_{m}}(\varepsilon)\|_n^2 \leq  \sum_{i=1}^{m}   \frac{\ell_{i}}{n} {\bf 1}_{\ell_{i} \leq n p_i/2}(\bar \varepsilon_{i})^2 + \sum_{i=1}^{m}   \frac{\ell_{i}}{n} {\bf 1}_{\ell_{i} > n p_i/2}(\bar \varepsilon_{i})^2.
\end{equation}
For the second term on right hand in  \eqref{dec1_improved}, we get 
$$
\sum_{i=1}^{m}   \frac{\ell_{i}}{n} {\bf 1}_{\ell_{i} > n p_i/2}(\bar \varepsilon_{i})^2 \leq \sum_{i=1}^{m} \frac{1}{\ell_i} {\bf 1}_{\ell_{i} > n p_i/2}\left( \frac{1}{\sqrt n}\sum_{k=1}^{n}  \varepsilon_{k} 
   {\bf 1}_{X_{k} \in Q_{i}}\right)^2 \leq 
\frac{2}{n} \sum_{i=1}^{m} \frac{1}{p_i} \left ( \frac {1} {\sqrt n}  \sum_{k=1}^{n}  \varepsilon_{k} 
   {\bf 1}_{X_{k} \in Q_{i}}\right )^2 .
$$
Consequently
$$
{\mathbb E} \left (  \sum_{i=1}^{m}   \frac{\ell_{i}}{n} {\bf 1}_{\ell_{i} > n p_i/2}(\bar \varepsilon_{i})^2 \right )
\leq \frac{2}{n} \sum_{i=1}^{m} \frac{1}{p_i} {\mathbb E}\left (\left ( \frac {1} {\sqrt n}  \sum_{k=1}^{n}  \varepsilon_{k} 
   {\bf 1}_{X_{k} \in Q_{i}}\right )^2 \right).
$$
The sequence  $(\varepsilon_{k} 
   {\bf 1}_{X_{k} \in Q_{i}})_{k \geq 1}$ is strictly stationary with auto-covariance coefficients $(\kappa_i(k))_{k \geq 0}$. Let us give some bounds on these coefficients. Since the sequences ${\mathbf X}_n$ and ${\mathbf e}_n$ are independent, $\kappa_i(0)={\mathbb E}(\varepsilon_1^2 {\bf 1}_{X_{1} \in Q_{i}})= \gamma_\varepsilon (0)p_i$ and 
 \begin{align}
 |\kappa_i(k)|=|{\mathbb E}(\varepsilon_1 \varepsilon_{k+1} {\bf 1}_{X_{1} \in Q_{i}}{\bf 1}_{X_{k+1} \in Q_{i}})|  &\leq 
  |{\mathbb E}(\varepsilon_1 \varepsilon_{k+1})|p_i^2 +  |{\mathbb E}(\varepsilon_1 \varepsilon_{k+1} ({\bf 1}_{X_{1} \in Q_{i}}{\bf 1}_{X_{k+1} \in Q_{i}}-p_i^2))|
  \nonumber \\ 
  &\leq  |\gamma_\varepsilon(k)|p_i^2+|\gamma_\varepsilon(k)\rho(k)| p_i \, . 
 \end{align}
 Let $A=\gamma_\varepsilon(0) + 2 \sum_{k \geq 1}|\gamma_\varepsilon(k)\rho(k)| $. Since $\gamma_\varepsilon$ satisfies \eqref{gammak}, we infer that, for some positive constant $K_{\gamma, c}$, 
 $$
 {\mathbb E}\left (\left ( \frac {1} {\sqrt n}  \sum_{k=1}^{n}  \varepsilon_{k} 
   {\bf 1}_{X_{k} \in Q_{i}}\right )^2 \right) \leq \kappa_{i}(0) + 2 \sum_{k=1}^n |\kappa_{i}(k)| \leq 
   A p_i + K_{\gamma, c} n^{1-\gamma} p_i^2 \, , 
 $$
 in such a way that 
 \begin{equation}
 \label{pen_fin_ld_improved}
 {\mathbb E} \left (  \sum_{i=1}^{m}   \frac{\ell_{i}}{n} {\bf 1}_{\ell_{i} > n p_i/2}(\bar \varepsilon_{i})^2 \right )
\leq 2A\frac{m}{n} + \frac{2K_{\gamma, c}}{n^{\gamma}} \, .
 \end{equation}

 Let us now deal with the first term on right hand in \eqref{dec1_improved}.  We first note that
 $$
 {\mathbb E}\left ((\bar \varepsilon_{i})^2 \right |{\mathbf X}_n) \leq \gamma_\varepsilon(0) \quad \text{almost surely,}
 $$
 and consequently
 $$
 {\mathbb E}\left ( \sum_{i=1}^{m}   \frac{\ell_{i}}{n} {\bf 1}_{\ell_{i} \leq n p_i/2}(\bar \varepsilon_{i})^2 \right )
 \leq \frac{\gamma_\varepsilon (0)}{2} \sum_{i=1}^{m} p_i {\mathbb P}(\ell_{i} \leq n p_i/2) \leq \frac{\gamma_\varepsilon (0)}{2} \sum_{i=1}^{m} p_i {\mathbb P}\left(\left |\frac{\ell_{i}}{n} -p_i \right | \geq  \frac{p_i}2 \right ).
 $$
 Applying Markov Inequality at order 2, we get that 
 \begin{equation}
 \label{last_ineq_improved}
 {\mathbb E}\left ( \sum_{i=1}^{m}   \frac{\ell_{i}}{n} {\bf 1}_{\ell_{i} \leq n p_i/2}(\bar \varepsilon_{i})^2 \right )
 \leq \frac{\gamma_\varepsilon(0)}{2} \sum_{i=1}^{m}\frac{4}{p_i n} p_i(1-p_i) \leq 2\gamma_\varepsilon (0)\frac{m}{n}\, .
 \end{equation}
 Let $B=2A+2\gamma_\varepsilon(0) $. Combining \eqref{dec1_improved}, \eqref{pen_fin_ld_improved} and \eqref{last_ineq_improved}, we obtain that
 \begin{equation}
 {\mathbb E} \left (  \left \|\text{Proj}_{S_{m}}(\varepsilon) \right \|_n^2 \right ) \leq B\frac{m}{n} + \frac{2K_{\gamma, c}}{n^{\gamma}} \, ,
 \end{equation}
 and the proof of Proposition \ref{prop_ld_improved} is complete.
 



\end{document}